\documentclass[12pt,thmsa]{article}
\usepackage{amssymb}

\usepackage{epsfig}
\usepackage{sw20jart}

\input{tcilatex}

\makeatletter

\def\diagram{\leftwidth=\z@ \rightwidth=\z@ \topheight=\z@
\botheight=\z@ \setbox\@picbox\hbox\bgroup}

\def\enddiagram{\egroup\wd\@picbox\rightwidth\unitlength
\ht\@picbox\topheight\unitlength \dp\@picbox\botheight\unitlength
\hskip\leftwidth\unitlength\box\@picbox}

\def\bfig{\begin{diagram}}
\def\efig{\end{diagram}}
\newcount\wideness \newcount\leftwidth \newcount\rightwidth
\newcount\highness \newcount\topheight \newcount\botheight

\def\ratchet#1#2{\ifnum#1<#2 \global #1=#2 \fi}

\def\putbox(#1,#2)#3{%
\horsize{\wideness}{#3} \divide\wideness by 2
{\advance\wideness by #1 \ratchet{\rightwidth}{\wideness}}
{\advance\wideness by -#1 \ratchet{\leftwidth}{\wideness}}
\vertsize{\highness}{#3} \divide\highness by 2
{\advance\highness by #2 \ratchet{\topheight}{\highness}}
{\advance\highness by -#2 \ratchet{\botheight}{\highness}}
\put(#1,#2){\makebox(0,0){$#3$}}}

\def\putlbox(#1,#2)#3{%
\horsize{\wideness}{#3}
{\advance\wideness by #1 \ratchet{\rightwidth}{\wideness}}
{\ratchet{\leftwidth}{-#1}}
\vertsize{\highness}{#3} \divide\highness by 2
{\advance\highness by #2 \ratchet{\topheight}{\highness}}
{\advance\highness by -#2 \ratchet{\botheight}{\highness}}
\put(#1,#2){\makebox(0,0)[l]{$#3$}}}

\def\putrbox(#1,#2)#3{%
\horsize{\wideness}{#3}
{\ratchet{\rightwidth}{#1}}
{\advance\wideness by -#1 \ratchet{\leftwidth}{\wideness}}
\vertsize{\highness}{#3} \divide\highness by 2
{\advance\highness by #2 \ratchet{\topheight}{\highness}}
{\advance\highness by -#2 \ratchet{\botheight}{\highness}}
\put(#1,#2){\makebox(0,0)[r]{$#3$}}}

\def\adjust[#1]{} 

\newcount \coefa
\newcount \coefb
\newcount \coefc
\newcount\tempcounta
\newcount\tempcountb
\newcount\tempcountc
\newcount\tempcountd
\newcount\xext
\newcount\yext
\newcount\xoff
\newcount\yoff
\newcount\gap%
\newcount\arrowtypea
\newcount\arrowtypeb
\newcount\arrowtypec
\newcount\arrowtyped
\newcount\arrowtypee
\newcount\height
\newcount\width
\newcount\xpos
\newcount\ypos
\newcount\run
\newcount\rise
\newcount\arrowlength
\newcount\halflength
\newcount\arrowtype
\newdimen\tempdimen
\newdimen\xlen
\newdimen\ylen
\newsavebox{\tempboxa}%
\newsavebox{\tempboxb}%
\newsavebox{\tempboxc}%

\newdimen\w@dth

\def\setw@dth#1#2{\setbox\z@\hbox{$#1$}\w@dth=\wd\z@
\setbox\@ne\hbox{$#2$}\ifnum\w@dth<\wd\@ne \w@dth=\wd\@ne \fi
\advance\w@dth by 1.2em}

\def\t@^#1_#2{\def\n@one{#1}\def\n@two{#2}\mathrel{\setw@dth{#1}{#2}
\mathop{\hbox to \w@dth{\rightarrowfill}}\limits
\ifx\n@one\empty\else ^{\box\z@}\fi
\ifx\n@two\empty\else _{\box\@ne}\fi}}
\def\t@@^#1{\@ifnextchar_ {\t@^{#1}}{\t@^{#1}_{}}}
\def\to{\@ifnextchar^ {\t@@}{\t@@^{}}}

\def\t@left^#1_#2{\def\n@one{#1}\def\n@two{#2}\mathrel{\setw@dth{#1}{#2}
\mathop{\hbox to \w@dth{\leftarrowfill}}\limits
\ifx\n@one\empty\else ^{\box\z@}\fi
\ifx\n@two\empty\else _{\box\@ne}\fi}}
\def\t@@left^#1{\@ifnextchar_ {\t@left^{#1}}{\t@left^{#1}_{}}}
\def\toleft{\@ifnextchar^ {\t@@left}{\t@@left^{}}}

\def\two@^#1_#2{\def\n@one{#1}\def\n@two{#2}\mathrel{\setw@dth{#1}{#2}
\mathop{\vcenter{\hbox to \w@dth{\rightarrowfill}\kern-1.7ex
                 \hbox to \w@dth{\rightarrowfill}}%
       }\limits
\ifx\n@one\empty\else ^{\box\z@}\fi
\ifx\n@two\empty\else _{\box\@ne}\fi}}
\def\tw@@^#1{\@ifnextchar_ {\two@^{#1}}{\two@^{#1}_{}}}
\def\two{\@ifnextchar^ {\tw@@}{\tw@@^{}}}

\def\tofr@^#1_#2{\def\n@one{#1}\def\n@two{#2}\mathrel{\setw@dth{#1}{#2}
\mathop{\vcenter{\hbox to \w@dth{\rightarrowfill}\kern-1.7ex
                 \hbox to \w@dth{\leftarrowfill}}%
       }\limits
\ifx\n@one\empty\else ^{\box\z@}\fi
\ifx\n@two\empty\else _{\box\@ne}\fi}}
\def\t@fr@^#1{\@ifnextchar_ {\tofr@^{#1}}{\tofr@^{#1}_{}}}
\def\tofro{\@ifnextchar^ {\t@fr@}{\t@fr@^{}}}

\def\mon{\mathop{\m@th\hbox to
      14.6\P@{\lasyb\char'51\hskip-2.1\P@$\arrext$\hss
$\mathord\rightarrow$}}\limits} 
\def\leftmono{\mathrel{\m@th\hbox to
14.6\P@{$\mathord\leftarrow$\hss$\arrext$\hskip-2.1\P@\lasyb\char'50%
}}\limits} 
\mathchardef\arrext="0200       

\def\settypes(#1,#2,#3){\arrowtypea#1 \arrowtypeb#2 \arrowtypec#3}
\def\settoheight#1#2{\setbox\@tempboxa\hbox{#2}#1\ht\@tempboxa\relax}%
\def\settodepth#1#2{\setbox\@tempboxa\hbox{#2}#1\dp\@tempboxa\relax}%
\def\settokens[#1`#2`#3`#4]{%
     \def\tokena{#1}\def\tokenb{#2}\def\tokenc{#3}\def\tokend{#4}}
\def\setsqparms[#1`#2`#3`#4;#5`#6]{%
\arrowtypea #1
\arrowtypeb #2
\arrowtypec #3
\arrowtyped #4
\width #5
\height #6
}
\def\setpos(#1,#2){\xpos=#1 \ypos#2}

\def\settriparms[#1`#2`#3;#4]{\settripairparms[#1`#2`#3`1`1;#4]}%

\def\settripairparms[#1`#2`#3`#4`#5;#6]{%
\arrowtypea #1
\arrowtypeb #2
\arrowtypec #3
\arrowtyped #4
\arrowtypee #5
\width #6
\height #6
}

\def\resetparms{\settripairparms[1`1`1`1`1;500]\width 500}

\resetparms

\def\mvector(#1,#2)#3{
\put(0,0){\vector(#1,#2){#3}}%
\put(0,0){\vector(#1,#2){26}}%
}
\def\evector(#1,#2)#3{{
\arrowlength #3
\put(0,0){\vector(#1,#2){\arrowlength}}%
\advance \arrowlength by-30
\put(0,0){\vector(#1,#2){\arrowlength}}%
}}

\def\horsize#1#2{%
\settowidth{\tempdimen}{$#2$}%
#1=\tempdimen
\divide #1 by\unitlength
}

\def\vertsize#1#2{%
\settoheight{\tempdimen}{$#2$}%
#1=\tempdimen
\settodepth{\tempdimen}{$#2$}%
\advance #1 by\tempdimen
\divide #1 by\unitlength
}

\def\putvector(#1,#2)(#3,#4)#5#6{{%
\ifnum3<\arrowtype
\putdashvector(#1,#2)(#3,#4)#5\arrowtype
\else
\ifnum\arrowtype<-3
\putdashvector(#1,#2)(#3,#4)#5\arrowtype
\else
\xpos=#1
\ypos=#2
\run=#3
\rise=#4
\arrowlength=#5
\ifnum \arrowtype<0
    \ifnum \run=0
        \advance \ypos by-\arrowlength
    \else
        \tempcounta \arrowlength
        \multiply \tempcounta by\rise
        \divide \tempcounta by\run
        \ifnum\run>0
            \advance \xpos by\arrowlength
            \advance \ypos by\tempcounta
        \else
            \advance \xpos by-\arrowlength
            \advance \ypos by-\tempcounta
        \fi
    \fi
    \multiply \arrowtype by-1
    \multiply \rise by-1
    \multiply \run by-1
\fi
\ifcase \arrowtype
\or \put(\xpos,\ypos){\vector(\run,\rise){\arrowlength}}%
\or \put(\xpos,\ypos){\mvector(\run,\rise)\arrowlength}%
\or \put(\xpos,\ypos){\evector(\run,\rise){\arrowlength}}%
\fi\fi\fi
}}

\def\putsplitvector(#1,#2)#3#4{
\xpos #1
\ypos #2
\arrowtype #4
\halflength #3
\arrowlength #3
\gap 140
\advance \halflength by-\gap
\divide \halflength by2
\ifnum\arrowtype>0
   \ifcase \arrowtype
   \or \put(\xpos,\ypos){\line(0,-1){\halflength}}%
       \advance\ypos by-\halflength
       \advance\ypos by-\gap
       \put(\xpos,\ypos){\vector(0,-1){\halflength}}%
   \or \put(\xpos,\ypos){\line(0,-1)\halflength}%
       \put(\xpos,\ypos){\vector(0,-1)3}%
       \advance\ypos by-\halflength
       \advance\ypos by-\gap
       \put(\xpos,\ypos){\vector(0,-1){\halflength}}%
   \or \put(\xpos,\ypos){\line(0,-1)\halflength}%
       \advance\ypos by-\halflength
       \advance\ypos by-\gap
       \put(\xpos,\ypos){\evector(0,-1){\halflength}}%
   \fi
\else \arrowtype=-\arrowtype
   \ifcase\arrowtype
   \or \advance \ypos by-\arrowlength
       \put(\xpos,\ypos){\line(0,1){\halflength}}%
       \advance\ypos by\halflength
       \advance\ypos by\gap
       \put(\xpos,\ypos){\vector(0,1){\halflength}}%
   \or \advance \ypos by-\arrowlength
       \put(\xpos,\ypos){\line(0,1)\halflength}%
       \put(\xpos,\ypos){\vector(0,1)3}%
       \advance\ypos by\halflength
       \advance\ypos by\gap
       \put(\xpos,\ypos){\vector(0,1){\halflength}}%
   \or \advance \ypos by-\arrowlength
       \put(\xpos,\ypos){\line(0,1)\halflength}%
       \advance\ypos by\halflength
       \advance\ypos by\gap
       \put(\xpos,\ypos){\evector(0,1){\halflength}}%
   \fi
\fi
}

\def\putmorphism(#1)(#2,#3)[#4`#5`#6]#7#8#9{{%
\run #2
\rise #3
\ifnum\rise=0
  \puthmorphism(#1)[#4`#5`#6]{#7}{#8}#9%
\else\ifnum\run=0
  \putvmorphism(#1)[#4`#5`#6]{#7}{#8}#9%
\else
\setpos(#1)%
\arrowlength #7
\arrowtype #8
\ifnum\run=0
\else\ifnum\rise=0
\else
\ifnum\run>0
    \coefa=1
\else
   \coefa=-1
\fi
\ifnum\arrowtype>0
   \coefb=0
   \coefc=-1
\else
   \coefb=\coefa
   \coefc=1
   \arrowtype=-\arrowtype
\fi
\width=2
\multiply \width by\run
\divide \width by\rise
\ifnum \width<0  \width=-\width\fi
\advance\width by60
\if l#9 \width=-\width\fi
\putbox(\xpos,\ypos){#4}
{\multiply \coefa by\arrowlength
\advance\xpos by\coefa
\multiply \coefa by\rise
\divide \coefa by\run
\advance \ypos by\coefa
\putbox(\xpos,\ypos){#5} }%
{\multiply \coefa by\arrowlength
\divide \coefa by2
\advance \xpos by\coefa
\advance \xpos by\width
\multiply \coefa by\rise
\divide \coefa by\run
\advance \ypos by\coefa
\if l#9%
   \putrbox(\xpos,\ypos){#6}%
\else\if r#9%
   \putlbox(\xpos,\ypos){#6}%
\fi\fi }%
{\multiply \rise by-\coefc
\multiply \run by-\coefc
\multiply \coefb by\arrowlength
\advance \xpos by\coefb
\multiply \coefb by\rise
\divide \coefb by\run
\advance \ypos by\coefb
\multiply \coefc by70
\advance \ypos by\coefc
\multiply \coefc by\run
\divide \coefc by\rise
\advance \xpos by\coefc
\multiply \coefa by140
\multiply \coefa by\run
\divide \coefa by\rise
\advance \arrowlength by\coefa
\ifcase\arrowtype
\or \put(\xpos,\ypos){\vector(\run,\rise){\arrowlength}}%
\or \put(\xpos,\ypos){\mvector(\run,\rise){\arrowlength}}%
\or \put(\xpos,\ypos){\evector(\run,\rise){\arrowlength}}%
\fi}\fi\fi\fi\fi}}

\newcount\numbdashes \newcount\lengthdash \newcount\increment

\def\howmanydashes{
\numbdashes=\arrowlength \lengthdash=40
\divide\numbdashes by \lengthdash
\lengthdash=\arrowlength
\divide\lengthdash by \numbdashes
\increment=\lengthdash
\multiply\lengthdash by 3
\divide\lengthdash by 5
}

\def\putdashvector(#1)(#2,#3)#4#5{%
\ifnum#3=0 \putdashhvector(#1){#4}#5
\else
\ifnum#2=0
\putdashvvector(#1){#4}#5\fi\fi}

\def\putdashhvector(#1,#2)#3#4{{%
\arrowlength=#3 \howmanydashes
\multiput(#1,#2)(\increment,0){\numbdashes}%
{\vrule height .4pt width \lengthdash\unitlength}
\arrowtype=#4 \xpos=#1
\ifnum\arrowtype<0 \advance\arrowtype by 7 \fi
\ifcase\arrowtype
\or \advance\xpos by 10
    \put(\xpos,#2){\vector(-1,0){\lengthdash}}
    \advance\xpos by 40
    \put(\xpos,#2){\vector(-1,0){\lengthdash}}
\or \advance \xpos by 10
    \put(\xpos,#2){\vector(-1,0){\lengthdash}}
    \advance\xpos by  \arrowlength
    \advance\xpos by  -50
    \put(\xpos,#2){\vector(-1,0){\lengthdash}}
\or \advance\xpos by 10
    \put(\xpos,#2){\vector(-1,0){\lengthdash}}
\or \advance\xpos by \arrowlength
    \advance\xpos by -\lengthdash
    \put(\xpos,#2){\vector(1,0){\lengthdash}}
\or {\advance\xpos by 10
    \put(\xpos,#2){\vector(1,0){\lengthdash}}}
    \advance\xpos by \arrowlength
    \advance\xpos by -\lengthdash
    \put(\xpos,#2){\vector(1,0){\lengthdash}}
\or \advance\xpos by \arrowlength
    \advance\xpos by -\lengthdash
    \put(\xpos,#2){\vector(1,0){\lengthdash}}
    \advance\xpos by -40
    \put(\xpos,#2){\vector(1,0){\lengthdash}}
   \fi
}}

\def\putdashvvector(#1,#2)#3#4{{%
\arrowlength=#3 \howmanydashes
\ypos=#2 \advance\ypos by -\arrowlength
\multiput(#1,#2)(0,\increment){\numbdashes}%
    {\vrule width .4pt height \lengthdash\unitlength}
\arrowtype=#4 \ypos=#2
\ifnum\arrowtype<0 \advance\arrowtype by 7 \fi
\ifcase\arrowtype
\or \advance\ypos by \arrowlength \advance\ypos by -40
    \put(#1,\ypos){\vector(0,1){\lengthdash}}
    \advance\ypos by -40
    \put(#1,\ypos){\vector(0,1){\lengthdash}}
\or \advance\ypos by 10
    \put(#1,\ypos){\vector(0,1){\lengthdash}}
    \advance\ypos by \arrowlength \advance\ypos by -40
    \put(#1,\ypos){\vector(0,1){\lengthdash}}
\or \advance\ypos by \arrowlength \advance\ypos by -40
    \put(#1,\ypos){\vector(0,1){\lengthdash}}
\or \advance\ypos by 10
    \put(#1,\ypos){\vector(0,-1){\lengthdash}}
\or \advance\ypos by 10
    \put(#1,\ypos){\vector(0,-1){\lengthdash}}
    \advance\ypos by \arrowlength \advance\ypos by -40
    \put(#1,\ypos){\vector(0,-1){\lengthdash}}
\or \advance\ypos by 10
    \put(#1,\ypos){\vector(0,-1){\lengthdash}}
    \advance\ypos by 40
    \put(#1,\ypos){\vector(0,-1){\lengthdash}}
\fi
}}

\def\puthmorphism(#1,#2)[#3`#4`#5]#6#7#8{{%
\xpos #1
\ypos #2
\width #6
\arrowlength #6
\arrowtype=#7
\putbox(\xpos,\ypos){#3\vphantom{#4}}%
{\advance \xpos by\arrowlength
\putbox(\xpos,\ypos){\vphantom{#3}#4}}%
\horsize{\tempcounta}{#3}%
\horsize{\tempcountb}{#4}%
\divide \tempcounta by2
\divide \tempcountb by2
\advance \tempcounta by30
\advance \tempcountb by30
\advance \xpos by\tempcounta
\advance \arrowlength by-\tempcounta
\advance \arrowlength by-\tempcountb
\putvector(\xpos,\ypos)(1,0)\arrowlength\arrowtype
\divide \arrowlength by2
\advance \xpos by\arrowlength
\vertsize{\tempcounta}{#5}%
\divide\tempcounta by2
\advance \tempcounta by20
\if a#8 %
   \advance \ypos by\tempcounta
   \putbox(\xpos,\ypos){#5}%
\else
   \advance \ypos by-\tempcounta
   \putbox(\xpos,\ypos){#5}%
\fi}}

\def\putvmorphism(#1,#2)[#3`#4`#5]#6#7#8{{%
\xpos #1
\ypos #2
\arrowlength #6
\arrowtype #7
\settowidth{\xlen}{$#5$}%
\putbox(\xpos,\ypos){#3}%
{\advance \ypos by-\arrowlength
\putbox(\xpos,\ypos){#4}}%
{\advance\arrowlength by-140
\advance \ypos by-70
\ifdim\xlen>0pt
   \if m#8%
      \putsplitvector(\xpos,\ypos)\arrowlength\arrowtype
   \else
   \putvector(\xpos,\ypos)(0,-1)\arrowlength\arrowtype
   \fi
\else
   \putvector(\xpos,\ypos)(0,-1)\arrowlength\arrowtype
\fi}%
\ifdim\xlen>0pt
   \divide \arrowlength by2
   \advance\ypos by-\arrowlength
   \if l#8%
      \advance \xpos by-40
      \putrbox(\xpos,\ypos){#5}%
   \else\if r#8%
      \advance \xpos by40
      \putlbox(\xpos,\ypos){#5}%
   \else
      \putbox(\xpos,\ypos){#5}%
   \fi\fi
\fi
}}

\def\putsquarep<#1>(#2)[#3;#4`#5`#6`#7]{{%
\setsqparms[#1]%
\setpos(#2)%
\settokens[#3]%
\puthmorphism(\xpos,\ypos)[\tokenc`\tokend`{#7}]{\width}{\arrowtyped}b%
\advance\ypos by \height
\puthmorphism(\xpos,\ypos)[\tokena`\tokenb`{#4}]{\width}{\arrowtypea}a%
\putvmorphism(\xpos,\ypos)[``{#5}]{\height}{\arrowtypeb}l%
\advance\xpos by \width
\putvmorphism(\xpos,\ypos)[``{#6}]{\height}{\arrowtypec}r%
}}

\def\putsquare{\@ifnextchar <{\putsquarep}{\putsquarep%
   <\arrowtypea`\arrowtypeb`\arrowtypec`\arrowtyped;\width`\height>}}
\def\square{\@ifnextchar< {\squarep}{\squarep
   <\arrowtypea`\arrowtypeb`\arrowtypec`\arrowtyped;\width`\height>}}
\def\squarep<#1>[#2`#3`#4`#5;#6`#7`#8`#9]{{
\setsqparms[#1]
\diagram
\putsquarep<\arrowtypea`\arrowtypeb`\arrowtypec`
\arrowtyped;\width`\height>
(0,0)[#2`#3`#4`{#5};#6`#7`#8`{#9}]
\enddiagram
}}                                                 
\def\putptrianglep<#1>(#2,#3)[#4`#5`#6;#7`#8`#9]{{%
\settriparms[#1]%
\xpos=#2 \ypos=#3
\advance\ypos by \height
\puthmorphism(\xpos,\ypos)[#4`#5`{#7}]{\height}{\arrowtypea}a%
\putvmorphism(\xpos,\ypos)[`#6`{#8}]{\height}{\arrowtypeb}l%
\advance\xpos by\height
\putmorphism(\xpos,\ypos)(-1,-1)[``{#9}]{\height}{\arrowtypec}r%
}}

\def\putptriangle{\@ifnextchar <{\putptrianglep}{\putptrianglep
   <\arrowtypea`\arrowtypeb`\arrowtypec;\height>}}
\def\ptriangle{\@ifnextchar <{\ptrianglep}{\ptrianglep
   <\arrowtypea`\arrowtypeb`\arrowtypec;\height>}}
\def\ptrianglep<#1>[#2`#3`#4;#5`#6`#7]{{
\settriparms[#1]
\diagram
\putptrianglep<\arrowtypea`\arrowtypeb`
\arrowtypec;\height>
(0,0)[#2`#3`#4;#5`#6`{#7}]
\enddiagram
}}                                            

\def\putqtrianglep<#1>(#2,#3)[#4`#5`#6;#7`#8`#9]{{%
\settriparms[#1]%
\xpos=#2 \ypos=#3
\advance\ypos by\height
\puthmorphism(\xpos,\ypos)[#4`#5`{#7}]{\height}{\arrowtypea}a%
\putmorphism(\xpos,\ypos)(1,-1)[``{#8}]{\height}{\arrowtypeb}l%
\advance\xpos by\height
\putvmorphism(\xpos,\ypos)[`#6`{#9}]{\height}{\arrowtypec}r%
}}

\def\putqtriangle{\@ifnextchar <{\putqtrianglep}{\putqtrianglep
   <\arrowtypea`\arrowtypeb`\arrowtypec;\height>}}
\def\qtriangle{\@ifnextchar <{\qtrianglep}{\qtrianglep
   <\arrowtypea`\arrowtypeb`\arrowtypec;\height>}}
\def\qtrianglep<#1>[#2`#3`#4;#5`#6`#7]{{
\settriparms[#1]
\width=\height                                
\diagram
\putqtrianglep<\arrowtypea`\arrowtypeb`
\arrowtypec;\height>
(0,0)[#2`#3`#4;#5`#6`{#7}]
\enddiagram
}}

\def\putdtrianglep<#1>(#2,#3)[#4`#5`#6;#7`#8`#9]{{%
\settriparms[#1]%
\xpos=#2 \ypos=#3
\puthmorphism(\xpos,\ypos)[#5`#6`{#9}]{\height}{\arrowtypec}b%
\advance\xpos by \height \advance\ypos by\height
\putmorphism(\xpos,\ypos)(-1,-1)[``{#7}]{\height}{\arrowtypea}l%
\putvmorphism(\xpos,\ypos)[#4``{#8}]{\height}{\arrowtypeb}r%
}}

\def\putdtriangle{\@ifnextchar <{\putdtrianglep}{\putdtrianglep
   <\arrowtypea`\arrowtypeb`\arrowtypec;\height>}}
\def\dtriangle{\@ifnextchar <{\dtrianglep}{\dtrianglep
   <\arrowtypea`\arrowtypeb`\arrowtypec;\height>}}
\def\dtrianglep<#1>[#2`#3`#4;#5`#6`#7]{{
\settriparms[#1]
\width=\height                                
\diagram
\putdtrianglep<\arrowtypea`\arrowtypeb`
\arrowtypec;\height>
(0,0)[#2`#3`#4;#5`#6`{#7}]
\enddiagram
}}

\def\putbtrianglep<#1>(#2,#3)[#4`#5`#6;#7`#8`#9]{{%
\settriparms[#1]%
\xpos=#2 \ypos=#3
\puthmorphism(\xpos,\ypos)[#5`#6`{#9}]{\height}{\arrowtypec}b%
\advance\ypos by\height
\putmorphism(\xpos,\ypos)(1,-1)[``{#8}]{\height}{\arrowtypeb}r%
\putvmorphism(\xpos,\ypos)[#4``{#7}]{\height}{\arrowtypea}l%
}}

\def\putbtriangle{\@ifnextchar <{\putbtrianglep}{\putbtrianglep
   <\arrowtypea`\arrowtypeb`\arrowtypec;\height>}}
\def\btriangle{\@ifnextchar <{\btrianglep}{\btrianglep
   <\arrowtypea`\arrowtypeb`\arrowtypec;\height>}}
\def\btrianglep<#1>[#2`#3`#4;#5`#6`#7]{{
\settriparms[#1]
\width=\height                               
\diagram
\putbtrianglep<\arrowtypea`\arrowtypeb`
\arrowtypec;\height>
(0,0)[#2`#3`#4;#5`#6`{#7}]
\enddiagram
}}

\def\putAtrianglep<#1>(#2,#3)[#4`#5`#6;#7`#8`#9]{{%
\settriparms[#1]%
\xpos=#2 \ypos=#3
{\multiply \height by2
\puthmorphism(\xpos,\ypos)[#5`#6`{#9}]{\height}{\arrowtypec}b}%
\advance\xpos by\height \advance\ypos by\height
\putmorphism(\xpos,\ypos)(-1,-1)[#4``{#7}]{\height}{\arrowtypea}l%
\putmorphism(\xpos,\ypos)(1,-1)[``{#8}]{\height}{\arrowtypeb}r%
}}

\def\putAtriangle{\@ifnextchar <{\putAtrianglep}{\putAtrianglep
   <\arrowtypea`\arrowtypeb`\arrowtypec;\height>}}
\def\Atriangle{\@ifnextchar <{\Atrianglep}{\Atrianglep
   <\arrowtypea`\arrowtypeb`\arrowtypec;\height>}}
\def\Atrianglep<#1>[#2`#3`#4;#5`#6`#7]{{
\settriparms[#1]
\width=\height                                     
\diagram
\putAtrianglep<\arrowtypea`\arrowtypeb`
\arrowtypec;\height>
(0,0)[#2`#3`#4;#5`#6`{#7}]
\enddiagram
}}

\def\putAtrianglepairp<#1>(#2)[#3;#4`#5`#6`#7`#8]{{%
\settripairparms[#1]%
\setpos(#2)%
\settokens[#3]%
\puthmorphism(\xpos,\ypos)[\tokenb`\tokenc`{#7}]{\height}{\arrowtyped}b%
\advance\xpos by\height
\puthmorphism(\xpos,\ypos)[\phantom{\tokenc}`\tokend`{#8}]%
{\height}{\arrowtypee}b%
\advance\ypos by\height
\putmorphism(\xpos,\ypos)(-1,-1)[\tokena``{#4}]{\height}{\arrowtypea}l%
\putvmorphism(\xpos,\ypos)[``{#5}]{\height}{\arrowtypeb}m%
\putmorphism(\xpos,\ypos)(1,-1)[``{#6}]{\height}{\arrowtypec}r%
}}

\def\putAtrianglepair{\@ifnextchar <{\putAtrianglepairp}{\putAtrianglepairp%
   <\arrowtypea`\arrowtypeb`\arrowtypec`\arrowtyped`\arrowtypee;\height>}}
\def\Atrianglepair{\@ifnextchar <{\Atrianglepairp}{\Atrianglepairp%
   <\arrowtypea`\arrowtypeb`\arrowtypec`\arrowtyped`\arrowtypee;\height>}}

\def\Atrianglepairp<#1>[#2;#3`#4`#5`#6`#7]{{
\settripairparms[#1]
\settokens[#2]
\width=\height                                
\diagram
\putAtrianglepairp                            
<\arrowtypea`\arrowtypeb`\arrowtypec`
\arrowtyped`\arrowtypee;\height>
(0,0)[{#2};#3`#4`#5`#6`{#7}]
\enddiagram
}}

\def\putVtrianglep<#1>(#2,#3)[#4`#5`#6;#7`#8`#9]{{%
\settriparms[#1]%
\xpos=#2 \ypos=#3
\advance\ypos by\height
{\multiply\height by2
\puthmorphism(\xpos,\ypos)[#4`#5`{#7}]{\height}{\arrowtypea}a}%
\putmorphism(\xpos,\ypos)(1,-1)[`#6`{#8}]{\height}{\arrowtypeb}l%
\advance\xpos by\height
\advance\xpos by\height
\putmorphism(\xpos,\ypos)(-1,-1)[``{#9}]{\height}{\arrowtypec}r%
}}

\def\putVtriangle{\@ifnextchar <{\putVtrianglep}{\putVtrianglep
   <\arrowtypea`\arrowtypeb`\arrowtypec;\height>}}
\def\Vtriangle{\@ifnextchar <{\Vtrianglep}{\Vtrianglep
   <\arrowtypea`\arrowtypeb`\arrowtypec;\height>}}
\def\Vtrianglep<#1>[#2`#3`#4;#5`#6`#7]{{
\settriparms[#1]
\width=\height                                 
\diagram
\putVtrianglep<\arrowtypea`\arrowtypeb`
\arrowtypec;\height>
(0,0)[#2`#3`#4;#5`#6`{#7}]
\enddiagram
}}

\def\putVtrianglepairp<#1>(#2)[#3;#4`#5`#6`#7`#8]{{
\settripairparms[#1]%
\setpos(#2)%
\settokens[#3]%
\advance\ypos by\height
\putmorphism(\xpos,\ypos)(1,-1)[`\tokend`{#6}]{\height}{\arrowtypec}l%
\puthmorphism(\xpos,\ypos)[\tokena`\tokenb`{#4}]{\height}{\arrowtypea}a%
\advance\xpos by\height
\puthmorphism(\xpos,\ypos)[\phantom{\tokenb}`\tokenc`{#5}]%
{\height}{\arrowtypeb}a%
\putvmorphism(\xpos,\ypos)[``{#7}]{\height}{\arrowtyped}m%
\advance\xpos by\height
\putmorphism(\xpos,\ypos)(-1,-1)[``{#8}]{\height}{\arrowtypee}r%
}}

\def\putVtrianglepair{\@ifnextchar <{\putVtrianglepairp}{\putVtrianglepairp%
    <\arrowtypea`\arrowtypeb`\arrowtypec`\arrowtyped`\arrowtypee;\height>}}
\def\Vtrianglepair{\@ifnextchar <{\Vtrianglepairp}{\Vtrianglepairp%
    <\arrowtypea`\arrowtypeb`\arrowtypec`\arrowtyped`\arrowtypee;\height>}}
\def\Vtrianglepairp<#1>[#2;#3`#4`#5`#6`#7]{{
\settripairparms[#1]
\settokens[#2]
\diagram
\putVtrianglepairp                             
<\arrowtypea`\arrowtypeb`\arrowtypec`
\arrowtyped`\arrowtypee;\height>
(0,0)[{#2};#3`#4`#5`#6`{#7}]
\enddiagram
}}

\def\putCtrianglep<#1>(#2,#3)[#4`#5`#6;#7`#8`#9]{{%
\settriparms[#1]%
\xpos=#2 \ypos=#3
\advance\ypos by\height
\putmorphism(\xpos,\ypos)(1,-1)[``{#9}]{\height}{\arrowtypec}l%
\advance\xpos by\height
\advance\ypos by\height
\putmorphism(\xpos,\ypos)(-1,-1)[#4`#5`{#7}]{\height}{\arrowtypea}l%
{\multiply\height by 2
\putvmorphism(\xpos,\ypos)[`#6`{#8}]{\height}{\arrowtypeb}r}%
}}

\def\putCtriangle{\@ifnextchar <{\putCtrianglep}{\putCtrianglep
    <\arrowtypea`\arrowtypeb`\arrowtypec;\height>}}
\def\Ctriangle{\@ifnextchar <{\Ctrianglep}{\Ctrianglep
    <\arrowtypea`\arrowtypeb`\arrowtypec;\height>}}
\def\Ctrianglep<#1>[#2`#3`#4;#5`#6`#7]{{
\settriparms[#1]
\width=\height                               
\diagram
\putCtrianglep<\arrowtypea`\arrowtypeb`
\arrowtypec;\height>
(0,0)[#2`#3`#4;#5`#6`{#7}]
\enddiagram
}}                                           
\def\putDtrianglep<#1>(#2,#3)[#4`#5`#6;#7`#8`#9]{{%
\settriparms[#1]%
\xpos=#2 \ypos=#3
\advance\xpos by\height \advance\ypos by\height
\putmorphism(\xpos,\ypos)(-1,-1)[``{#9}]{\height}{\arrowtypec}r%
\advance\xpos by-\height \advance\ypos by\height
\putmorphism(\xpos,\ypos)(1,-1)[`#5`{#8}]{\height}{\arrowtypeb}r%
{\multiply\height by 2
\putvmorphism(\xpos,\ypos)[#4`#6`{#7}]{\height}{\arrowtypea}l}%
}}

\def\putDtriangle{\@ifnextchar <{\putDtrianglep}{\putDtrianglep
    <\arrowtypea`\arrowtypeb`\arrowtypec;\height>}}
\def\Dtriangle{\@ifnextchar <{\Dtrianglep}{\Dtrianglep
   <\arrowtypea`\arrowtypeb`\arrowtypec;\height>}}
\def\Dtrianglep<#1>[#2`#3`#4;#5`#6`#7]{{
\settriparms[#1]
\width=\height                              
\diagram
\putDtrianglep<\arrowtypea`\arrowtypeb`
\arrowtypec;\height>
(0,0)[#2`#3`#4;#5`#6`{#7}]
\enddiagram
}}                                          
\def\setrecparms[#1`#2]{\width=#1 \height=#2}%

\def\recursep<#1`#2>[#3;#4`#5`#6`#7`#8]{{%
\width=#1 \height=#2
\settokens[#3]
\settowidth{\tempdimen}{$\tokena$}
\ifdim\tempdimen=0pt
  \savebox{\tempboxa}{\hbox{$\tokenb$}}%
  \savebox{\tempboxb}{\hbox{$\tokend$}}%
  \savebox{\tempboxc}{\hbox{$#6$}}%
\else
  \savebox{\tempboxa}{\hbox{$\hbox{$\tokena$}\times\hbox{$\tokenb$}$}}%
  \savebox{\tempboxb}{\hbox{$\hbox{$\tokena$}\times\hbox{$\tokend$}$}}%
  \savebox{\tempboxc}{\hbox{$\hbox{$\tokena$}\times\hbox{$#6$}$}}%
\fi
\ypos=\height
\divide\ypos by 2
\xpos=\ypos
\advance\xpos by \width
\bfig
\putCtrianglep<-1`1`1;\ypos>(0,0)[`\tokenc`;#5`#6`{#7}]%
\puthmorphism(\ypos,0)[\tokend`\usebox{\tempboxb}`{#8}]{\width}{-1}b%
\puthmorphism(\ypos,\height)[\tokenb`\usebox{\tempboxa}`{#4}]{\width}{-1}a%
\advance\ypos by \width
\putvmorphism(\ypos,\height)[``\usebox{\tempboxc}]{\height}1r%
\efig
}}

\def\recurse{\@ifnextchar <{\recursep}{\recursep<\width`\height>}}

\def\puttwohmorphisms(#1,#2)[#3`#4;#5`#6]#7#8#9{{%
\puthmorphism(#1,#2)[#3`#4`]{#7}0a
\ypos=#2
\advance\ypos by 20
\puthmorphism(#1,\ypos)[\phantom{#3}`\phantom{#4}`#5]{#7}{#8}a
\advance\ypos by -40
\puthmorphism(#1,\ypos)[\phantom{#3}`\phantom{#4}`#6]{#7}{#9}b
}}

\def\puttwovmorphisms(#1,#2)[#3`#4;#5`#6]#7#8#9{{%
\putvmorphism(#1,#2)[#3`#4`]{#7}0a
\xpos=#1
\advance\xpos by -20
\putvmorphism(\xpos,#2)[\phantom{#3}`\phantom{#4}`#5]{#7}{#8}l
\advance\xpos by 40
\putvmorphism(\xpos,#2)[\phantom{#3}`\phantom{#4}`#6]{#7}{#9}r
}}

\def\puthcoequalizer(#1)[#2`#3`#4;#5`#6`#7]#8#9{{%
\setpos(#1)%
\puttwohmorphisms(\xpos,\ypos)[#2`#3;#5`#6]{#8}11%
\advance\xpos by #8
\puthmorphism(\xpos,\ypos)[\phantom{#3}`#4`#7]{#8}1{#9}
}}

\def\putvcoequalizer(#1)[#2`#3`#4;#5`#6`#7]#8#9{{%
\setpos(#1)%
\puttwovmorphisms(\xpos,\ypos)[#2`#3;#5`#6]{#8}11%
\advance\ypos by -#8
\putvmorphism(\xpos,\ypos)[\phantom{#3}`#4`#7]{#8}1{#9}
}}

\def\putthreehmorphisms(#1)[#2`#3;#4`#5`#6]#7(#8)#9{{%
\setpos(#1) \settypes(#8)
\if a#9 %
     \vertsize{\tempcounta}{#5}%
     \vertsize{\tempcountb}{#6}%
     \ifnum \tempcounta<\tempcountb \tempcounta=\tempcountb \fi
\else
     \vertsize{\tempcounta}{#4}%
     \vertsize{\tempcountb}{#5}%
     \ifnum \tempcounta<\tempcountb \tempcounta=\tempcountb \fi
\fi
\advance \tempcounta by 60
\puthmorphism(\xpos,\ypos)[#2`#3`#5]{#7}{\arrowtypeb}{#9}
\advance\ypos by \tempcounta
\puthmorphism(\xpos,\ypos)[\phantom{#2}`\phantom{#3}`#4]{#7}{\arrowtypea}{#9}
\advance\ypos by -\tempcounta \advance\ypos by -\tempcounta
\puthmorphism(\xpos,\ypos)[\phantom{#2}`\phantom{#3}`#6]{#7}{\arrowtypec}{#9}
}}

\def\setarrowtoks[#1`#2`#3`#4`#5`#6]{%
\def\toka{#1}
\def\tokb{#2}
\def\tokc{#3}
\def\tokd{#4}
\def\toke{#5}
\def\tokf{#6}
}
\def\hex{\@ifnextchar <{\hexp}{\hexp<1000`400>}}
\def\hexp<#1`#2>[#3`#4`#5`#6`#7`#8;#9]{%
\setarrowtoks[#9]
\yext=#2 \advance \yext by #2
\xext=#1 \advance\xext by \yext
\bfig
\putCtriangle<-1`0`1;#2>(0,0)[`#5`;\tokb``\tokd]
\xext=#1 \yext=#2 \advance \yext by #2
\putsquare<1`0`0`1;\xext`\yext>(#2,0)[#3`#4`#7`#8;\toka```\tokf]
\advance \xext by #2
\putDtriangle<0`1`-1;#2>(\xext,0)[`#6`;`\tokc`\toke]
\efig
}

{\catcode`\ =13\global\let =\ \catcode`\^^M=13
\gdef^^M{\par\noindent}}
\newtheorem{thm_remark}[theorem]{Remark}
\begin{document}

\author{R. F. Picken \\
Departamento de Matem\'{a}tica and\\
CEMAT,\\
Centro de Matem\'{a}tica e Aplica\c{c}\~{o}es\\
Instituto Superior T\'{e}cnico,\\
Avenida Rovisco Pais\\
1049-001 Lisboa, Portugal\\
e-mail: picken@math.ist.utl.pt\\
and\\
P. A. Semi\~{a}o \\
\'{A}rea Departamental de Matem\'{a}tica and\\
CEMAT, \\
Centro de Matem\'{a}tica e Aplica\c{c}\~{o}es\\
Universidade do Algarve,\\
Faculdade de Ci\^{e}ncias e Tecnologia,\\
8000-062 Faro, Portugal\\
e-mail: psemiao@ualg.pt}
\title{{\Large \textbf{A\ classical approach to TQFT's\thanks{%
This work was supported by the programme ``\emph{Programa Operacional
Ci\^{e}ncia, Tecnologia, Inova\c{c}\~{a}o}'' (POCTI) of the \emph{Funda\c{c}%
\~{a}o para a Ci\^{e}ncia e Tecnologia} (FCT), cofinanced by the European
Community fund FEDER.}}}}
\maketitle

\begin{abstract}
We present a general framework for TQFT and related constructions using the
language of monoidal categories. We construct a topological category $%
\mathcal{C}$ and an algebraic category $\mathcal{D}$, both monoidal, and a
TQFT functor is then defined as a certain type of monoidal functor from $%
\mathcal{C}$ to $\mathcal{D}$. In contrast with the cobordism approach, this
formulation of TQFT is closer in spirit to the classical functors of
algebraic topology, like homology. The fundamental operation of gluing is
incorporated at the level of the morphisms in the topological category
through the notion of a gluing morphism, which we define. It allows not only
the gluing together of two separate objects, but also the self-gluing of a
single object to be treated in the same fashion. As an example of our
framework we describe TQFT's for oriented $2\mathrm{D}$-manifolds, and
classify a family of them in terms of a pair of tensors satisfying some
relations.\vspace{0.5cm}

\begin{flushright}
IST/DM/xx/02\\
math.AT/0212310\\
\end{flushright}
\end{abstract}

\tableofcontents

\section{Introduction}

The term Topological Quantum Field Theory, or TQFT for short, was introduced
by Witten in \cite{Witten2} to describe a class of quantum field theories
whose action is diffeomorphism invariant. To capture the common features of
a number of examples that appeared, Atiyah, in a seminal article \cite
{Atiyah}, formulated a set of axioms for TQFT that were modelled on similar
axioms for conformal field theory, due to Segal \cite{Segal}. A $(d+1)$%
-dimensional TQFT assigns to every closed oriented $d$-dimensional manifold $%
A$ a finite-dimensional vector space $V_{A}$, and to every $(d+1)$%
-dimensional oriented manifold with boundary $X$ an element $Z_{X}$ of the
vector space assigned to its boundary, subject to certain rules. The notion
of TQFT has had a pervasive influence in several areas of mathematics, in
particular in low-dimensional topology. Amongst the many constructions
coming out of TQFT, we would single out the $3$-manifold invariants due to
Witten \cite{Witten1} and Reshetikhin-Turaev \cite{Reshetikhin-Turaev}, and
Turaev-Viro \cite{Turaev-Viro}, the combinatorial formula for the signature
of $4$-manifolds due to Crane-Yetter \cite{Crane-Yetter}, and the
Dijkgraaf-Witten invariants for manifolds of any dimension, coming from a
discretized path integral for gauge theories with finite group \cite
{Dijkgraaf-Witten}. For a review of TQFT aimed at mathematicians, see \cite
{Sawin2}.

As recognized in Atiyah's article, there are many possible modifications of
the basic axioms. Indeed the Witten invariant from \cite{Witten1} does not
fit into the axioms, since it is a combination of a 3-manifold invariant and
an invariant for an embedded $1$-manifold. Representations of braids \cite
{Birman} and tangles \cite{Turaev} also fit loosely into the TQFT framework,
with the difference that they are trivial as $1$-dimensional manifolds, but
non-trivial due to their embedding in $3\mathrm{D}$ space. Generalizations
arise when one allows the manifolds to have extra geometrical structures. An
example is Homotopy Quantum Field Theory (HQFT) \cite
{TuraevHQFT1,TuraevHQFT2} (see also \cite{Brightwell-Turner} and \cite
{Goncalo}), where the manifolds come equipped with a map to a fixed target
space.

Turning to the algebraic side, many interesting ideas have appeared about
the nature of the vector spaces and elements to be assigned to the
manifolds. These ideas constitute, in a sense, extra structure, similar to
the possible extra structure on the topological side considered above.
Whilst not wishing to summarize all constructions, we would like to mention
as examples, representations of the quantum group $SL(2)_{q}$ at roots of
unity for $3$-dimensional TQFT's \cite{Reshetikhin-Turaev}, spherical
categories for $3$-dimensional TQFT's \cite{Barrett-Westbury}, Hopf
categories \cite{Crane-Frenkel,Crane-Yetter2} and spherical $2$-categories
\cite{Mackaay1,Mackaay2} for $4$-dimensional TQFT's, constructions coming
from operator algebras \cite{Evans-KawBook}, and the very interesting
higher-dimensional algebra programme \cite{Baez-Dolan,Baez-LangfordIV} in
which notions of higher algebra ($n$-categories) enter simultaneously on the
topological and algebraic side.

To complete this brief survey of variations of TQFT, there are cases where
both the topological and algebraic side come from geometry, namely parallel
transport for vector bundles with connection. Here the underlying structure
is a vector bundle $(E,\pi ,M)$ with connection, and the $1$-dimensional
TQFT assigns to each point $a$ of the base space $M$ its fibre $\pi ^{-1}(a)$%
, and to each path from $a$ to $b$ the parallel transport operator from the
fibre over $a$ to the fibre over $b$, regarded as an element of $\pi
^{-1}(a)^{*}\otimes \pi ^{-1}(b)$. This is essentially the functorial
viewpoint of connections given in \cite{Segal}. An interpretation in terms
of HQFT was recently given for parallel transport for gerbes with connection
\cite{Bunke-Turner-Willerton}.

As stated earlier, in his article on the axioms of TQFT Atiyah already
envisaged various modifications of the axioms to enlarge the scope, so it
was natural that mathematicians would attempt to formulate the theory in a
more general fashion. There are two such general definitions of TQFT that we
would like to mention in this context, due to Quinn \cite{QuinnInBook,Quinn2}
and Turaev \cite{TuraevBook}. Quinn achieved generality by giving a very
general definition of a ``domain category'', endowed with a collection of
structures which are abstractions of topological notions, such as boundary,
cylinder or gluing. Turaev achieved generality by defining a ``space
structure'' on a topological space, which includes as special cases a choice
of orientation, a differentiable structure or a structure of a $\mathrm{CW}$%
-complex. Both definitions are rather abstract, and both are in the
so-called cobordism framework, which has some limitations as we now explain.

The cobordism approach to TQFT's arises by taking an algebraic idea and
transporting it to the topological side. Suppose the manifold $X$ has an
in-boundary $A_{1}^{-}$ and an out-boundary $A_{2}^{+}$, like the starting
and end point of the curve on $M$ in the parallel transport example
mentioned above. To this boundary one assigns the tensor product $%
V_{A_{1}^{-}}\otimes V_{A_{2}^{+}}$, which is assumed to be isomorphic to $%
V_{A_{1}^{+}}^{*}\otimes V_{A_{2}^{+}}$. Thus $Z_{X}$ can be regarded as a
linear map from $V_{A_{1}^{+}}$ to $V_{A_{2}^{+}}$, i.e. a morphism in the
algebraic category. In the cobordism approach, $X$ is then also taken to be
a morphism from $A_{1}$ to $A_{2}$, in a category on the topological side
called the cobordism category. More precisely, appropriate equivalence class
of manifolds $X$ constitute the morphisms in the cobordism category, so as
to ensure the associativity of composition and the existence of identity
morphisms. An advantage of this shift of viewpoint is that the TQFT
assignments become functorial both at the level of $X$'s and $A$'s, and thus
a TQFT becomes a functor from the cobordism category to the category of
vector spaces. The cobordism viewpoint is very natural for $1$-dimensional
TQFT's describing braids and tangles (these have been termed ``embedded
TQFT's'' in \cite{Roger}), since morphisms are intuitively associated to $1$%
-dimensional topology, but is less natural, or rather, has to be amended to
higher cobordism categories in the higher-dimensional algebra setting, when
one is dealing with dimensions greater than one. Despite the attractive
functorial feature of the cobordism approach, a disadvantage is that making
the $X$'s into morphisms on the topological side leaves no room for other
morphisms between the $A$'s, in particular isomorphisms between them. (In
the context of HQFT, Rodrigues in \cite{Goncalo} showed a nice way to get
round this problem, by means of some quotient constructions, allowing both
cobordisms and isomorphisms to appear as topological morphisms on the same
footing.) A second disadvantage of the cobordism approach is that, although
the gluing together of two manifolds $X_{1}$ and $X_{2}$ along a shared
boundary comes in naturally as the composition of topological morphisms, the
notion of self-gluing of a single manifold $X$ by identifying two of its
boundary components has no such natural interpretation, and has to be dealt
with in theorems. A third disadvantage of the cobordism approach is that the
algebraic notion of duality, which was necessary to pass from elements of a
linear space to linear maps, requires a rather subtle treatment when it is
transported to the topological side - see \cite{Goncalo} for a definition in
the context of HQFT and further references on the subject of duality in
categories. Finally, a disadvantage, to our mind, of the cobordism approach
is that the same topological object acquires a plethora of algebraic guises.
Thus, for example, disregarding orientations, the cylinder can be viewed in
three different ways as a morphism from the circle to the circle, or a
morphism from two circles to the empty set, or a morphism from the empty set
to two circles.

Thus we were motivated to return to Atiyah's original formulation where $%
Z_{X}$ is simply an element of $V_{A}$. Now Atiyah's axioms in this form do
not describe a functor, although they do contain functorial ingredients. One
of our main goals was to define TQFT's as functors, analogously to the
classical functors of algebraic topology such as homology and homotopy
(hence our slightly tongue-in-cheek title). The approach we found worked
best was to take objects on the topological side to be essentially pairs $%
(X,A)$ (together with an ``inclusion'' morphism from $A$ to $X$), and to
allow only a very restricted class of morphisms between these objects,
namely isomorphisms and gluing morphisms, where the latter are essentially
morphisms from an object to a copy of the object that results after gluing
some components of $A$ together. A TQFT is then a functor from this category
to an algebraic category whose objects are also pairs, essentially a vector
space and an element of that space. Furthermore, the TQFT functors preserve
structures, namely a monoidal structure (which is essentially the disjoint
union on the topological side) as well as an endofunctor (which is an
abstraction of the operation of changing orientation on the topological
side). Thus we use the language of monoidal categories and monoidal
functors. Generality in describing a wide range of topological situations is
achieved by building up the topological category, in a series of steps, out
of a fairly arbitrary topological starting category.

Our approach is rather detailed, since we wanted to arrive at a very
concrete and explicit formulation, allowing calculations to be performed
efficiently. The insight that is gained, after working through the
formalism, is that TQFT is in essence a ``topological tensor calculus'',
where manifolds get tensors assigned to them, with the number of indices
equal to the number of connected components of $A$, and where gluing
manifolds together corresponds to contracting the indices. This picture is
made explicit in our illustrative example of $2$-dimensional TQFT's, where
we obtain a result precisely of this nature: all $2$-manifolds with boundary
are represented by a tensor built up from just two basic tensors, one with
one index and one with three indices, representing the disk and the
pair-of-pants, respectively. This result should be contrasted with Abrams'
\cite{LowellAbrams} characterization of $2$-dimensional TQFT's in terms of
Frobenius algebras, using the cobordism approach - see also earlier work by
Sawin \cite{Sawin1}. There the topological category ($2$-Cobord) is
generated by five surfaces, two different disks, two different
pairs-of-pants and one cylinder. The first four surfaces correspond to the
(co)unit and (co)multiplication of the Frobenius algebra structure on the
vector space assigned to the circle, which is again an example of
topological objects acquiring multiple algebraic guises. Our result for $2$%
-dimensional TQFT's is much closer in spirit to that of R. Lawrence \cite
{RuthLawrence}, who gets a classification in terms of three tensors,
corresponding to the disk, pair-of-pants and cylinder, respectively.

Rather than give an outline of the whole construction here, we refer to the
introduction to each section for a description of the various stages
involved. Since the whole setting is that of monoidal categories, we have
included some relevant definitions in the introduction - see below. The
general theory is interspersed with a detailed presentation of an example,
namely $2$-dimensional TQFT's, which are characterized in Theorems \ref
{thm44} and \ref{thm56}. A preliminary version of this work appeared in \cite
{PickSemiao} and it was the subject of the Ph.D. thesis \cite{PSemiaoPhD}.
We are very grateful to Louis Crane for his suggestions and comments on an
earlier version of this article.

To complete this introduction we recall some definitions concerning monoidal
categories. For the background material on monoidal categories we refer the
reader to the textbooks \cite[pp. 161-170]{MacLane1} and \cite[pp. 281-288]
{Kassel}:

\begin{definition}
A \emph{monoidal category} is a sextuple $(\mathbf{C},\otimes ,I,a,l,r)$
consisting of a category $\mathbf{C}$, a bifunctor $\otimes :\mathbf{C}%
\times \mathbf{C}\rightarrow \mathbf{C}$, called the monoidal product, an
object $I\in \func{Ob}(\mathbf{C})$, called the unit object, and natural
isomorphisms $a:(-\otimes -)\otimes -\rightarrow -\otimes (-\otimes -)$
(associativity), $l:I\otimes -\rightarrow \limfunc{Id}\nolimits_{\mathbf{C}}$
(left unit) and $r:-\otimes I\rightarrow \limfunc{Id}\nolimits_{\mathbf{C}}$
(right unit), such that for all objects $A,B,C$ and $D$ the following
diagrams:%
%
\begin{center}
\xext=2600
\yext=700
\begin{picture}(\xext,\yext)(\xoff,\yoff)
\putmorphism(0,500)(1,0)[((A\otimes B)\otimes C)\otimes D`(A\otimes B)\otimes (C\otimes D)`a_{A\otimes B,C,D}]{1300}1a
\putmorphism(1300,500)(1,0)[\phantom{(A\otimes B)\otimes (C\otimes D)}`A\otimes (B\otimes (C\otimes D))`a_{A,B,C\otimes D}]{1300}1a
\putmorphism(0,500)(0,-1)[\phantom{((A\otimes B)\otimes C)\otimes D}`\phantom{(A\otimes B)\otimes (C\otimes D}`a_{A,B,C}\otimes id_{D}]{500}1l
\putmorphism(2400,500)(0,-1)[\phantom{A\otimes (B\otimes (C\otimes D))}`\phantom{(A\otimes B)\otimes (C\otimes D}`id_{A}\otimes a_{B,C,D}]{500}{-1}r
\putmorphism(0,0)(1,0)[(A\otimes (B\otimes C))\otimes D`A\otimes ((B\otimes C)\otimes D)`a_{A,B\otimes C,D}]{2400}1b
\end{picture}
\end{center}

%
and%
\begin{center}
\Vtriangle[(A\otimes I)\otimes B`A\otimes(I\otimes B)`A\otimes B;a_{A,I,B}`r_{A}\otimes id_{B}`id_{A}\otimes l_{B}]
\end{center}
%
commute. A monoidal category is often denoted simply by $\mathbf{C}$.\newline
A monoidal category is said to be \emph{strict}, if the associativity, left
and right unit are all identities of the category. A strict monoidal
category is denoted by $(\mathbf{C},\otimes ,I)$ or $\mathbf{C}$.
\end{definition}

\begin{thm_remark}
According to a result by Mac Lane, any monoidal category is monoidal
equivalent to a strict monoidal category. (See \cite[pp. 256-257]{MacLane1}
for the definition of monoidal equivalence and the proof. See also
\cite[p. 291]{Kassel}.)\newline
In the constructions of Sections \ref{section2}-\ref{section5} below, only
strict monoidal categories will appear.
\end{thm_remark}

\begin{definition}
\label{Def1.3}Let $(\mathbf{C},\otimes ,I,a,l,r)$ and $(\mathbf{D},\otimes
^{\prime },I^{\prime },a^{\prime },l^{\prime },r^{\prime })$ be monoidal
categories. A \emph{monoidal functor} from $\mathbf{C}$ to $\mathbf{D}$ is a
triple $(F,\varphi _{2},\varphi _{0})$, which consists of a functor $F:%
\mathbf{C}\rightarrow \mathbf{D}$, a natural transformation $\varphi
_{2}:F(-)\otimes ^{\prime }F(-)\rightarrow F(-\otimes -)$ and a morphism $%
\varphi _{0}:I^{\prime }\rightarrow F(I)$, such that the following diagrams%
%
\begin{center}
\xext=1700
\yext=1150
\begin{picture}(\xext,\yext)(\xoff,\yoff)
\putmorphism(0,1000)(1,0)[(F(A)\otimes 'F(B))\otimes ' F(C)`F(A)\otimes ' (F(B)\otimes ' F(C))`a'_{F(A),F(B),F(C)}]{1800}1a
\putmorphism(0,1000)(0,-1)[\phantom{(F(A)\otimes 'F(B))\otimes ' F(C)}`F(A\otimes B)\otimes ' F(C)`\varphi_{2 (A,B)}\otimes ' id_{F(C)}]{500}1l
\putmorphism(0,500)(0,-1)[\phantom{A}`\phantom{A}`\varphi_{2 (A\otimes B,C)}]{500}1l
\putmorphism(1750,1000)(0,-1)[\phantom{F(A)\otimes ' (F(B)\otimes ' F(C))}`F(A)\otimes ' F(B\otimes C)`id_{F(A)}\otimes ' \varphi_{2 (B,C)}]{500}1r
\putmorphism(1750,500)(0,-1)[\phantom{A}`\phantom{A}`\varphi_{2 (A,B\otimes C)}]{500}1r
\putmorphism(0,0)(1,0)[F((A\otimes B)\otimes C)`F(A\otimes (B\otimes C))`F(a_{A,B,C})]{1800}1b
\end{picture}
\end{center}

%
and%
\begin{center}
\setsqparms[1`1`1`1;1000`700]
\square[I'\otimes ' F(A)`F(A)`F(I)\otimes' F(A)`F(I\otimes A);l'_{F(A)}`\varphi_{0}\otimes' id_{F(A)}`F(l_{A})`\varphi_{2 (I,A)}]
\end{center}
%
and%
\begin{center}
\setsqparms[1`1`1`1;1000`700]
\square[F(A)\otimes ' I'`F(A)`F(A)\otimes ' F(I)`F(A\otimes I);r'_{F(A)}`id_{F(A)}\otimes' \varphi_{0}`F(r_{A})`\varphi_{2 (A,I)}]
\end{center}
%
commute. We will frequently write a monoidal functor $(F,\varphi
_{2},\varphi _{0})$ simply as $F$.

The monoidal functor is said to be \emph{strong}, if the natural
transformation $\varphi _{2}$ is a natural isomorphism, and the morphism $%
\varphi _{0}$ is an isomorphism, and is said to be \emph{strict} if they are
identities.
\end{definition}

\begin{definition}
A \emph{monoidal natural transformation} $\eta :(F,\varphi _{2},\varphi
_{0})\rightarrow (G,\psi _{2},\psi _{0})$ is a natural transformation $\eta
:F\rightarrow G$ such that the following diagrams, for all objects $A,B$%
\begin{center}
\setsqparms[1`1`1`1;1000`700]
\square[F(A)\otimes ' F(B)`G(A)\otimes ' G(B)`F(A\otimes B)`G(A\otimes B);\eta_{A}\otimes ' \eta_{B}`\varphi_{2 (A,B)}`\psi_{2 (A,B)}`\eta_{A\otimes B}]
\end{center}
%
and%
\begin{center}
\settriparms[-1`1`1;400]
\Ctriangle[F(I)`I'`G(I);\varphi_{0}`\eta_{I}`\psi_{0}]
\end{center}%
commute.
\end{definition}

\begin{definition}
Let $(\mathbf{C},\otimes ,I,a,l,r)$ be a monoidal category and $\tau :%
\mathbf{C}\times \mathbf{C}\rightarrow \mathbf{C}\times \mathbf{C}$ the flip
functor, defined for all pairs of objects and morphisms by:
\[
\tau ((A,B)\stackrel{(f,g)}{\rightarrow }(A^{\prime },B^{\prime }))=(B,A)%
\stackrel{(g,f)}{\rightarrow }(B^{\prime },A^{\prime })\text{.}
\]
A \emph{braiding} is a natural isomorphism
\[
c:-\otimes -\rightarrow (-\otimes -)\circ \tau \text{,}
\]
such that, for all triples $(A,B,C)$ of objects of the category, the
following two hexagonal diagrams%
%
\begin{center}
\xext=2100
\yext=1150
\begin{picture}(\xext,\yext)(\xoff,\yoff)
\putmorphism(0,1000)(1,0)[A\otimes (B\otimes C)`(B\otimes C)\otimes A`c_{(A,B\otimes C)}]{1800}1a
\putmorphism(0,1000)(0,-1)[\phantom{A\otimes (B\otimes C)}`(A\otimes B)\otimes C`a_{A,B,C}]{500}{-1}l
\putmorphism(0,500)(0,-1)[\phantom{(A\otimes B)\otimes C}`\phantom{(B\otimes A)\otimes C}`c_{(A,B)}\otimes id_{C}]{500}1l
\putmorphism(1750,1000)(0,-1)[\phantom{(B\otimes C)\otimes A}`B\otimes (C\otimes A)`a_{B,C,A}]{500}1r
\putmorphism(1750,500)(0,-1)[\phantom{A}`\phantom{A}`id_{B}\otimes c_{(A,C)}]{500}{-1}r
\putmorphism(0,0)(1,0)[(B\otimes A)\otimes C`B\otimes (A\otimes C)`a_{B,A,C}]{1800}1b
\end{picture}
\end{center}

%
and%
%
\begin{center}
\xext=2100
\yext=1150
\begin{picture}(\xext,\yext)(\xoff,\yoff)
\putmorphism(0,1000)(1,0)[(A\otimes B)\otimes C`C\otimes (A\otimes B)`c_{(A\otimes B,C)}]{1800}1a
\putmorphism(0,1000)(0,-1)[\phantom{(A\otimes B)\otimes C}`A\otimes (B\otimes C)`a^{-1}_{A,B,C}]{500}{-1}l
\putmorphism(0,500)(0,-1)[\phantom{A\otimes (B\otimes C)}`\phantom{A\otimes (C\otimes B)}`id_{A}\otimes c_{(B,C)}]{500}1l
\putmorphism(1750,1000)(0,-1)[\phantom{C\otimes (A\otimes B)}`(C\otimes A)\otimes B`a^{-1}_{C,A,B}]{500}1r
\putmorphism(1750,500)(0,-1)[\phantom{(C\otimes A)\otimes B}`\phantom{(A\otimes C)\otimes B}`c_{(A,C)}\otimes id_{B}]{500}{-1}r
\putmorphism(0,0)(1,0)[A\otimes (C\otimes B)`(A\otimes C)\otimes B`a^{-1}_{A,B,C}]{1800}1b
\end{picture}
\end{center}

%
commute.\newline
A \emph{braided monoidal category} $\mathbf{C}$ is a monoidal category
endowed with a braiding, also denoted $(\mathbf{C},\otimes ,I,a,l,r,c)$, or $%
(\mathbf{C},\otimes ,I,c)$ if $\mathbf{C}$ is strict monoidal.\newline
A braided monoidal category $(\mathbf{C},\otimes ,I,a,l,r,c)$ is called a
\emph{symmetric monoidal category} if the braiding $c$ satisfies:
\[
c_{(B,A)}\circ c_{(A,B)}=\limfunc{id}\nolimits_{A\otimes B}\text{,}
\]
for all pairs $(A,B)$ of objects.
\end{definition}

\begin{thm_remark}
We note that the condition $c_{(B,A)}\circ c_{(A,B)}=\limfunc{id}%
\nolimits_{A\otimes B}$ is equivalent to $c_{(A,B)}=c_{(B,A)}^{-1}$, which
implies that the second hexagonal diagram is equal to the first. Therefore a
symmetric monoidal category is a monoidal category in which the first
hexagonal diagram commutes and the symmetry condition holds.
\end{thm_remark}

Finally, we present the definition of a braided monoidal functor.

\begin{definition}
Let $(\mathbf{C},\otimes ,I,a,l,r,c)$ and $(\mathbf{D},\otimes ^{\prime
},I^{\prime },a^{\prime },l^{\prime },r^{\prime },c^{\prime })$ be braided
monoidal categories. A monoidal functor $(F,\varphi _{2},\varphi _{0})$ from
$\mathbf{C}$ to $\mathbf{D}$ is said to be a \emph{braided monoidal functor}
if, for all pairs of objects $(A,B)$, the following diagram%
\begin{center}
\setsqparms[1`1`1`1;1000`600]
\square[F(A)\otimes ' F(B)`F(A\otimes B)`F(B)\otimes ' F(A)`F(B\otimes A);\varphi_{2 (A,B)}`c'_{(F(A),F(B))}`F(c_{(A,B)})`\varphi_{2 (B,A)}]
\end{center}
%
commutes.

A \emph{symmetric monoidal functor} is a braided monoidal functor between
symmetric monoidal categories.
\end{definition}

\section{The topological category\label{section2}}

The construction of the topological category will proceed in stages, which
we will first describe informally. The starting point is a category which is
large enough to contain all the objects one wishes to describe, called the
topological starting category. This might be, for instance, the category of
all smooth oriented manifolds with boundary. Within this category we then
focus attention on a class of ``larger objects'' and a class of
``subobjects'', by choosing a class of ``inclusion'' morphisms whose domains
are the subobjects and whose codomains are the larger objects. For instance,
the larger objects could be $d$-dimensional manifolds with boundary, and the
subobjects ($d-1$)-dimensional manifolds with empty boundary, with the
inclusion morphisms mapping the subobjects into boundary components of the
larger objects. The subobjects form a category, called the category of
subobjects, with isomorphisms between subobjects as its morphisms. The
objects of the topological category itself are triples consisting of a
larger object, a subobject and an inclusion map between them. The morphisms
between these objects are restricted to be of two types only: isomorphisms
and gluing morphisms. The isomorphisms correspond to the respective larger
objects and subobjects being isomorphic in a compatible way. The gluing
morphisms correspond to gluing larger objects together along one or more
pairs of subobjects, and are best described as morphisms from the objects
before gluing to a copy of the objects after gluing. The main result that
has to be established is that these morphisms close under suitably defined
composition. There is an extra piece of structure which plays a crucial role
throughout the construction, namely an endofunctor on the topological
starting category, which generalizes the operation of changing orientation
in the example of oriented manifolds.

We now proceed with the details of the construction.

\begin{definition}
A \emph{topological starting category} is a triple $(\mathbf{C},F,P)$, where
$\mathbf{C}$ is a symmetric, strict monoidal category $(\mathbf{C},\sqcup
,E,c)$, $F$ is a symmetric, strict monoidal forgetful functor from $\mathbf{C%
}$ to $\mathrm{Top}$ (the category of topological spaces and continuous maps
with its standard monoidal structure and braiding) and $P$ is a symmetric,
strong monoidal endofunctor $(P,\pi _{2},\pi _{0})$ on $\mathbf{C}$, such
that
\[
F\circ P=F\hspace{0.2in}\text{and}\hspace{0.2in}F(\pi _{2})=\limfunc{id}%
\text{.}
\]
Instead of $(\mathbf{C},F,P)$ we will frequently write simply $\mathbf{C}$.
\end{definition}

\begin{thm_remark}
The forgetful functor $F$ justifies the adjective topological. It sends each
object $X$ of $\mathbf{C}$ to its underlying topological space $F(X)$, sends
each morphism of $\mathbf{C}$ to the underlying continuous map, sends $%
\sqcup $ to the disjoint union (or topological sum) of topological spaces,
and sends $E$ to the empty space $\emptyset $.\newline
One may think of the objects of $\mathbf{C}$ as topological spaces with some
extra structure, e.g. oriented topological manifolds. The endofunctor $P$
acts on the extra structure, e.g. by changing the orientation, but leaves
the underlying topological space unchanged. We do not require $P^{2}=%
\limfunc{Id}$, although the only examples we have in mind do have this
property.\medskip
\end{thm_remark}

We will illustrate the constructions in Sections 2-5 with a simple example,
which will keep returning.\medskip

{\noindent \bf Example: }%
We take $\mathbf{C}$ to be the category whose objects are $1$-dimensional,
compact, oriented, topological manifolds without boundary or $2$%
-dimensional, compact, oriented, topological manifolds with boundary%
\footnote{%
The category of manifolds with boundary is more accurately described in
terms of a category whose objects are pairs $(X,\partial X)$, where $%
\partial X$ is the boundary of the manifold $X$. In our example, $S^{2}$,
for instance, is regarded as a manifold with boundary corresponding to the
pair $(S^{2},\emptyset )$.}, together with finite disjoint unions of these
objects. We introduce the following symbols and a preferred presentation for
some familiar objects:

\begin{itemize}
\item  $C_{+}$ denotes the circle, identified with $\left\{ z\in \Bbb{C}%
:\left| z\right| =1\right\} $ with anticlockwise orientation.

\item  $D_{+}$ denotes the disk, identified with $\left\{ z\in \Bbb{C}%
:\left| z\right| \leq 1\right\} $ with the orientation induced by the
standard orientation of the complex plane (given by the choice of
coordinates $(x,y)$, where $z=x+iy$).

\item  $A_{+}$ denotes the annulus, identified with $\left\{ z\in \Bbb{C}%
:1\leq \left| z\right| \leq 2\right\} $ with the orientation induced by the
standard orientation of the complex plane.

\item  $P_{+}$ denotes the pair-of-pants (i.e. the three holed sphere or
trinion) identified with a disk in the complex plane with two holes removed,
with the orientation induced by the standard orientation of the complex
plane.

\item  $S_{+}$ denotes the $2$-sphere, $S^{2}$, identified with the unit
sphere in $\Bbb{R}^{3}$ and with the orientation given by the choice of
coordinates $(x,y)$ at $(0,0,1)$.

\item  Finally, $T_{+}$ denotes the torus $S^{1}\times S^{1}$, identified
with\linebreak $\left\{ (e^{i\theta },e^{i\varphi })\in \Bbb{C}^{2}:0\leq
\theta ,\varphi <2\pi \right\} $ and with the orientation given by the
choice of coordinates $(\theta ,\varphi )$.
\end{itemize}

For each of these objects replacing the ``$+$'' suffix by a ``$-$'' denotes
the same object with the opposite orientation. We will abbreviate a disjoint
union of the form $X_{+}\sqcup X_{-}\sqcup \cdots \sqcup X_{+}$ by $%
X_{+-\cdots +}$. We will equate $X_{1}\sqcup X_{2}\sqcup \cdots \sqcup X_{n}$
with $X_{1}\times \left\{ 1\right\} \cup X_{2}\times \left\{ 2\right\} \cup
\cdots \cup X_{n}\times \left\{ n\right\} $. The objects:
\[
\emptyset ,\,\,C_{\pm },D_{\pm },P_{\pm }\text{,}
\]
play a fundamental role, since, up to isomorphism, all other objects can be
obtained from these by disjoint union and gluing, as we shall see later.

The morphisms of $\mathbf{C}$ are orientation-preserving maps between the
above objects, for instance $f:C_{+}\rightarrow C_{+}$, given by $f(z)=z$,
or $f:C_{-}\rightarrow C_{+}$, given by $f(z)=\overline{z}$. Since the
orientation of a manifold with boundary induces an orientation on its
boundary, the morphisms of $\mathbf{C}$ include orientation-preserving maps
from the oriented circle to oriented $2$-dimensional manifolds with boundary
whose image is contained in the boundary of the $2$-dimensional manifold,
such as $f:C_{+}\rightarrow D_{+}$, given by $f(z)=z$. The monoidal product $%
\sqcup $ is the disjoint union, taken to be strictly associative, the unit
object $E$ is $\emptyset $, also regarded as strict, i.e. $X\sqcup \emptyset
=\emptyset \sqcup X=X$, and the braiding $c_{(X,Y)}:X\sqcup Y\rightarrow
Y\sqcup X$ is the usual flip map which sends $(x,1)$ to $(x,2)$ and $(y,2)$
to $(y,1)$. The forgetful functor $F$ maps each oriented manifold to its
underlying topological space and each orientation-preserving map to the
underlying continuous map. Finally, the endofunctor $P$ reverses the
orientation of the objects and leaves the morphisms unchanged in the sense
that $P(f)$ and $f$ have the same transformation law, for any morphism $f$.
Thus the endofunctor $P$ is involutive, i.e. $P^{2}=\limfunc{Id}$. The
endofunctor $P$ is strict in this example, meaning that $\pi _{2}$ and $\pi
_{0}$ are identities. Thus, in particular, $P$ leaves the unit object $%
\emptyset $ unchanged.$\blacktriangle $\medskip

Having chosen a topological starting category $\mathbf{C}$, our next goal is
to introduce a subcategory of $\mathbf{C}$, which we call the category of
subobjects of $\mathbf{C}$. The idea is to separate the objects into
``larger'' objects and ``smaller'' objects, the latter behaving like the
boundary components of the larger objects. For this purpose we first
formalize the notion of a subobject following \cite{Adamek}.

\begin{definition}
\label{Def23}Let $\mathcal{M}$ be a non-empty class of monomorphisms of $%
\mathbf{C}$. An $\mathcal{M}$-\emph{subobject} of an object $X$ is a pair $%
(A,m)$, where $m:A\rightarrow X$ belongs to $\mathcal{M}$.
\end{definition}

\noindent The class $\mathcal{M}$ plays the role of a class of ``inclusion''
maps from ``smaller'', or ``boundary'', objects into ``larger'' objects.

In Definition \ref{Def25} we will be formulating properties for an
appropriate class $\mathcal{M}$, but first we need the following concepts:

\begin{definition}
Let $(\mathbf{C},F,P\mathbf{)}$ be a topological starting category.

\begin{enumerate}
\item[\emph{(i)}]  An object $A$ of $\mathbf{C}$ is said to be \emph{%
irreducible}, if $F(A)$ is non-empty and connected.

\item[\emph{(ii)}]  An object of $\mathbf{C}$ is said to be \emph{factorized}%
, if it is of the form $P^{n}(E)$, with $n\in \Bbb{N}$ (where $P^{0}:=%
\limfunc{Id}$) or a finite monoidal product of irreducible objects.\smallskip
\end{enumerate}
\end{definition}

\begin{definition}
\label{Def25}We say that $\mathcal{M}$ of Definition \ref{Def23} is an \emph{%
appropriate} class if it satisfies the following properties:
\end{definition}

\begin{enumerate}
\item  For any $m\in \mathcal{M}$, $F(m)$ is an embedding in $\mathrm{Top}$,
i.e. a homeomorphism onto its image.

\item  \label{monoidclosed}$\mathcal{M}$ is closed under the monoidal
product.

\item  $\mathcal{M}$ is closed under the endofunctor $P$.

\item  \label{prop4M}The domain of any $m\in \mathcal{M}$ is either $%
P^{n}(E) $ with $n\in \Bbb{N}$ or any object that can be generated from
irreducible objects by taking finite monoidal products and applying $P$.

\item  \label{permiso}(isomorphism-closure property)\newline
If $m\in \mathcal{M}$ and $f\in \limfunc{Iso}(\mathbf{C})$\footnote{%
The class $\limfunc{Iso}(\mathbf{C})$ is the class of all isomorphisms of $%
\mathbf{C}$.} such that $m\circ f$ exists, then $m\circ f\in \mathcal{M}$.

\item  \label{factprop}(subdivision property)\newline
Let $(A_{i})_{i\in I}$ be a family of irreducible objects, where $I$ is a
finite ordered index set. We denote by $A_{I}:=\sqcup _{i\in I}A_{i}$ the
corresponding factorized object and set $A_{\emptyset }:=E$.$\smallskip $%
\newline
a)\label{factpropa} For any $\mathcal{M}$-subobject $(A_{I},m)$ of $X$ and
any proper ordered subset $J$ of $I$ there exists a morphism $%
s_{J,I}:A_{J}\rightarrow A_{I}$ (independent of $m$) such that $%
m_{J}:=m\circ s_{J,I}$ belongs to $\mathcal{M}$, i.e. the following diagram%
\begin{center}
\settriparms[1`1`-1;600]
\Vtriangle[A_{J}`X`A_{I};m_{J}`s_{J,I}`m]
\end{center}
%
commutes.\newline
Under the forgetful functor the morphisms $F(s_{J,I}):F(A)_{J}\rightarrow
F(A)_{I}$ are the canonical monomorphisms.$\smallskip $\newline
b) For $K$ a proper ordered subset of $J$ one has:
\[
s_{J,I}\circ s_{K,J}=s_{K,I}\text{.}
\]

c) For any finite ordered sets $J\varsubsetneq I$, $L\varsubsetneq K$, with $%
I$ and $K$ disjoint, one has:
\[
s_{J,I}\sqcup s_{L,K}=s_{JL,IK}\text{,}
\]
where $IK$ denotes the finite ordered set consisting of the elements of $I$
followed by the elements of $K$.$\smallskip $\newline
d)\label{factpropc} For any $\mathcal{M}$-subobjects $(A_{I},m)$ and $%
(A_{I^{\prime }}^{\prime },m^{\prime })$ with $\func{card}(I^{\prime })=%
\func{card}(I)$, isomorphism $\alpha :A_{I}\rightarrow A_{I^{\prime
}}^{\prime }$, and a proper ordered subset $J^{\prime }$ of $I^{\prime }$,
there exists a proper ordered subset $J$ of $I$ with $\func{card}(J)=\func{%
card}(J^{\prime })$ and an isomorphism $\alpha _{J,J^{\prime }}:=\sqcup
_{i\in J}\alpha _{i}:A_{J}\rightarrow A_{J^{\prime }}^{\prime }$
(independent of $m$ and $m^{\prime }$), such that the following diagram%
\begin{center}
\setsqparms[1`-1`-1`1;800`700]
\square[A_{I}`A'_{I'}`A_{J}`A'_{J'};\alpha`s_{J,I}`s_{J',I'}`\alpha _{J,J'}]
\end{center}
%
commutes.
\end{enumerate}

\begin{thm_remark}
\label{mnotation}The intuitive content of the subdivision property is that
when a collection of objects form a subobject of $X$ any subcollection of
them also do, with the respective elements of $\mathcal{M}$ related in the
appropriate way.\newline
The condition d) says that an isomorphism $\alpha $ induces a one-to-one
correspondence between connected components of the domain and codomain
subobjects, and isomorphisms $\alpha _{i}$ between corresponding components.
Since $\alpha _{J,J^{\prime }}$ is unique by the commutativity of the
diagram in $\mathcal{M}$-\ref{factprop}d) and the fact that $s_{J^{\prime
},I^{\prime }}$ is a monomorphism (which in its turn follows from $\mathcal{M%
}$-\ref{factprop}a)), for $K\varsubsetneq J$ and $K^{\prime }\varsubsetneq
J^{\prime }$ we have, from $\mathcal{M}$-\ref{factprop}b), $(\alpha
_{J,J^{\prime }})_{K,K^{\prime }}=\alpha _{K,K^{\prime }}$ and we use the
simpler expression on the right hand side when it occurs.\smallskip
\end{thm_remark}

{\noindent \bf Example: }%
We choose $\mathcal{M}$ to be the class of monomorphisms whose domain is
either $\emptyset $ or a finite disjoint union of oriented circles and whose
codomain is either $\emptyset $ or an oriented surface with boundary, such
that each circle in the domain is mapped isomorphically to a boundary
component of the codomain surface. Examples of morphisms in $\mathcal{M}$
are the empty maps from $\emptyset $ to $\emptyset $, to $T_{\pm }$ or to $%
S_{\pm }$; the maps $C_{\pm }\rightarrow D_{\pm }$ sending $z$ to $z$, and
the map $C_{-+}\rightarrow A_{+}$, sending $(z,1)$ to $z$ and $(z,2)$ to $2z$%
. The class $\mathcal{M}$ clearly satisfies the properties 1)-6) above. The
morphisms $s_{J,I}$ in the subdivision property are the canonical
monomorphisms, e.g. for any $m$ with domain $C_{+-+-}$ and $I=\left\{
1,2,3,4\right\} $, $J=\left\{ 2,3\right\} $ one has $s_{J,I}:C_{-+}%
\rightarrow C_{+-+-}$ given by $(z,1)\mapsto (z,2)$ and $(z,2)\mapsto (z,3)$.%
$\blacktriangle $\medskip

\begin{definition}
\label{DefN27}Given a topological starting category $(\mathbf{C,}F,P)$ and
an appropriate class of monomorphisms $\mathcal{M}$, we define a symmetric,
strict monoidal category $(\mathcal{S}(\mathbf{C}),\sqcup ,E,c)$, also
written $\mathcal{S}(\mathbf{C})$ for short, as follows:

\begin{enumerate}
\item[\emph{a)}]  the class of objects of $\mathcal{S}(\mathbf{C})$ is the
subclass of objects of $\mathbf{C}$ which are the domain of some $m\in
\mathcal{M}$,

\item[\emph{b)}]  the morphisms of $\mathcal{S}(\mathbf{C})$ are the
isomorphisms of $\mathbf{C}$ restricted to the objects of $\mathcal{S}(%
\mathbf{C})$,

\item[\emph{c)}]  the composition, identity morphisms, monoidal product $%
\sqcup $, unit object $E$ and braiding $c$ are inherited from $\mathbf{C}$.
\end{enumerate}

We define a symmetric, strong monoidal endofunctor $P$ on $\mathcal{S}(%
\mathbf{C})$ as the restriction to $\mathcal{S}(\mathbf{C})$ of $P$ defined
on $\mathbf{C}$. The \emph{category of subobjects}, is the pair $(\mathcal{S}%
(\mathbf{C}),P)$, where $\mathcal{S}(\mathbf{C})$ denotes the category with
its symmetric, monoidal structure and $P$ denotes the above monoidal
endofunctor on $\mathcal{S}(\mathbf{C})$.
\end{definition}

\begin{theorem}
\label{Thm26}$(\mathcal{S}(\mathbf{C}),\sqcup ,E,c)$ is indeed a symmetric,
strict monoidal category, and $P$ restricts to a symmetric, strong monoidal
endofunctor on $\mathcal{S}(\mathbf{C})$.
\end{theorem}

\proof%
One just has to check that the composition, identity morphisms, monoidal
structure and braiding indeed restrict to $\mathcal{S}(\mathbf{C})$. For the
composition and identity morphisms this is clear. The monoidal product is
defined in $\mathcal{S}(\mathbf{C})$ because of property \ref{monoidclosed}
of $\mathcal{M}$. The unit object $E$ of $\mathbf{C}$ is an object of $%
\mathcal{S}(\mathbf{C})$, due to property $\mathcal{M}$-6a), since $\mathcal{%
M}$ is non-empty. The braiding isomorphisms restricted to the objects of $%
\mathcal{S}(\mathbf{C})$ are morphisms of $\mathcal{S}(\mathbf{C})$ because
of property \ref{permiso} (isomorphism-closure property) of $\mathcal{M}$.
The functor $P$ closes on the objects of $\mathcal{S}(\mathbf{C})$ because
of $\mathcal{M}$-3), and on the morphisms of $\mathcal{S}(\mathbf{C})$
because, being a functor, $P$ maps isomorphisms to isomorphisms. Thus $P$
defines an endofunctor on $\mathcal{S}(\mathbf{C})$. For any $A,B\in \func{Ob%
}(\mathcal{S}(\mathbf{C}))$, $P(A)\sqcup P(B)$ and $P(A\sqcup B)$ are
objects of $\mathcal{S}(\mathbf{C})$ by $\mathcal{M}$-2) and $\mathcal{M}$%
-3), and thus $\pi _{2(A,B)}$, which is an isomorphism since $P$ is strong
monoidal, is a morphism of $\mathcal{S}(\mathbf{C})$. Likewise $\pi _{0}$ is
a morphism of $\mathcal{S}(\mathbf{C})$ and, looking at Definition (\ref
{Def1.3}), it is clear that $(P,\pi _{2},\pi _{0})$ is strong monoidal on $%
\mathcal{S}(\mathbf{C})$. It is also clearly symmetric.%
\endproof%
\medskip

\begin{thm_remark}
Using the isomorphism-closure property of $\mathcal{M}$, any object of $%
\mathcal{S}(\mathbf{C})$ is isomorphic to a factorized object of $\mathcal{S}%
(\mathbf{C})$, via a string of morphisms which are monoidal products of
identity morphisms and morphisms of the form $P^{k}(\pi _{2(X,X^{\prime
})}^{-1})$, with $k\in \Bbb{N}$. We will proceed from now on to write all
objects of $\mathcal{S}(\mathbf{C})$ as if they were factorized and let the
necessary isomorphisms to achieve this be understood.
\end{thm_remark}

We will distinguish two special cases of morphisms of $\mathcal{S}(\mathbf{C}%
)$.

\begin{definition}
A morphism of $\mathcal{S}(\mathbf{C})$ is called:

\begin{enumerate}
\item  a \emph{permuting isomorphism}, if it is of the form
\[
p_{\sigma }:A_{1}\sqcup \cdots \sqcup A_{n}\rightarrow A_{\sigma (1)}\sqcup
\cdots \sqcup A_{\sigma (n)}\text{,}
\]
where $A_{i}$ are irreducible subobjects and $\sigma $ is a permutation of $%
\left\{ 1,2,\ldots ,n\right\} \subseteq \Bbb{N}$, and it is generated, via
the monoidal product and composition, by identity morphisms and braidings
only. Because of naturality, a permuting isomorphism may always be written
as a composition of elementary permuting isomorphisms of the form:
\[
(\tbigsqcup\limits_{j=1}^{i-1}\limfunc{id}\nolimits_{A_{j}})\sqcup
c_{(A_{i},A_{i+1})}\sqcup (\tbigsqcup\limits_{j=i+2}^{n}\limfunc{id}%
\nolimits_{A_{j}})\text{.}
\]

\item  an \emph{order-preserving isomorphism}, if it is of the form
\[
\alpha _{1}\sqcup \cdots \sqcup \alpha _{n}:A_{1}\sqcup A_{2}\sqcup \cdots
\sqcup A_{n}\rightarrow A_{1}^{\prime }\sqcup A_{2}^{\prime }\sqcup \cdots
\sqcup A_{n}^{\prime }\text{,}
\]
where $A_{i},A_{i}^{\prime }$ are irreducible subobjects and each $\alpha
_{i}:A_{i}\rightarrow A_{i}^{\prime }$ is an isomorphism of $\mathbf{C}$%
.\smallskip
\end{enumerate}
\end{definition}

{\noindent \bf Example: }%
The category $\mathcal{S}(\mathbf{C})$ has as objects the empty space $%
\emptyset $, the circles $C_{\pm }$ and finite disjoint unions of these
circles. The morphisms of $\mathcal{S}(\mathbf{C})$, apart from the identity
morphism for $\emptyset $, are orientation-preserving automorphisms of
circles and finite disjoint unions of these morphisms, as well as
isomorphisms such as $\alpha :C_{-+}\rightarrow C_{++}$, given by $%
(z,1)\mapsto (\overline{z},2)$, $(z,2)\mapsto (z,1)$, which is the
composition of a braiding (i.e. a permuting isomorphism) and an
order-preserving isomorphism.$\blacktriangle $\medskip

We are now ready to construct the topological category $\mathcal{C}$, which
will be the domain category for the TQFT functor.

\begin{definition}
\label{Def213}Let $\mathcal{O}$ be a class of triples $(X,A,m)$, where $%
(A,m) $ is an $\mathcal{M}$-subobject of $X$, satisfying:

\begin{enumerate}
\item[1)]  $(E,E,\mathrm{id}_{E})$ $\in \mathcal{O}$.

\item[2)]  $(X,A,m)$, $(X^{\prime },A^{\prime },m^{\prime })\in \mathcal{O}$
implies $(X\sqcup X^{\prime },A\sqcup A^{\prime },m\sqcup m^{\prime })\in
\mathcal{O}$.

\item[3)]  $(X,A,m)\in \mathcal{O}$ implies $(P(X),P(A),P(m))\in \mathcal{O}$%
.

\item[4)]  $(X\sqcup X^{\prime },A\sqcup A^{\prime },m\sqcup m^{\prime })\in
\mathcal{O}$ implies $(X^{\prime }\sqcup X,A^{\prime }\sqcup A,m^{\prime
}\sqcup m)\in \mathcal{O}$.
\end{enumerate}

We define $\mathcal{C}$ to be the category with $\mathrm{Ob}(\mathcal{C})=%
\mathcal{O}$, morphisms consisting of two classes, isomorphisms and gluing
morphisms, defined below (Definition \ref{Defisomorphisms} and Definition
\ref{DefGluings}), and composition and identity morphisms defined below
(Definition \ref{def27} and \ref{Defidentities}).\medskip
\end{definition}

{\noindent \bf Example: }%
For our example we choose the class of objects of $\mathcal{C}$ to be all
triples $(X,A,m)$, where $A$ is $\emptyset $ or a finite disjoint union of
oriented circles isomorphic to $\partial X$ and the image of $m\in \mathcal{M%
}$ is $\partial X$. Thus for example we have the following objects:

\begin{itemize}
\item  $X=A=\emptyset $ and $m=\limfunc{id}\nolimits_{\emptyset }$,

\item  $X=S_{\pm }$ or $T_{\pm }$ (the sphere or the torus, with either
orientation), $A=\emptyset $, and $m$ is the corresponding empty map,

\item  $X=A_{+}$ (the annulus with positive orientation), $A=C_{-+}$ and $m$
is given by $m(z,1)=z$, $m(z,2)=2z$.
\end{itemize}

\noindent Note that, in this example, we do not admit objects having, e.g. $%
X=A_{+}$ and $A=C_{+}$, where $A$ only corresponds to part of the boundary
of $X$, although the general theory does not forbid this.$\blacktriangle $%
\medskip

The morphisms of $\mathcal{C}$ are restricted to be of two types,
isomorphisms and gluing morphisms. We start by defining the isomorphisms:

\begin{definition}
\label{Defisomorphisms}Let $(X,A,m)$ and $(X^{\prime },A^{\prime },m^{\prime
})$ be two objects of $\mathcal{C}$. An isomorphism between them is a pair $%
(f,\alpha )$, where $f:X\rightarrow X^{\prime }$ is an isomorphism in $%
\mathbf{C}$ and $\alpha :A\rightarrow A^{\prime }$ is a morphism of $%
\mathcal{S}(\mathbf{C})$, such that the following diagram%
\begin{center}
\setsqparms[1`-1`-1`1;800`600]
\square[X`X'`A`A';f`m`m'`\alpha]
\end{center}
%
is commutative.\medskip
\end{definition}

\begin{thm_remark}
Two general types of isomorphisms are:
\end{thm_remark}

\begin{itemize}
\item  $(\limfunc{id}\nolimits_{X},\limfunc{id}\nolimits_{A})$ for any $%
(X,A,m)\in \func{Ob}(\mathcal{C})$.

\item  $(c_{(X,X^{\prime })},c_{(A,A^{\prime })}):(X\sqcup X^{\prime
},A\sqcup A^{\prime },m\sqcup m^{\prime })\rightarrow (X^{\prime }\sqcup
X,A^{\prime }\sqcup A,m^{\prime }\sqcup m)$ for any $(X,A,m),(X^{\prime
},A^{\prime },m^{\prime })\in \func{Ob}(\mathcal{C})$.
\end{itemize}

{\noindent \bf Example: }%
\label{examplepg8} Some examples of isomorphisms of $\mathcal{C}$, apart
from the ones in the above remark, are:

\begin{itemize}
\item  the reverse map
\[
(r,\alpha ):(D_{+},C_{+},m_{+})\rightarrow (D_{-},C_{-},m_{-})\text{,}
\]
where $m_{\pm }:C_{\pm }\rightarrow D_{\pm }$ are given by $m_{\pm }(z)=z$, $%
r:D_{+}\rightarrow D_{-}$ is given by $r(z)=\overline{z}$ and $\alpha
:C_{+}\rightarrow C_{-}$ is given by $\alpha (z)=\overline{z}$.

\item  the following morphism, which can be interpreted as a reordering of
the two boundary components of the annulus:
\[
(\limfunc{id}\nolimits_{A_{+}},c_{(C_{-},C_{+})}):(A_{+},C_{-+},m)%
\rightarrow (A_{+},C_{+-},m^{\prime })\text{,}
\]
where $m(z,1)=z=m^{\prime }(z,2)$ and $m(z,2)=2z=m^{\prime }(z,1)$.$%
\blacktriangle $\medskip
\end{itemize}

Next we define the other class of morphisms of $\mathcal{C}$, which we call
gluing morphisms. We recall the notation $A_{I}:=\sqcup _{i\in I}A_{i}$ for
a finite ordered index set, from property \ref{factprop} of $\mathcal{M}$.
Also we set $P(A)_{I}:=\sqcup _{i\in I}P(A_{i})$.\medskip

\begin{definition}
\label{DefGluings}Let $(X,A,m)$, $(X^{\prime },A^{\prime },m^{\prime })$ be
objects of $\mathcal{C}$, where $A:=A_{N}$ for some finite ordered index set
$N$. Let $I$, $J$ and $R$ be disjoint ordered subsets of $N$, with $I$, $J$
non-empty and of the same cardinality, and such that, as sets, $N=I\cup
J\cup R$. A gluing morphism from $(X,A,m)$ to $(X^{\prime },A^{\prime
},m^{\prime })$ is a triple $(f,\varphi ,\alpha )$, where $X\stackrel{f}{%
\rightarrow }X^{\prime }$ is a morphism of $\mathbf{C}$, and $\varphi
:A_{I}\rightarrow P(A)_{J}$ and $\alpha :A_{R}\longrightarrow A^{\prime }$
are morphisms of $\mathcal{S}(\mathbf{C)}$, such that:

\begin{enumerate}
\item[\emph{i)}]  $f\circ m_{R}=m^{\prime }\circ \alpha $,

\item[\emph{ii)}]  $F(f\circ m_{I})$ is a (topological) embedding, the
images of $F(f\circ m_{I})$ and $F(f\circ m_{R})$ are separated in $%
F(X^{\prime })$, and the image under $F(f)$ of $\func{Cls}(\func{Im}%
(F(m_{I})))$ is closed in $F(X^{\prime })$, where $\func{Cls}$ is the
closure.

\item[\emph{iii)}]  $F(f\circ m_{J})\circ F(\varphi )=F(f\circ m_{I})$,

\item[\emph{iv)}]  $F(X)\setminus \func{Im}(F(m_{IJ}))$ is isomorphic to $%
F(X^{\prime })\setminus \func{Im}(F(f\circ m_{I}))$, via $F(f)$.
\end{enumerate}
\end{definition}

\begin{thm_remark}
Let us refer to the connected components of $A$, $A^{\prime }$, etc. as
boundaries, although they are not necessarily boundary components of
topological manifolds. Intuitively we are gluing the boundaries $A_{I}$ to
the boundaries $A_{J}$ and the remaining boundaries $A_{R}$ are not glued.
The gluing morphism itself is really the morphism $f$ from $X$ ``before
gluing'' to $X^{\prime }$, which is a copy of the quotient space ``after
gluing'' (see figure below). Note that the gluing can involve gluing
together disconnected spaces as well as self-gluing of a connected space.
The morphisms $\varphi $ and $\alpha $ give supplementary information about
how the boundaries are glued together, and how the boundaries of $X$ which
are not glued map to the boundary of $X^{\prime }$.\textbf{\ }We note that%
\textbf{\ }different choices of $I$, $J$ and $\varphi $ may be associated
with the same morphisms $f$ and $\alpha $, since we may reorder $A_{I}$ and $%
A_{J}$ using permuting isomorphisms and compose $\varphi $ with the
respective inverse isomorphisms to make a different $\varphi $. Condition i)
of the definition says that the boundaries of $X$ which are not glued are
isomorphic to the boundaries of $X^{\prime }$, via $\alpha $. Conditions
ii)-iv) are purely topological and capture the intuitive idea of gluing.
Conditions ii) and iii) say that the glued boundaries map $2:1$ onto their
image in $X^{\prime }$, which is isomorphic to each half of the glued
boundaries. This image has to be separated from the boundaries of $X^{\prime
}$ to avoid problems when composing two gluing morphisms. Note that $%
F(\varphi )$ in iii) is defined in the obvious way by regarding $\varphi $
as a morphism of $\mathbf{C}$. (Indeed the forgetful functor $F$ on $\mathbf{%
C}$ restricts to the subcategory $\mathcal{S}(\mathbf{C})$.) Finally,
condition iv) says that, topologically, everything that is not glued in $X$
is isomorphic to everything in $X^{\prime }$ that is not the image of the
glued boundaries.

\begin{figure}[h]
\centerline{\psfig{figure=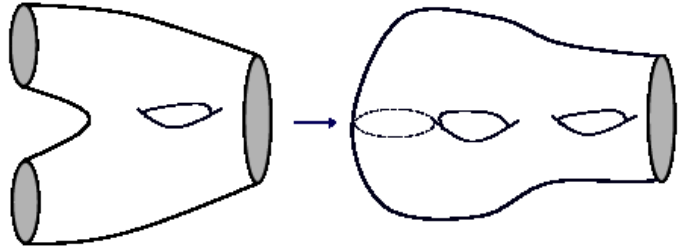,height=3.5cm,
}}
\caption{Self-gluing of an object}
\label{selfglue}
\end{figure}

\vspace{0.3in}
\end{thm_remark}

{\noindent \bf Example: }%
\label{examplepg10} Gluing morphisms can arise from the gluing together of
separate objects, or the self-gluing of a single object, as well as a
combination of both types of gluing.

\begin{enumerate}
\item[a)]  The gluing together of two discs to give a sphere $S$ may be
described by the gluing morphism:
\[
(f,\varphi ,\alpha ):(D_{-+},C_{-+},m)\rightarrow (S_{+},\emptyset
,m^{\prime })\text{,}
\]
where $m(z,1)=(z,1)$, $m(z,2)=(z,2)$, $m^{\prime }$ is the empty map, and
\[
I=\left\{ 1\right\} \text{, }J=\left\{ 2\right\} \text{, }R=\emptyset \text{.%
}
\]
The morphisms $f$, $\varphi $ and $\alpha $ are given by (writing $z=x+iy\in
D$):
\begin{eqnarray*}
f(z,1) &=&(x,y,-\sqrt{1-\left| z\right| ^{2}})\text{\hspace{0.15in},\hspace{%
0.15in}}f(z,2)=(x,y,\sqrt{1-\left| z\right| ^{2}}) \\
\varphi &=&\limfunc{id}\nolimits_{C_{-}}\text{\hspace{0.15in},\hspace{0.15in}%
}\alpha =\limfunc{id}\nolimits_{\emptyset }\text{.}
\end{eqnarray*}

\item[b)]  A similar example is the gluing of a disc into one of the holes
of a pair-of-pants $P$ to give an annulus $A$. For convenience we take $P$
to be
\[
P=A\setminus \left\{ z:\left| z+\tfrac{3}{2}\right| <\tfrac{1}{4}\right\}
\text{.}
\]
The gluing morphism is then:
\[
(f,\varphi ,\alpha ):(D_{+}\sqcup P_{+},C_{+--+},m)\rightarrow
(A_{+},C_{-+},m^{\prime })\text{,}
\]
where $m(z,1)=(z,1)$, $m(z,2)=(-\frac{3}{2}+\frac{z}{4},2)$, $m(z,3)=(z,2)$,
$m(z,4)=(2z,2)$, $m^{\prime }(z,1)=z$, $m^{\prime }(z,2)=2z$ and
\[
I=\left\{ 1\right\} \hspace{0.1in}\text{,}\hspace{0.1in}J=\left\{ 2\right\}
\hspace{0.1in}\text{,}\hspace{0.1in}R=\left\{ 3,4\right\} \text{.}
\]
The morphisms $f$, $\varphi $ and $\alpha $ are given by:
\[
f(z,1)=-\frac{3}{2}+\frac{z}{4}\hspace{0.1in}\text{,}\hspace{0.1in}f(z,2)=z%
\hspace{0.1in}\text{,}\hspace{0.1in}\varphi =\limfunc{id}\nolimits_{C_{+}}%
\hspace{0.1in}\text{,}\hspace{0.1in}\alpha =\alpha _{1}\sqcup \alpha _{2}=%
\limfunc{id}\nolimits_{C_{-+}}\text{.}
\]

\item[c)]  An example of a self-gluing is the gluing morphism from the
annulus to the torus:
\[
(g,\psi ,\beta ):(A_{+},C_{-+},m^{\prime })\rightarrow (T_{+},\emptyset
,m^{\prime \prime })\text{,}
\]
where $m^{\prime }$ is as in b) and $m^{\prime \prime }$ is the empty map,
and
\[
I=\left\{ 1\right\} \hspace{0.1in}\text{,}\hspace{0.1in}J=\left\{ 2\right\}
\hspace{0.1in}\text{,}\hspace{0.1in}R=\emptyset \text{.}
\]
The morphisms $g$, $\psi $ and $\beta $ are given by:
\[
g(z)=(\frac{z}{\left| z\right| },e^{2\pi i(2-\left| z\right| )})\hspace{0.1in%
}\text{,}\hspace{0.1in}\psi =\limfunc{id}\nolimits_{C_{-}}\hspace{0.1in}%
\text{,}\hspace{0.1in}\beta =\limfunc{id}\nolimits_{\emptyset }\text{.\label%
{examplepg17}}
\]

\item[d)]  Finally we can glue a disc to a pair-of-pants simultaneously with
the self-gluing of its remaining boundary components via the following
gluing morphism:
\[
(h,\chi ,\gamma ):(D_{+}\sqcup P_{+},C_{+--+},m)\rightarrow (T_{+},\emptyset
,m^{\prime \prime })\text{,}
\]
where $m$, $m^{\prime \prime }$ are as in b) and c), and
\[
I=\left\{ 1,3\right\} \hspace{0.1in}\text{,}\hspace{0.1in}J=\left\{
2,4\right\} \hspace{0.1in}\text{,}\hspace{0.1in}R=\emptyset \text{.}
\]
The morphisms $h$, $\chi $ and $\gamma $ are given by:
\begin{eqnarray*}
h(z,1) &=&(\frac{-\frac{3}{2}+\frac{z}{4}}{\left| -\frac{3}{2}+\frac{z}{4}%
\right| },e^{2\pi i(2-\left| -\frac{3}{2}+\frac{z}{4}\right| )})\hspace{0.1in%
}\text{,}\hspace{0.1in}h(z,2)=(\frac{z}{\left| z\right| },e^{2\pi i(2-\left|
z\right| )})\text{,} \\
&&\chi :=\chi _{1}\sqcup \chi _{2}=\limfunc{id}\nolimits_{C_{+-}}\hspace{%
0.1in}\text{,}\hspace{0.1in}\gamma =\limfunc{id}\nolimits_{\emptyset }\text{.%
}\blacktriangle
\end{eqnarray*}
\medskip
\end{enumerate}

Next we define the composition of morphisms in $\mathcal{C}$:

\begin{definition}
\label{def27}Let $(X,A,m)$, $(X^{\prime },A^{\prime },m^{\prime })$ and $%
(X^{\prime \prime },A^{\prime \prime },m^{\prime \prime })$ be objects of $%
\mathcal{C}$.

\begin{itemize}
\item  (isomorphism-isomorphism)

Let $(f,\alpha ):(X,A,m)\rightarrow (X^{\prime },A^{\prime },m^{\prime })$
and $(g,\beta ):(X^{\prime },A^{\prime },m^{\prime })\rightarrow (X^{\prime
\prime },A^{\prime \prime },m^{\prime \prime })$ be isomorphisms. Then we
define:
\[
(g,\beta )\circ (f,\alpha )=(g\circ f,\beta \circ \alpha )\text{.}
\]

\item  (gluing - isomorphism)

Let $(f,\varphi ,\alpha ):(X,A,m)\rightarrow (X^{\prime },A^{\prime
},m^{\prime })$ be a gluing morphism with $\varphi :A_{I}\rightarrow
P(A)_{J} $, $\alpha :A_{R}\rightarrow A^{\prime }$ and let $(g,\beta
):(X^{\prime },A^{\prime },m^{\prime })\rightarrow (X^{\prime \prime
},A^{\prime \prime },m^{\prime \prime })$ be an isomorphism. Then we define:
\[
(g,\beta )\circ (f,\varphi ,\alpha )=(g\circ f,\varphi ,\beta \circ \alpha )%
\text{.}
\]

\item  (isomorphism - gluing)

Let $(f,\alpha ):(X,A,m)\rightarrow (X^{\prime },A^{\prime },m^{\prime })$
be an isomorphism and let $(g,\psi ,\beta ):(X^{\prime },A^{\prime
},m^{\prime })\rightarrow (X^{\prime \prime },A^{\prime \prime },m^{\prime
\prime })$ be a gluing morphism with $\psi :A_{I^{\prime }}^{\prime
}\rightarrow P(A^{\prime })_{J^{\prime }}$ and $\beta :A_{R^{\prime
}}^{\prime }\rightarrow A^{\prime \prime }$. By property $\mathcal{M}$-\ref
{factprop}d), corresponding to $I^{\prime }$, $J^{\prime }$ and $R^{\prime }$%
, there exist ordered sets $I$, $J$, $R$ and isomorphisms
\begin{eqnarray*}
&&\alpha _{I,I^{\prime }}:A_{I}\rightarrow A_{I^{\prime }}^{\prime }\text{,}
\\
&&\alpha _{J,J^{\prime }}:A_{J}\rightarrow A_{J^{\prime }}^{\prime }\text{,}
\\
&&\alpha _{R,R^{\prime }}:A_{R}\rightarrow A_{R^{\prime }}^{\prime }\text{.}
\end{eqnarray*}
Then we define:
\[
(g,\psi ,\beta )\circ (f,\alpha )=(g\circ f,P(\alpha ^{-1})_{J,J^{\prime
}}\circ \psi \circ \alpha _{I,I^{\prime }},\beta \circ \alpha _{R,R^{\prime
}})\text{,}
\]
where $P(\alpha ^{-1})_{J,J^{\prime }}$ is defined to be the isomorphism $%
\sqcup _{i\in J}P(\alpha _{i}^{-1})$, and the $\alpha _{i}$'s are the factor
isomorphisms of property $\mathcal{M}$-\ref{factprop}d).

\item  (gluing - gluing)

Let $(f,\varphi ,\alpha ):(X,A,m)\rightarrow (X^{\prime },A^{\prime
},m^{\prime })$ and $(g,\psi ,\beta ):(X^{\prime },A^{\prime },m^{\prime
})\rightarrow (X^{\prime \prime },A^{\prime \prime },m^{\prime \prime })$ be
gluing morphisms with
\begin{eqnarray*}
&&\varphi :A_{I}\rightarrow P(A)_{J}\text{,} \\
&&\alpha :A_{R}\rightarrow A^{\prime }\text{,} \\
&&\psi :A_{I^{\prime }}^{\prime }\rightarrow P(A^{\prime })_{J^{\prime }}%
\text{,} \\
&&\beta :A_{R^{\prime }}^{\prime }\rightarrow A^{\prime \prime }\text{,}
\end{eqnarray*}
By property $\mathcal{M}$-\ref{factprop}d), corresponding to $I^{\prime }$, $%
J^{\prime }$ and $R^{\prime }$, there exist $\widetilde{I}$, $\widetilde{J}$
and $\widetilde{R}$, ordered subsets of $R$, and isomorphisms:
\[
\alpha _{\widetilde{I},I^{\prime }}:A_{\widetilde{I}}\rightarrow
A_{I^{\prime }}^{\prime }\text{\hspace{0.15in},\hspace{0.15in}}\alpha _{%
\widetilde{J},J^{\prime }}:A_{\widetilde{J}}\rightarrow A_{J^{\prime
}}^{\prime }\text{\hspace{0.15in}and\hspace{0.15in}}\alpha _{\widetilde{R}%
,R^{\prime }}:A_{\widetilde{R}}\rightarrow A_{R^{\prime }}^{\prime }\text{.}
\]
Then we define:
\[
(g,\psi ,\beta )\circ (f,\varphi ,\alpha )=(g\circ f,\varphi \sqcup
(P(\alpha ^{-1})_{\widetilde{J},J^{\prime }}\circ \psi \circ \alpha _{%
\widetilde{I},I^{\prime }}),\beta \circ \alpha _{\widetilde{R},R^{\prime }})%
\text{.}
\]
\smallskip
\end{itemize}
\end{definition}

{\noindent \bf Example: }%
The composition
\[
(g,\psi ,\beta )\circ (f,\varphi ,\alpha )=(g\circ f,\varphi \sqcup P(\alpha
_{2}^{-1})\circ \psi \circ \alpha _{1},\limfunc{id}\nolimits_{\emptyset })
\]
from $(D_{+}\sqcup P_{+},C_{+--+},m)$ to $(T_{+},\emptyset ,m^{\prime \prime
})$ from our previous example b) and c) is equal to $(h,\chi ,\gamma )$ from
the same example d), as is easily checked.$\blacktriangle $\medskip

The main result concerning the morphisms of $\mathcal{C}$ is:

\begin{theorem}
The class of morphisms of $\mathcal{C}$ is closed under the above
composition.
\end{theorem}

\proof%
We will just show this for the gluing-gluing case of Definition \ref{def27}.
The other combinations are proved analogously. Let $(f,\varphi ,\alpha
):(X,A,m)\rightarrow (X^{\prime },A^{\prime },m^{\prime })$ and $(g,\psi
,\beta ):(X^{\prime },A^{\prime },m^{\prime })\rightarrow (X^{\prime \prime
},A^{\prime \prime },m^{\prime \prime })$ be gluing morphisms, where $%
\varphi :A_{I}\rightarrow P(A)_{J}$, $\alpha :A_{R}\rightarrow A^{\prime }$,
$\psi :A_{I^{\prime }}^{\prime }\rightarrow P(A^{\prime })_{J^{\prime }}$
and $\beta :A_{R^{\prime }}^{\prime }\rightarrow A^{\prime \prime }$%
.\smallskip \newline
i) Set $A=A_{N}$ and $A^{\prime }=A_{N^{\prime }}^{\prime }$, where $N$, $%
N^{\prime }$ are ordered index sets. By property $\mathcal{M}$-\ref
{factpropc}d) there exist an ordered set $\widetilde{R}\varsubsetneq R$ and
an isomorphism $\alpha _{\widetilde{R},R^{\prime }}$ such that $\alpha \circ
s_{\widetilde{R},R}=s_{R^{\prime },N^{\prime }}\circ \alpha _{\widetilde{R}%
,R^{\prime }}$. Thus
\begin{eqnarray*}
f\circ m_{\widetilde{R}} &=&f\circ m_{R}\circ s_{\widetilde{R},R} \\
&=&m^{\prime }\circ \alpha \circ s_{\widetilde{R},R} \\
&=&m^{\prime }\circ s_{R^{\prime },N^{\prime }}\circ \alpha _{\widetilde{R}%
,R^{\prime }} \\
&=&m_{R^{\prime }}^{\prime }\circ \alpha _{\widetilde{R},R^{\prime }}\text{,}
\end{eqnarray*}
where we use i) for the gluing morphism $(f,\varphi ,\alpha )$ in the second
equality. Here and below we use $\mathcal{M}$-\ref{factpropc}a) and the
simplified notation in Remark \ref{mnotation} without comment. Composing
this equation with $g$ and using i) for the gluing morphism $(g,\psi ,\beta
) $ we have:
\[
(g\circ f)\circ m_{\widetilde{R}}=g\circ m_{R^{\prime }}^{\prime }\circ
\alpha _{\widetilde{R},R^{\prime }}=m^{\prime \prime }\circ (\beta \circ
\alpha _{\widetilde{R},R^{\prime }})\text{.}
\]
\newline
To simplify the notation for the proof of the topological conditions, we
will omit writing the forgetful functor $F$ and consider all objects and
morphisms to be in $\mathrm{Top}$.\smallskip \newline
ii) First we show that $(g\circ f)\circ m_{I}:A_{I}\rightarrow X^{\prime
\prime }$ and $(g\circ f)\circ m_{\widetilde{I}}:A_{\widetilde{I}%
}\rightarrow X^{\prime \prime }$ are topological embeddings and then we show
that $\func{Im}((g\circ f)\circ m_{I})$ and $\func{Im}((g\circ f)\circ m_{%
\widetilde{I}})$ are separated in $X^{\prime \prime }$, thus concluding that
$(g\circ f)\circ m_{I\widetilde{I}}$ is a topological embedding.\newline
It is straightforward to show $\func{Im}(m_{I\widetilde{I}})=\func{Im}%
(m_{I})\cup \func{Im}(m_{\widetilde{I}})$, using the commutativities $m_{I%
\widetilde{I}}\circ s_{I,I\widetilde{I}}=m_{I}$ and $m_{I\widetilde{I}}\circ
s_{\widetilde{I},I\widetilde{I}}=m_{\widetilde{I}}$ and bearing in mind that
$s_{I,I\widetilde{I}}$ and $s_{\widetilde{I},I\widetilde{I}}$ are the
canonical monomorphisms in $\mathrm{Top}$. Thus we have:
\[
\func{Im}((g\circ f)\circ m_{I\widetilde{I}})=\func{Im}((g\circ f)\circ
m_{I})\cup \func{Im}((g\circ f)\circ m_{\widetilde{I}})\text{.}
\]
Now $(g\circ f)\circ m_{I}$ is an embedding since $f\circ m_{I}$ is an
embedding by ii) for $(f,\varphi ,\alpha )$ and $g$ restricted to $\func{Im}%
(f\circ m_{I})\subseteq X^{\prime }\setminus \func{Im}(m_{I^{\prime
}J^{\prime }}^{\prime })$ is an isomorphism by iv) for $(g,\psi ,\beta )$.
Furthermore $(g\circ f)\circ m_{\widetilde{I}}$ is an embedding, since by
the property $\mathcal{M}$-\ref{factprop}d) there exist an ordered set $%
\widetilde{I}\varsubsetneq R$ and an isomorphism $\alpha _{\widetilde{I}%
,I^{\prime }}$ such that $\alpha \circ s_{\widetilde{I},R}=s_{I^{\prime
},N^{\prime }}\circ \alpha _{\widetilde{I},I^{\prime }}$. Thus we derive
\begin{eqnarray*}
f\circ m_{\widetilde{I}} &=&f\circ m_{R}\circ s_{\widetilde{I},R} \\
&=&m^{\prime }\circ \alpha \circ s_{\widetilde{I},R} \\
&=&m^{\prime }\circ s_{I^{\prime },N^{\prime }}\circ \alpha _{\widetilde{I}%
,I^{\prime }} \\
&=&m_{I^{\prime }}^{\prime }\circ \alpha _{\widetilde{I},I^{\prime }}\text{,}
\end{eqnarray*}
which on composing with $g$, gives $g\circ f\circ m_{\widetilde{I}}=g\circ
m_{I^{\prime }}^{\prime }\circ \alpha _{\widetilde{I},I^{\prime }}$. Now
since $g\circ m_{I^{\prime }}^{\prime }$ is an embedding from ii) for $%
(g,\psi ,\beta )$, and $\alpha _{\widetilde{I},I^{\prime }}$ is an
isomorphism, we conclude that $g\circ f\circ m_{\widetilde{I}}$ is an
embedding.

Next we show that $\func{Im}((g\circ f)\circ m_{I})$ and $\func{Im}((g\circ
f)\circ m_{\widetilde{I}})$ are separated in $X^{\prime \prime }$, which,
together with the previous results proves that $(g\circ f)\circ m_{I%
\widetilde{I}}$ is an embedding. For this we set $C:=\func{Im}(f\circ m_{I})$
and $D:=\func{Im}(f\circ m_{\widetilde{I}})=\func{Im}(m_{I^{\prime
}}^{\prime }\circ \alpha _{\widetilde{I},I^{\prime }})=\func{Im}%
(m_{I^{\prime }}^{\prime })$. We need to show that $\func{Cls}(g(C))\cap
g(D)=\emptyset $ and $g(C)\cap \func{Cls}(g(D))=\emptyset $. We have, by iv)
for $(g,\psi ,\beta )$, $\func{Cls}(g(C))=g(\func{Cls}(C))\subseteq
X^{\prime \prime }\setminus g(D)$, since $\func{Cls}(C)\subseteq X^{\prime
}\setminus \func{Im}(m_{I^{\prime }J^{\prime }}^{\prime })$, by ii) for $%
(f,\varphi ,\alpha )$. Thus $\func{Cls}(g(C))\cap g(D)=\emptyset $. To show $%
g(C)\cap \func{Cls}(g(D))=\emptyset $, we have:
\begin{eqnarray*}
g(C)\cap \func{Cls}(g(D)) &=&g(C)\cap g(\func{Cls}(D)) \\
&=&(g(C)\cap g(D))\cup (g(C)\cap g(\func{Cls}(D)\setminus D)) \\
&=&g(C)\cap g(\func{Cls}(D)\setminus D) \\
&=&g(C\cap (\func{Cls}(D)\setminus D)) \\
&=&\emptyset \text{.}
\end{eqnarray*}
The first equality follows from ii) for $(g,\psi ,\beta )$, which says that $%
g(\func{Cls}(D))$ is closed in $X^{\prime \prime }$, since $g(D)\subseteq g(%
\func{Cls}(D))$, hence $\func{Cls}(g(D))\subseteq g(\func{Cls}(D))$ and
hence, because $g$ is continuous, $\func{Cls}(g(D))=g(\func{Cls}(D))$. The
second follows from properties of operations on sets. The third is clear
since we have $g(C)\cap g(D)=\emptyset $. The fourth follows from the fact
that $g$ restricted to $X^{\prime }\setminus \func{Im}(m_{I^{\prime
}J^{\prime }}^{\prime })$ is an isomorphism, by iv) for $(g,\psi ,\beta )$
and by the fact that $\func{Cls}(D)\setminus D\subseteq X^{\prime }\setminus
\func{Im}(m_{I^{\prime }J^{\prime }}^{\prime })$. This inclusion comes from
the fact that $(\func{Cls}(D)\setminus D)\cap \func{Im}(m_{I^{\prime
}J^{\prime }}^{\prime })$ is indeed empty, since $A_{I^{\prime }}^{\prime }$
and $A_{J^{\prime }}^{\prime }$ are separated in $A_{I^{\prime }}^{\prime
}\sqcup A_{J^{\prime }}^{\prime }$ and $m_{I^{\prime }J^{\prime }}^{\prime }$
is an embedding by $\mathcal{M}$-1). The final equality follows from ii) for
$(f,\varphi ,\alpha )$, since $C=\func{Im}(f\circ m_{I})$ and $D=\func{Im}%
(f\circ m_{\widetilde{I}})\varsubsetneq \func{Im}(f\circ m_{R})$ are
separated in $X^{\prime }$.

Now we show that $\func{Im}(g\circ f\circ m_{I\widetilde{I}})$ and $\func{Im}%
(m^{\prime \prime })$ are separated in $X^{\prime \prime }$. This holds
since $C$ and $\func{Im}(m_{R^{\prime }}^{\prime })\subseteq \func{Im}%
(m^{\prime })$ are separated in $X^{\prime }$, by ii) for $(f,\varphi
,\alpha )$, and $g$ is an isomorphism restricted to $X^{\prime }\setminus
\func{Im}(m_{I^{\prime }J^{\prime }}^{\prime })$, so that $g(C)$ and $\func{%
Im}(g\circ m_{R^{\prime }}^{\prime })=\func{Im}(m^{\prime \prime }\circ
\beta )=\func{Im}(m^{\prime \prime })$ are separated, by iv) for $(g,\psi
,\beta )$. Furthermore $\func{Im}(g\circ f\circ m_{\widetilde{I}})=g(D)$ and
$\func{Im}(m^{\prime \prime })$ are separated in $X^{\prime \prime }$, by
ii) for $(g,\psi ,\beta )$.

It remains to show that $(g\circ f)(\func{Cls}(\func{Im}(m_{I})))$ and $%
(g\circ f)(\func{Cls}(\func{Im}(m_{\widetilde{I}})))$ are closed. For this
we have:
\[
(g\circ f)(\func{Cls}(\func{Im}(m_{I})))=g(\func{Cls}(C))=\func{Cls}(g(C))%
\text{,}
\]
which is closed. The first equality follows from the fact that $f(\func{Cls}(%
\func{Im}(m_{I})))$ is closed, by ii) for $(f,\varphi ,\alpha )$ and the
fact that $C=\func{Im}(f\circ m_{I})\subseteq f(\func{Cls}(\func{Im}%
(m_{I}))) $, hence, because $f$ is continuous, $\func{Cls}(C)=f(\func{Cls}(%
\func{Im}(m_{I})))$. The second follows from the fact that $\func{Cls}%
(C)\subseteq X^{\prime }\setminus \func{Im}(m_{I^{\prime }J^{\prime
}}^{\prime })$ and from iv) for $(g,\psi ,\beta )$. For $(g\circ f)(\func{Cls%
}(\func{Im}(m_{\widetilde{I}})))$ we have:
\[
(g\circ f)(\func{Cls}(\func{Im}(m_{\widetilde{I}})))=g(\func{Cls}(f(\func{Im}%
(m_{\widetilde{I}}))))=g(\func{Cls}(D))\text{,}
\]
which is closed by ii) for $(g,\psi ,\beta )$. The first equality follows
from the fact that $\func{Cls}(\func{Im}(m_{\widetilde{I}}))\subseteq
X\setminus \func{Im}(m_{IJ})$ and from iv) for $(f,\varphi ,\alpha )$%
.\smallskip \newline
iii) We must show $(g\circ f\circ m_{J\widetilde{J}})\circ (\varphi \sqcup
(\alpha _{\widetilde{J}J^{\prime }}^{-1}\circ \psi \circ \alpha _{\widetilde{%
I}I^{\prime }}))=g\circ f\circ m_{I\widetilde{I}}$. Let $x\in A_{I\widetilde{%
I}}$. Then either there exists $y\in A_{I}$ such that $x=s_{I,I\widetilde{I}%
}(y)$ or there exists $z\in A_{\widetilde{I}}$ such that $x=s_{\widetilde{I}%
,I\widetilde{I}}(z)$, since the morphisms $s_{I,I\widetilde{I}}$ and $s_{%
\widetilde{I},I\widetilde{I}}$ are the canonical monomorphisms in $\mathrm{%
Top}$. In the first case we have:
\begin{eqnarray*}
m_{J\widetilde{J}}\circ (\varphi \sqcup (\alpha _{\widetilde{J}J^{\prime
}}^{-1}\circ \psi \circ \alpha _{\widetilde{I}I^{\prime }}))(x) &=&m_{J%
\widetilde{J}}\circ (\varphi \sqcup (\alpha _{\widetilde{J}J^{\prime
}}^{-1}\circ \psi \circ \alpha _{\widetilde{I}I^{\prime }}))\circ s_{I,I%
\widetilde{I}}(y) \\
&=&m_{J\widetilde{J}}\circ s_{J,J\widetilde{J}}\circ \varphi (y) \\
&=&m_{J}\circ \varphi (y)
\end{eqnarray*}
using $\mathcal{M}$-\ref{factprop}d) in the second equality. Thus
\begin{eqnarray*}
(g\circ f\circ m_{J\widetilde{J}})\circ (\varphi \sqcup (\alpha _{\widetilde{%
J}J^{\prime }}^{-1}\circ \psi \circ \alpha _{\widetilde{I}I^{\prime }}))(x)
&=&(g\circ f\circ m_{J}\circ \varphi )(y) \\
&=&(g\circ f\circ m_{I})(y) \\
&=&(g\circ f\circ m_{I\widetilde{I}}\circ s_{I,I\widetilde{I}})(y) \\
&=&(g\circ f\circ m_{I\widetilde{I}})(x)\text{,}
\end{eqnarray*}
using (iii) for $(f,\varphi ,\alpha )$ in the second equality. In the second
case we have:
\begin{eqnarray*}
m_{J\widetilde{J}}\circ (\varphi \sqcup (\alpha _{\widetilde{J}J^{\prime
}}^{-1}\circ \psi \circ \alpha _{\widetilde{I}I^{\prime }}))(x) &=&(m_{J%
\widetilde{J}}\circ (\varphi \sqcup (\alpha _{\widetilde{J}J^{\prime
}}^{-1}\circ \psi \circ \alpha _{\widetilde{I}I^{\prime }}))\circ s_{%
\widetilde{I},I\widetilde{I}})(z) \\
&=&(m_{J\widetilde{J}}\circ s_{\widetilde{J},J\widetilde{J}}\circ \alpha _{%
\widetilde{J}J^{\prime }}^{-1}\circ \psi \circ \alpha _{\widetilde{I}%
I^{\prime }})(z) \\
&=&(m_{\widetilde{J}}\circ \alpha _{\widetilde{J}J^{\prime }}^{-1}\circ \psi
\circ \alpha _{\widetilde{I}I^{\prime }})(z)\text{,}
\end{eqnarray*}
using property $\mathcal{M}$-\ref{factprop}d) in the second equality. Thus
\begin{eqnarray*}
(g\circ f\circ m_{J\widetilde{J}})\circ (\varphi \sqcup (\alpha _{\widetilde{%
J}J^{\prime }}^{-1}\circ \psi \circ \alpha _{\widetilde{I}I^{\prime }}))(x)
&=&(g\circ f\circ m_{\widetilde{J}}\circ \alpha _{\widetilde{J}J^{\prime
}}^{-1}\circ \psi \circ \alpha _{\widetilde{I}I^{\prime }})(z) \\
&=&(g\circ m_{J^{\prime }}^{\prime }\circ \psi \circ \alpha _{\widetilde{I}%
I^{\prime }})(z) \\
&=&(g\circ m_{I^{\prime }}^{\prime }\circ \alpha _{\widetilde{I}I^{\prime
}})(z) \\
&=&(g\circ f\circ m_{\widetilde{I}})(z) \\
&=&(g\circ f\circ m_{I\widetilde{I}}\circ s_{\widetilde{I},I\widetilde{I}%
})(z) \\
&=&(g\circ f\circ m_{I\widetilde{I}})(x)\text{,}
\end{eqnarray*}
using (i) for $(f,\varphi ,\alpha )$ and $\mathcal{M}$-\ref{factprop}d) in
the second and fourth equalities and (iii) for $(g,\psi ,\beta )$ in the
third equality.\smallskip \newline
iv) Finally for this condition we have:
\begin{eqnarray*}
X\setminus \func{Im}(m_{IJ\widetilde{I}\widetilde{J}}) &\cong &(X^{\prime
}\setminus \func{Im}(f\circ m_{IJ}))\setminus \func{Im}(f\circ m_{\widetilde{%
I}\widetilde{J}}) \\
&=&(X^{\prime }\setminus \func{Im}(f\circ m_{I}))\setminus \func{Im}%
(m_{I^{\prime }J^{\prime }}^{\prime }) \\
&=&(X^{\prime }\setminus \func{Im}(m_{I^{\prime }J^{\prime }}^{\prime
}))\setminus \func{Im}(f\circ m_{I}) \\
&\cong &(X^{\prime \prime }\setminus \func{Im}(g\circ m_{I^{\prime
}}^{\prime }))\setminus \func{Im}(g\circ f\circ m_{I}) \\
&=&(X^{\prime \prime }\setminus \func{Im}(g\circ f\circ m_{I}))\setminus
\func{Im}(g\circ f\circ m_{\widetilde{I}}) \\
&=&(X^{\prime \prime }\setminus \func{Im}(g\circ f\circ m_{I}))\setminus
\func{Im}(g\circ f\circ m_{\widetilde{I}}) \\
&=&X^{\prime \prime }\setminus (\func{Im}(g\circ f\circ m_{I})\cup \func{Im}%
(g\circ f\circ m_{\widetilde{I}})) \\
&=&X^{\prime \prime }\setminus \func{Im}(g\circ f\circ m_{I\widetilde{I}})%
\text{.}
\end{eqnarray*}
These equalities and isomorphisms are proved by using iv) for the gluing
morphisms $(f,\varphi ,\alpha )$ and $(g,\psi ,\beta )$ and general
properties of the operation of difference of sets.%
\endproof%

Furthermore we have:

\begin{theorem}
The composition in $\mathcal{C}$ is associative.
\end{theorem}

\proof%
We will show this for the case of three gluing morphisms. The other
combinations are easily checked. We will consider the gluing morphisms
\[
(f,\varphi ,\alpha ):(X,A,m)\rightarrow (X^{\prime },A^{\prime },m^{\prime })%
\text{,}
\]
where $\varphi :A_{I}\rightarrow P(A)_{J}$, $\alpha :A_{R}\rightarrow
A^{\prime }$,
\[
(g,\psi ,\beta ):(X^{\prime },A^{\prime },m^{\prime })\rightarrow (X^{\prime
\prime },A^{\prime \prime },m^{\prime \prime })\text{,}
\]
where $\psi :A_{I^{\prime }}^{\prime }\rightarrow P(A^{\prime })_{J^{\prime
}}$, $\beta :A_{R^{\prime }}^{\prime }\rightarrow A^{\prime \prime }$ and
\[
(h,\rho ,\gamma ):(X^{\prime \prime },A^{\prime \prime },m^{\prime \prime
})\rightarrow (X^{\prime \prime \prime },A^{\prime \prime \prime },m^{\prime
\prime \prime })\text{,}
\]
where $\rho :A_{I^{\prime \prime }}^{\prime \prime }\rightarrow P(A^{\prime
\prime })_{J^{\prime \prime }}$ and $\gamma :A_{R^{\prime \prime }}^{\prime
\prime }\rightarrow A^{\prime \prime \prime }$.\newline
By the property $\mathcal{M}$-\ref{factprop}d) applied to $\alpha $,
corresponding to the ordered subsets $I^{\prime }$, $J^{\prime }$, $%
R^{\prime }$, there exist ordered subsets $\widetilde{I},\widetilde{J},%
\widetilde{R}\varsubsetneq R$, and isomorphisms $\alpha _{\widetilde{I}%
,I^{\prime }}$, $\alpha _{\widetilde{J},J^{\prime }}$, $\alpha _{\widetilde{R%
},R^{\prime }}$. By the same property applied to $\beta $, corresponding to $%
I^{\prime \prime }$, $J^{\prime \prime }$, $R^{\prime \prime }$ there exist
ordered subsets $\widetilde{I}^{\prime },\widetilde{J}^{\prime },\widetilde{R%
}^{\prime }\varsubsetneq R^{\prime }$ and isomorphisms $\beta _{\widetilde{I}%
^{\prime },I^{\prime \prime }}$, $\beta _{\widetilde{J}^{\prime },J^{\prime
\prime }}$, $\beta _{\widetilde{R}^{\prime },R^{\prime \prime }}$. Moreover,
applying property $\mathcal{M}$-\ref{factprop}d) to $\alpha _{\widetilde{R}%
,R^{\prime }}$ corresponding to the subsets $\widetilde{I}^{\prime },%
\widetilde{J}^{\prime },\widetilde{R}^{\prime }\varsubsetneq R^{\prime }$
there exist ordered subsets $\widehat{I},\widehat{J},\widehat{R}%
\varsubsetneq R$ and isomorphisms $\alpha _{\widehat{I},\widetilde{I}%
^{\prime }}$, $\alpha _{\widehat{J},\widetilde{J}^{\prime }}$ and $\alpha _{%
\widehat{R},\widetilde{R}^{\prime }}$. For one side of the associativity
equation, we have (denoting composition by juxtaposition for simplicity of
notation):
\begin{eqnarray*}
(h,\rho ,\gamma )\left( (g,\psi ,\beta )(f,\varphi ,\alpha )\right)
&=&(h,\rho ,\gamma )\left( gf,\varphi \sqcup (P(\alpha ^{-1})_{\widetilde{J}%
,J^{\prime }}\psi \alpha _{\widetilde{I},I^{\prime }}),\beta \alpha _{%
\widetilde{R},R^{\prime }}\right) \\
&=&(h(gf),\varphi \sqcup (P(\alpha ^{-1})_{_{\widetilde{J},J^{\prime }}}\psi
\alpha _{\widetilde{I},I^{\prime }})\sqcup (P(\beta \alpha _{\widetilde{R}%
,R^{\prime }})_{\widehat{J},J^{\prime \prime }}^{-1} \\
&&\rho (\beta \alpha _{\widetilde{R},R^{\prime }})_{\widehat{I},I^{\prime
\prime }}),\gamma (\beta \alpha _{\widetilde{R},R^{\prime }})_{\widehat{I}%
,I^{\prime \prime }}) \\
&=&(h(gf),\varphi \sqcup (P(\alpha ^{-1})_{_{\widetilde{J},J^{\prime }}}\psi
\alpha _{\widetilde{I},I^{\prime }})\sqcup (P(\alpha ^{-1})_{\widehat{J},%
\widetilde{J}^{\prime }} \\
&&P(\beta ^{-1})_{\widetilde{J}^{\prime },J^{\prime \prime }}\rho \beta _{%
\widetilde{I^{\prime }},I^{\prime \prime }}\alpha _{\widehat{I},\widetilde{I}%
^{\prime }}),\gamma (\beta _{\widetilde{R}^{\prime },R^{\prime \prime
}}\alpha _{\widehat{R},\widetilde{R}^{\prime }}))\text{.}
\end{eqnarray*}
For the other side we have:
\begin{eqnarray*}
\left( (h,\rho ,\gamma )(g,\psi ,\beta )\right) (f,\varphi ,\alpha )
&=&(hg,\psi \sqcup (P(\beta ^{-1})_{\widetilde{J^{\prime }},J^{\prime \prime
}}\rho \beta _{\widetilde{I^{\prime }},I^{\prime \prime }}),\gamma \beta _{%
\widetilde{R^{\prime }},R^{\prime \prime }})(f,\varphi ,\alpha ) \\
&=&((hg)f,\varphi \sqcup ((P(\alpha ^{-1})_{\widetilde{J}\widehat{J}%
,J^{\prime }\widetilde{J}^{\prime }}(\psi \sqcup P(\beta ^{-1})_{\widetilde{J%
}^{\prime },J^{\prime \prime }}\rho \beta _{\widetilde{I}^{\prime
},I^{\prime \prime }}) \\
&&\alpha _{\widetilde{I}\widehat{I},I^{\prime }\widetilde{I}^{\prime
}},(\gamma \beta _{\widetilde{R}^{\prime },R^{\prime \prime }})\alpha _{%
\widehat{R},\widetilde{R}^{\prime }}) \\
&=&((hg)f,\varphi \sqcup (P(\alpha ^{-1})_{\widetilde{J},J^{\prime }}\psi
\alpha _{\widetilde{I},I^{\prime }})\sqcup (P(\alpha ^{-1})_{\widehat{J},%
\widetilde{J}^{\prime }} \\
&&P(\beta ^{-1})_{\widetilde{J}^{\prime },J^{\prime \prime }}\rho \beta _{%
\widetilde{I}^{\prime },I^{\prime \prime }}\alpha _{\widehat{I},\widetilde{I}%
^{\prime }}),(\gamma \beta _{\widetilde{R}^{\prime },R^{\prime \prime
}})\alpha _{\widehat{R},\widetilde{R}^{\prime }})\text{,}
\end{eqnarray*}
using the interchange law in the third equality. The two expressions are
equal because of the associativity of composition in $\mathbf{C}$.%
\endproof%
\medskip

\begin{definition}
\label{Defidentities}Let $(X,A,m)\in \func{Ob}(\mathcal{C})$. The identity
morphism on $(X,A,m)$ is defined as follows:
\[
\limfunc{id}\nolimits_{(X,A,m)}=(\limfunc{id}\nolimits_{X},\limfunc{id}%
\nolimits_{A}):(X,A,m)\rightarrow (X,A,m)\text{.}
\]
\end{definition}

\begin{theorem}
\label{Thm215}$\mathcal{C}$ of Definition \ref{Def213} is a category.
\end{theorem}

\proof%
This follows from the results already shown above and the fact that the
morphisms of Definition \ref{Defidentities} obviously satisfy the
requirements to be identity morphisms.%
\endproof%
\medskip

\begin{definition}
\label{Def216}A monoidal product, unit object and braiding on $\mathcal{C}$
are defined as follows:\newline
monoidal product on objects:
\[
(X,A,m)\sqcup (X^{\prime },A^{\prime },m^{\prime })=(X\sqcup X^{\prime
},A\sqcup A^{\prime },m\sqcup m^{\prime })\text{,}
\]
\newline
monoidal product on morphisms:
\begin{eqnarray*}
(f,\alpha )\sqcup (g,\beta ) &=&(f\sqcup g,\alpha \sqcup \beta )\text{,} \\
(f,\varphi ,\alpha )\sqcup (g,\beta ) &=&(f\sqcup g,\varphi ,\alpha \sqcup
\beta )\text{,} \\
(f,\alpha )\sqcup (g,\psi ,\beta ) &=&(f\sqcup g,\psi ,\alpha \sqcup \beta )%
\text{,} \\
(f,\varphi ,\alpha )\sqcup (g,\psi ,\beta ) &=&(f\sqcup g,\varphi \sqcup
\psi ,\alpha \sqcup \beta )\text{,}
\end{eqnarray*}
unit object:
\[
\widetilde{E}=(E,E,\limfunc{id}\nolimits_{E})
\]
braiding:
\[
c_{((X,A,m),(X^{\prime },A^{\prime },m^{\prime }))}=(c_{(X,X^{\prime
})},c_{(A,A^{\prime })})\text{.}
\]
\end{definition}

\begin{theorem}
$(\mathcal{C},\sqcup ,\widetilde{E},c)$ is a symmetric, strict monoidal
category.
\end{theorem}

\proof%
The monoidal product on $\mathcal{C}$ closes on objects by 2) of Definition
\ref{Def213}, and clearly the four morphisms in Definition \ref{Def216} are
isomorphisms or gluing morphisms, so the monoidal product closes on
morphisms too. The monoidal product is strictly associative in $\mathcal{C}$%
, since it is strictly associative in $\mathbf{C}$ and $\mathcal{S}(\mathbf{C%
})$. The unit object $\widetilde{E}$ is an object of $\mathcal{C}$ because
of 1) of Definition \ref{Def213}, and is obviously a strict unit, since $E$
is a strict unit in $\mathbf{C}$ and $\mathcal{S}(\mathbf{C})$. The braiding
defined in Definition \ref{Def216} is well-defined because of 4) of
Definition \ref{Def213} and is a braiding, because $c$ is a braiding in $%
\mathbf{C}$ and $\mathcal{S}(\mathbf{C})$. It also clearly satisfies the
symmetry condition.%
\endproof%
\medskip

For some TQFT functors only the endofunctor $P$ on $\mathcal{S}(\mathbf{C})$
will play a role. For others we need a corresponding endofunctor on $%
\mathcal{C}$.

\begin{definition}
\label{Def218}The endofunctor $(P,\pi _{2},\pi _{0})$ on $\mathbf{C}$
extends to an endofunctor $(\mathbf{P},\mathbf{\pi }_{2},\mathbf{\pi }_{0})$
on $\mathcal{C}$ as follows:\newline
on objects:
\[
\mathbf{P}(X,A,m)=(P(X),P(A),P(m))\text{,}
\]
on morphisms
\begin{eqnarray*}
\mathbf{P}(f,\alpha ) &=&(P(f),P(\alpha ))\text{,} \\
\mathbf{P}(f,\varphi ,\alpha ) &=&(P(f),P(\varphi ),P(\alpha ))\text{,}
\end{eqnarray*}
and the natural isomorphism $\mathbf{\pi }_{2}$ and the isomorphism $\mathbf{%
\pi }_{0}$ are given by:
\begin{eqnarray*}
\mathbf{\pi }_{2((X,A,m),(X^{\prime },A^{\prime },m^{\prime }))} &=&(\pi
_{2(X,X^{\prime })},\pi _{2(A,A^{\prime })})\text{,} \\
\mathbf{\pi }_{0} &=&(\pi _{0},\pi _{0})\text{.}
\end{eqnarray*}
\end{definition}

\begin{theorem}
$(\mathbf{P},\mathbf{\pi }_{2},\mathbf{\pi }_{0})$ is a symmetric, strong
monoidal endofunctor on $\mathcal{C}$.
\end{theorem}

\proof%
The natural isomorphism $\mathbf{\pi }_{2((X,A,m),(X^{\prime },A^{\prime
},m^{\prime }))}$ and the isomorphism $\mathbf{\pi }_{0}$ have domain and
codomain in $\func{Ob}(\mathcal{C})$, because of $1)$, $2)$ and $3)$ of
Definition \ref{Def213}, and are isomorphisms of $\mathcal{C}$ since $P$ is
strong monoidal and because of the naturality of $\pi _{2}$. The naturality
of $\mathbf{\pi }_{2}$ also follows directly from the naturality of $\pi
_{2} $, by combining two commuting diagrams. The conditions on $\mathbf{\pi }%
_{2}$ and $\mathbf{\pi }_{0}$ follow directly from combining two
corresponding diagrams for $\pi _{2}$ and $\pi _{0}$.%
\endproof%
\medskip

We can now conclude this section by defining two types of topological
categories, corresponding to two types of algebraic categories and two types
of TQFT functors, as we shall see later.

\begin{definition}
Let $(\mathbf{C},F,P)$ be a topological starting category, $\mathcal{M}$ be
an appropriate class of monomorphisms and $\mathcal{O}$ be a class of
objects of $\mathcal{C}$ satisfying Definition \ref{Def213}.

The \emph{topological category} for these data is the pair $(\mathcal{C},P)$%
, where $\mathcal{C}$ denotes the symmetric, strict monoidal category $(%
\mathcal{C},\sqcup ,\widetilde{E},c)$ and $P$ is the monoidal endofunctor on
$\mathcal{S}(\mathbf{C})$ defined in Definition \ref{DefN27}.\newline
The \emph{full topological category} for these data is the pair $(\mathcal{C}%
,\mathbf{P})$, where $\mathcal{C}$ is as before and $\mathbf{P}$ is the
monoidal endofunctor on $\mathcal{C}$ defined in Definition \ref{Def218}%
.\medskip
\end{definition}

Thus we have achieved our objective for this section, namely to define a
general-purpose topological category capable of describing gluing and change
of orientation in a generalized sense. Although the notation for the objects
and morphisms of the topological category is a little cumbersome, it allows
all features of a gluing operation to be specified clearly. The topological
gluing operation lies at the heart of TQFT, and the topological category we
have constructed distills out precisely that feature, suppressing all
morphisms which do not fit into the scheme of gluing. Indeed the gluing is
described at the lowest categorical level, namely at the level of morphisms,
and the only extra structures on the category are the symmetric monoidal
structure and, of course, the endofunctor representing change of
orientation. This structure is about as simple as it can be, given the broad
scope of the framework.

\section{The algebraic category\label{section3}}

The construction of the algebraic category in this section again proceeds in
several stages, like in the previous section. The starting point is a
category, the algebraic starting category, whose objects are
(finite-dimensional) $K$-modules, where $K$ is a ring, possibly endowed with
further structures such as an inner product. Then we introduce a so-called
evaluation on each object, which behaves very much like a hermitian
structure on a $K$-module. This gives rise to a category called the
evaluation-preserving category, which has the same objects as the algebraic
starting category, and whose morphisms are isomorphisms preserving the
evaluations. The evaluation-preserving category is the algebraic counterpart
to the category of subobjects which appears in the construction of the
topological category. The algebraic category itself has as its objects pairs
consisting of an object of the algebraic starting category and an element of
its underlying $K$-module. The morphisms of the algebraic category are
morphisms of the algebraic starting category preserving the elements. (Thus
a morphism implies an equation, and these equations will be crucial later on
in determining TQFT functors.) Again there is an extra piece of structure
which plays a role throughout the construction, coming from an endofunctor
on the algebraic starting category. In the example, this is given by the
operation of replacing a $K$-module by the same module with the same
addition, but with a different scalar multiplication which uses the complex
conjugate.

We now proceed with the details of the construction.

\begin{definition}
An \emph{algebraic starting category} is a triple $(\mathbf{D},G,Q)$, where $%
\mathbf{D}$ is a symmetric, strict monoidal category $(\mathbf{D},\otimes
,I,c)$, $G$ is a symmetric, strict monoidal forgetful functor from $\mathbf{D%
}$ to $K$-$\mathrm{Mod}$ (the category whose objects are finitely generated $%
K$-modules, where $K$ is a fixed commutative ring with unit, and whose
morphisms are $K$-homomorphisms) with its standard monoidal structure
(tensor product) and braiding, and $Q$ is a symmetric, strong monoidal
endofunctor $(Q,\theta _{2},\theta _{0})$ on $\mathbf{D}$.
\end{definition}

\begin{thm_remark}
The forgetful functor $G$ justifies the adjective algebraic. It sends each
object $V$ of $\mathbf{D}$ to an underlying finitely generated $K$-module $%
G(V)$ and each morphism $f$ of $\mathbf{D}$ to an underlying $K$%
-homomorphism $G(f)$. We may think of $Q$ as assigning to each module $V$
its ``conjugate'' module, a notion which will become clear in the example.
\end{thm_remark}

{\noindent \bf Example: }%
Let $V$ be a finite-dimensional vector space over $\Bbb{C}$. We denote by $%
\overline{V}$ the vector space which has the same elements and addition
operation as $V$, but a different operation of scalar multiplication denoted
$\overline{\cdot }$, defined by $c\rule{0.05in}{0in}\overline{\cdot }\rule%
{0.05in}{0in}x=\overline{c}x$, where $c\in \Bbb{C}$, $x\in V$, $\overline{c}$
denotes the complex conjugate of $c$ and scalar multiplication in $V$ is
denoted by juxtaposition. (We remark that the example could equally well be
carried through using any commutative ring with an involution, instead of $%
\Bbb{C}$ and the operation of conjugation.) The objects of $\mathbf{D}$ are
finite-dimensional complex vector spaces generated by $\Bbb{C}$ and $V$ via
the operation of tensor products, and the unary operation `$^{-}$' just
introduced. Thus an example of an object of $\mathbf{D}$ is $\overline{%
V\otimes \overline{V}}$. The motivation for taking this class of objects for
$\mathbf{D}$, instead of simply all finite-dimensional vector spaces over $%
\Bbb{C}$, is that it is of minimal size for providing a good target category
for the TQFT's in our example\footnote{%
Roughly speaking, for the description of the TQFT functor we only need to
give one vector space for each type of irreducible subobject on the
topological side, and in this example there are just two types: positively
oriented circles and negatively oriented circles.}.\textbf{\ }The morphisms
of $\mathbf{D}$ are $\Bbb{C}$-linear maps between the objects of $\mathbf{D}$%
. The monoidal structure is given by the tensor product $\otimes $, taken to
be strict (thus, e.g. $V\otimes \Bbb{C}=V$), and the unit object $I$, also
strict, is $\Bbb{C}$. The braiding is given by the usual flip map
\[
c_{(W,Y)}:W\otimes Y\rightarrow Y\otimes W\text{,}\hspace{0.2in}w\otimes
y\mapsto y\otimes w\text{,}
\]
extended by linearity, which clearly satisfies the condition to be a
symmetry. The forgetful functor $G$ formally maps each object of $\mathbf{D}$
to the corresponding module over $\Bbb{C}$, regarded as a commutative ring,
instead of a field, but $G$ will be treated as the identity functor from now
on.

The endofunctor $Q$ acts on objects by
\[
Q(W)=\overline{W}\text{.}
\]
Before describing the action of $Q$ on morphisms it is useful to introduce
the following ``identity'' map,
\[
k_{W}:W\rightarrow \overline{W},
\]
defined by $k_{W}(x)=x$ for all $x$ in $W$. Note that, although it preserves
addition, $k_{W}$ is not a morphism of $\mathbf{D}$, since $k_{W}(cx)=cx=%
\overline{c}\stackrel{\_}{\cdot }x=\overline{c}\stackrel{\_}{\cdot }k_{W}(x)$%
. We now define
\[
Q(Y\stackrel{f}{\rightarrow }W)=\overline{Y}\stackrel{Q(f)}{\rightarrow }%
\overline{W}\text{,}
\]
where
\[
Q(f)=k_{W}\circ f\circ k_{Y}^{-1}
\]
is a linear map, i.e. a morphism of $\mathbf{D}$, as is easily checked.
Finally the corresponding isomorphism
\[
\theta _{0}:\Bbb{C}\rightarrow \overline{\Bbb{C}}
\]
is given by
\[
\theta _{0}(z)=k_{\Bbb{C}}(\overline{z})
\]
(note the conjugation in the argument which ensures that $\theta _{0}$ is
linear), and the natural isomorphism $\theta _{2}$ is given by isomorphisms
\[
\theta _{2(Y,W)}:\overline{Y}\otimes \overline{W}\rightarrow \overline{%
Y\otimes W},
\]
where
\[
\theta _{2(Y,W)}=k_{Y\otimes W}\circ (k_{Y}^{-1}\otimes k_{W}^{-1})\text{,}
\]
is again a morphism of $\mathbf{D}$. Here the tensor product of $k_{Y}^{-1}$
and $k_{W}^{-1}$ is defined by
\[
y\otimes w\mapsto k_{Y}^{-1}(y)\otimes k_{W}^{-1}(w)
\]
which is well-defined and preserves addition, but not scalar multiplication.
We omit the details of the simple check that $(Q,\theta _{2},\theta _{0})$
satisfies the definition of a monoidal endofunctor.$\blacktriangle \medskip $

We next introduce another category, which we denote by $\mathcal{S}(\mathbf{D%
})$, and which will correspond to the category $\mathcal{S}(\mathbf{C})$ of
the previous section, when we start to talk about TQFT functors in the next
section. $\mathcal{S}(\mathbf{D})$ has the same class of objects as $\mathbf{%
D}$, but its morphisms are restricted to a certain subclass of morphisms,
namely isomorphisms preserving the so-called evaluations, which we now
introduce.

\begin{definition}
A \emph{choice of evaluations} on an algebraic starting category $(\mathbf{D}%
,G,Q)$ is an assignment to each $V\in \func{Ob}(\mathbf{D})$ of a morphism $%
e_{V}:Q(V)\otimes V\rightarrow I$ such that:

\begin{itemize}
\item  (multiplicativity axiom)

for each pair of objects $(V,W)$ of $\mathcal{S}(\mathbf{D)}$%
\[
e_{V\otimes W}\circ (\theta _{2(V,W)}\otimes \limfunc{id}\nolimits_{V\otimes
W})\circ (\limfunc{id}\nolimits_{Q(V)}\otimes c_{(V,Q(W))}\otimes \limfunc{id%
}\nolimits_{W})=e_{V}\otimes e_{W}\text{,}
\]
i.e. the following diagram%
\begin{center}
\setsqparms[1`1`1`1;1600`700]
\square[Q(V)\otimes V\otimes Q(W)\otimes W`I\otimes I=I`Q(V\otimes W)\otimes V\otimes W`I;e_{V}\otimes e_{W}`t`\limfunc{id}\nolimits_{I}`e _{V\otimes W}]
\newline
\newline
(where $t=(\theta _{2(V,W)}\otimes \limfunc{id}\nolimits_{V\otimes W})\circ (\limfunc{id}\nolimits_{Q(V)}\otimes c_{(V,Q(W))}\otimes \limfunc{id}\nolimits_{W})$)
\end{center}%
is commutative.

\item  (conjugation axioms)

a) for each object $V$ of $\mathcal{S}(\mathbf{D)}$,
\[
Q(e_{V})\circ \theta _{2(Q(V),V)}=\theta _{0}\circ e_{Q(V)}\text{,}
\]
i.e. the following diagram%
\begin{center}
\setsqparms[1`1`1`1;1400`700]
\square[Q(Q(V))\otimes Q(V)`Q(Q(V)\otimes V)`I`Q(I);\theta _{2(Q(V),V)}`e_{Q(V)}`Q(e_{V})`\theta _{0}]
\end{center}%
is commutative.

b) for each pair $(V,W)$ of objects of $\mathcal{S}(\mathbf{D)}$,
\[
e_{Q(V\otimes W)}\circ (Q(\theta _{2(V,W)})\otimes \theta
_{2(V,W)})=e_{Q(V)\otimes Q(W)}\text{.}
\]

c) for the unit object $I$,
\[
e_{Q(I)}\circ (Q(\theta _{0})\otimes \theta _{0})=e_{I}\text{.}
\]
\smallskip
\end{itemize}
\end{definition}

{\noindent \bf Example: }%
\label{examplep30}Because the objects of $\mathbf{D}$ in the example are
generated from $V$ and $\Bbb{C}$ by taking tensor products and applying $Q$,
it is enough to specify $e_{V}$, since the other evaluations follow from the
axioms. Let $(e_{i})_{i}$ be a basis of $V$, and denote by $(\overline{e}%
_{i})_{i}$ the same basis regarded as a basis of $\overline{V}$. Then we set
\[
e_{V}(\overline{e}_{i}\otimes e_{j})=a_{ij}\text{,}
\]
where $a_{ij}\in \Bbb{C}$. Now, by the conjugation axiom a), $e_{\overline{V}%
}$ is given by the composition $\theta _{0}^{-1}\circ \overline{e_{V}}\circ
\theta _{2(\overline{V},V)}$. A short calculation using the previous
definitions gives:
\[
e_{\overline{V}}(e_{i}\otimes \overline{e}_{j})=\overline{a}_{ij}\text{.%
\label{page22}}
\]
Similar calculations using the multiplicativity axiom with $W=\Bbb{C}$ give
\[
e_{\Bbb{C}}(k_{\Bbb{C}}(r)\otimes s)=\overline{r}s,
\]
and for $e_{\overline{\Bbb{C}}}$ we find
\[
e_{\overline{\Bbb{C}}}(r\otimes k_{\Bbb{C}}(s))=r\overline{s}.
\]
Finally to obtain, e.g. $e_{V\otimes \overline{V}}$, we apply the
multiplicativity axiom to get

\[
e_{V\otimes \overline{V}}(\overline{e_{i}\otimes \overline{e}_{j}}\otimes
e_{k}\otimes \overline{e}_{l})=a_{ik}\overline{a}_{jl}\text{.}
\]
We still need to check that the conjugation axioms b) and c) are satisfied
by these choices. For c) we have:
\begin{eqnarray*}
e_{\overline{\Bbb{C}}}\circ (\overline{\theta }_{0}\otimes \theta _{0})(k_{%
\Bbb{C}}(r)\otimes s) &=&e_{\overline{\Bbb{C}}}((k_{\overline{\Bbb{C}}}\circ
\theta _{0}\circ k_{\Bbb{C}}^{-1})\otimes \theta _{0})(k_{\Bbb{C}}(r)\otimes
s) \\
&=&e_{\overline{\Bbb{C}}}(k_{\overline{\Bbb{C}}}k_{\Bbb{C}}(\overline{r}%
)\otimes k_{\Bbb{C}}(\overline{s})) \\
&=&e_{\overline{\Bbb{C}}}(\overline{r}\otimes k_{\Bbb{C}}(\overline{s})) \\
&=&\overline{r}s \\
&=&e_{\Bbb{C}}(k_{\Bbb{C}}(r)\otimes s)\text{,}
\end{eqnarray*}
and we leave the check of b) to the reader.$\blacktriangle \medskip $

\begin{definition}
Given a choice of evaluations $(e_{V})_{V\in \func{Ob}(\mathbf{D})}$ on an
algebraic starting category $(\mathbf{D},G,Q)$, we define a symmetric,
strict monoidal category $(\mathcal{S}(\mathbf{D}),\otimes ,I,c)$, also
written $\mathcal{S}(\mathbf{D})$ for short, as follows:

\begin{enumerate}
\item[\emph{a)}]  $\func{Ob}(\mathcal{S}(\mathbf{D}))=\func{Ob}(\mathbf{D})$,

\item[\emph{b)}]  the morphisms of $\mathcal{S}(\mathbf{D})$ are the
isomorphisms of $\mathbf{D}$ which preserve the evaluations, in the
following sense:
\[
e_{W}\circ (Q(f)\otimes f)=e_{V}\text{,}
\]
for $f:V\rightarrow W$ belonging to $\limfunc{Iso}(\mathbf{D)}$,

\item[\emph{c)}]  the composition, identity morphisms, monoidal product $%
\otimes $, unit object $I$ and the braiding $c$ are inherited from $\mathbf{D%
}$.
\end{enumerate}
\end{definition}

\begin{theorem}
$(\mathcal{S}(\mathbf{D}),\otimes ,I,c)$ is a symmetric, strict monoidal
category.
\end{theorem}

\proof%
$\mathcal{S}(\mathbf{D)}$ is a category, since clearly all identity
morphisms are morphisms of $\mathcal{S}(\mathbf{D)}$, and since $\limfunc{Mor%
}(\mathcal{S}(\mathbf{D))}$ is closed under composition: for $f:V\rightarrow
W$ and $g:W\rightarrow U$ morphisms of $\mathcal{S}(\mathbf{D)}$ we have:
\begin{eqnarray*}
e_{U}\circ (Q(g\circ f)\otimes (g\circ f)) &=&e_{U}\circ ((Q(g)\circ
Q(f))\otimes (g\circ f)) \\
&=&e_{U}\circ ((Q(g)\otimes g)\circ (Q(f)\otimes f)) \\
&=&e_{W}\circ (Q(f)\otimes f) \\
&=&e_{V}
\end{eqnarray*}
using the interchange law in the second equality. $\limfunc{Mor}(\mathcal{S}(%
\mathbf{D))}$ is closed under the monoidal product since, for two morphisms
of $\mathcal{S}(\mathbf{D)}$ $f:V\rightarrow W$ and $g:U\rightarrow Y$, we
have:
\begin{eqnarray*}
e_{W\otimes Y}\circ (Q(f\otimes g)\otimes f\otimes g) &=&(e_{W}\otimes
e_{Y})\circ (\limfunc{id}\nolimits_{Q(W)}\otimes c_{(Q(Y),W)}\otimes
\limfunc{id}\nolimits_{Y})\circ \\
&&(\theta _{2(W,Y)}^{-1}\otimes \limfunc{id}\nolimits_{W\otimes Y})\circ
(Q(f\otimes g)\otimes f\otimes g) \\
&=&(e_{W}\otimes e_{Y})\circ (\limfunc{id}\nolimits_{Q(W)}\otimes
c_{(Q(Y),W)}\otimes \limfunc{id}\nolimits_{Y})\circ \\
&&(Q(f)\otimes Q(g)\otimes f\otimes g)\circ (\theta _{2(V,U)}^{-1}\otimes
\limfunc{id}\nolimits_{V\otimes U}) \\
&=&(e_{W}\otimes e_{Y})\circ (Q(f)\otimes f\otimes Q(g)\otimes g)\circ \\
&&(\limfunc{id}\nolimits_{Q(V)}\otimes c_{(Q(U),V)}\otimes \limfunc{id}%
\nolimits_{U})\circ (\theta _{2(V,U)}^{-1}\otimes \limfunc{id}%
\nolimits_{V\otimes U}) \\
&=&((e_{W}\circ (Q(f)\otimes f))\otimes (e_{Y}\circ (Q(g)\otimes g)))\circ \\
&&(\limfunc{id}\nolimits_{Q(V)}\otimes c_{(Q(U),V)}\otimes \limfunc{id}%
\nolimits_{U})\circ (\theta _{2(V,U)}^{-1}\otimes \limfunc{id}%
\nolimits_{V\otimes U}) \\
&=&(e_{V}\otimes e_{U})\circ (\limfunc{id}\nolimits_{Q(V)}\otimes
c_{(Q(U),V)}\otimes \limfunc{id}\nolimits_{U})\circ \\
&&(\theta _{2(V,U)}^{-1}\otimes \limfunc{id}\nolimits_{V\otimes U}) \\
&=&e_{V\otimes U}\text{.}
\end{eqnarray*}
The first equality is the multiplicativity axiom, the second is the
interchange law and the naturality of $\theta _{2}$, the third is the
interchange law and the naturality of the braiding, the fourth is the
interchange law, the fifth is the preservation of the evaluations by $f$ and
$g$, and the last is again the multiplicativity axiom. Since we are assuming
strictness the structural isomorphisms are identities and thus $\mathcal{S}(%
\mathbf{D}\mathcal{)}$ is monoidal.

Finally $\mathcal{S}(\mathbf{D}\mathcal{)}$ is symmetric since for any pair
of objects $(V,W)$ of $\mathcal{S}(\mathbf{D})$, $c_{(V,W)}$ belongs to $%
\limfunc{Mor}(\mathcal{S}(\mathbf{D))}$. To show this we need the following
equation (writing composition as juxtaposition):
\begin{eqnarray}
&&c_{(Q(V)\otimes V,Q(W)\otimes W)}(\limfunc{id}\nolimits_{Q(V)}\otimes
c_{(Q(W),V)}\otimes \limfunc{id}\nolimits_{W})=\text{\label{eq31}} \\
&=&(((\limfunc{id}\nolimits_{Q(W)}\otimes
c_{(Q(V),W)})(c_{(Q(V),Q(W))}\otimes \limfunc{id}\nolimits_{W}))\otimes
\limfunc{id}\nolimits_{V})(\limfunc{id}\nolimits_{Q(V)\otimes Q(W)}\otimes
c_{(V,W)})\text{.}  \nonumber
\end{eqnarray}
which may be easily derived, using the commutativity of the second hexagonal
diagram of the braiding and then the commutativity of the first hexagonal
diagram, composing with the expression $\limfunc{id}\nolimits_{Q(V)}\otimes
c_{(Q(W),V)}\otimes \limfunc{id}\nolimits_{W}$, and applying the interchange
law twice.\newline
Then we have:
\begin{eqnarray*}
e_{V\otimes W} &=&(e_{V}\otimes e_{W})(\limfunc{id}\nolimits_{Q(V)}\otimes
c_{(Q(W),V)}\otimes \limfunc{id}\nolimits_{W})(\theta _{2(V,W)}^{-1}\otimes
\limfunc{id}\nolimits_{V\otimes W}) \\
&=&c_{(I,I)}(e_{V}\otimes e_{W})(\limfunc{id}\nolimits_{Q(V)}\otimes
c_{(Q(W),V)}\otimes \limfunc{id}\nolimits_{W})(\theta _{2(V,W)}^{-1}\otimes
\limfunc{id}\nolimits_{V\otimes W}) \\
&=&(e_{W}\otimes e_{V})c_{(Q(V)\otimes V,Q(W)\otimes W)}(\limfunc{id}%
\nolimits_{Q(V)}\otimes c_{(Q(W),V)}\otimes \limfunc{id}\nolimits_{W})(%
\theta _{2(V,W)}^{-1}\otimes \limfunc{id}\nolimits_{V\otimes W}) \\
&=&(e_{W}\otimes e_{V})(((\limfunc{id}\nolimits_{Q(W)}\otimes
c_{(Q(V),W)})(c_{(Q(V),Q(W))}\otimes \limfunc{id}\nolimits_{W}))\otimes
\limfunc{id}\nolimits_{V}) \\
&&(\limfunc{id}\nolimits_{Q(V)\otimes Q(W)}\otimes c_{(V,W)})(\theta
_{2(V,W)}^{-1}\otimes \limfunc{id}\nolimits_{V\otimes W}) \\
&=&(e_{W}\otimes e_{V})(\limfunc{id}\nolimits_{Q(W)}\otimes
c_{(Q(V),W)}\otimes \limfunc{id}\nolimits_{V})(c_{(Q(V),Q(W))}\otimes
\limfunc{id}\nolimits_{W\otimes V}) \\
&&(\limfunc{id}\nolimits_{Q(V)\otimes Q(W)}\otimes c_{(V,W)})(\theta
_{2(V,W)}^{-1}\otimes \limfunc{id}\nolimits_{V\otimes W}) \\
&=&(e_{W}\otimes e_{V})(\limfunc{id}\nolimits_{Q(W)}\otimes
c_{(Q(V),W)}\otimes \limfunc{id}\nolimits_{V})(c_{(Q(V),Q(W))}\otimes
c_{(V,W)}) \\
&&(\theta _{2(V,W)}^{-1}\otimes \limfunc{id}\nolimits_{V\otimes W}) \\
&=&(e_{W}\otimes e_{V})(\limfunc{id}\nolimits_{Q(W)}\otimes
c_{(Q(V),W)}\otimes \limfunc{id}\nolimits_{V})(c_{(Q(V),Q(W))}\theta
_{2(V,W)}^{-1})\otimes c_{(V,W)}) \\
&=&(e_{W}\otimes e_{V})(\limfunc{id}\nolimits_{Q(W)}\otimes
c_{(Q(V),W)}\otimes \limfunc{id}\nolimits_{V})((\theta
_{2(W,V)}^{-1}Q(c_{(V,W)}))\otimes c_{(V,W)}) \\
&=&(e_{W}\otimes e_{V})(\limfunc{id}\nolimits_{Q(W)}\otimes
c_{(Q(V),W)}\otimes \limfunc{id}\nolimits_{V})(\theta _{2(W,V)}^{-1}\otimes
\limfunc{id}\nolimits_{W\otimes V}) \\
&&(Q(c_{(V,W)})\otimes c_{(V,W)}) \\
&=&e_{W\otimes V}(Q(c_{(V,W)})\otimes c_{(V,W)})\text{.}
\end{eqnarray*}
The first equality is the multiplicativity condition of the evaluation, the
second is the insertion of the identity $c_{(I,I)}$, the third is the
naturality of the braiding, the fourth is Equation (\ref{eq31}), the fifth
is the insertion of $\limfunc{id}\nolimits_{V}$ and the interchange law, the
sixth and seventh is the interchange law, the eighth is the naturality of $%
\theta _{2}$, the ninth is the insertion of $\limfunc{id}\nolimits_{V\otimes
W}$ and the interchange law, and the last is the multiplicativity condition
of the evaluation.%
\endproof%
\medskip

{\noindent \bf Example: }%
\label{examplep33}Regarding the morphisms of $\mathcal{S}(\mathbf{D})$, we
will consider in detail the morphisms from $V$ to $V$. Let $f:V\rightarrow V$
be given by
\[
f(e_{j})=b_{ij}e_{i}\text{,}
\]
(using the summation convention over repeated indices). The condition for $f$
to be a morphism of $\mathcal{S}(\mathbf{D})$
\[
e_{V}\circ (\overline{f}\otimes f)=e_{V}
\]
corresponds to the matrix equation:
\[
\overline{B}^{T}AB=A\text{,}
\]
where $A=\left[ a_{ij}\right] $ is the matrix corresponding to the
evaluation and $B=\left[ b_{ij}\right] $ is the matrix corresponding to $f$.
Here we use
\[
\overline{f}(\overline{e}_{i})=\overline{b}_{ki}\stackrel{\_}{\cdot }%
\overline{e}_{k}
\]
and
\[
e_{V}(\overline{b}_{ki}\stackrel{\_}{\cdot }\overline{e}_{k}\otimes
b_{lj}e_{l})=\overline{b}_{ki}b_{lj}a_{kl}\text{.}
\]
In particular, if $A$ is the identity matrix then $B$ is a unitary matrix. $%
\blacktriangle \medskip $

\begin{definition}
\label{Def37}We define a symmetric, strong monoidal endofunctor $Q$ on $%
\mathcal{S}(\mathbf{D})$ as the restriction to $\mathcal{S}(\mathbf{D})$ of $%
Q$ defined on $\mathbf{D}$. The \emph{evaluation-preserving category }is the
pair\emph{\ }$(\mathcal{S}(\mathbf{D}),Q)$, where $\mathcal{S}(\mathbf{D})$
denotes the category with its symmetric, monoidal structure and $Q$ denotes
the above monoidal endofunctor on $\mathcal{S}(\mathbf{D})$.
\end{definition}

\begin{theorem}
$Q$ indeed restricts to a symmetric, strong monoidal endofunctor on $%
\mathcal{S}(\mathbf{D})$.
\end{theorem}

\proof%
The conjugation axioms b) and c) imply that $\theta _{2(V,W)}$ for any $V$, $%
W$ and $\theta _{0}$ are morphisms of $\mathcal{S}(\mathbf{D})$. Therefore
it is enough to show that $\limfunc{Mor}(\mathcal{S}(\mathbf{D}))$ is closed
under $Q$. Let $f:V\rightarrow W$ be a morphism of $\mathcal{S}(\mathbf{D})$%
. Since $f$ is an isomorphism, $Q(f)$ is an isomorphism too.

Then we have:
\begin{eqnarray*}
e_{Q(W)}(Q(Q(f))\otimes Q(f)) &=&\theta _{0}^{-1}Q(e_{W})\theta
_{2(Q(W),W)}(Q(Q(f))\otimes Q(f)) \\
&=&\theta _{0}^{-1}Q(e_{W})Q(Q(f)\otimes f)\theta _{2(Q(V),V)} \\
&=&\theta _{0}^{-1}Q(e_{V})\theta _{2(Q(V),V)} \\
&=&e_{Q(V)}\text{,}
\end{eqnarray*}
using the conjugation axiom a) in the first and last equalities and the
naturality of $\theta _{2}$ in the second equality.%
\endproof%
\medskip

The evaluation-preserving category will play an important role in
representing topological isomorphisms when we come to describing the TQFT
functor. However the actual codomain category $\mathcal{D}$ of the TQFT is
constructed directly from the algebraic starting category $\mathbf{D}$.

\begin{definition}
The category $\mathcal{D}$ is defined as follows:

\begin{enumerate}
\item[-]  the objects of $\mathcal{D}$ are all pairs of the form $(V,x)$,
where $V\in \func{Ob}(\mathbf{D})$ and $x\in G(V)$, the underlying $K$%
-module of $V$.

\item[-]  the morphisms of $\mathcal{D}$ from $(V,x)$ to $(W,y)$ are all
triples $(f,x,y)$, where $f\in \limfunc{Mor}\nolimits_{\mathbf{D}}(V,W)$, $%
x\in G(V)$ and $y\in G(W)$, such that $G(f)(x)=y$.\newline
(We will frequently write simply $f:(V,x)\rightarrow (W,y)$ for morphisms of
$\mathcal{D}$, and say that the morphisms of $\mathcal{D}$ preserve
elements.)

\item[-]  composition and identity morphisms are given by:
\begin{eqnarray*}
(g,y,z)\circ (f,x,y) &=&(g\circ f,x,z)\text{,} \\
\limfunc{id}\nolimits_{(V,x)} &=&(\limfunc{id}\nolimits_{V},x,x)\text{.}
\end{eqnarray*}
\end{enumerate}

A monoidal product, unit object and braiding on $\mathcal{D}$ are defined as
follows:\newline
monoidal product on objects:
\[
(V,x)\otimes (W,y)=(V\otimes W,x\otimes y)\text{,}
\]
\newline
monoidal product on morphisms:
\[
(f,x,y)\otimes (g,z,w)=(f\otimes g,x\otimes z,y\otimes w)\text{,}
\]
\newline
unit object:
\[
\widetilde{I}:=(I,i)\text{, where }i\text{ is the unit element of }G(I)\text{%
,}
\]
\newline
braiding:
\[
c_{((V,x),(W,y))}=(c_{(V,W)},x\otimes y,y\otimes x)\text{.}
\]
\end{definition}

It is straightforward to show:

\begin{theorem}
$(\mathcal{D},\otimes ,\widetilde{I},c)$ is a symmetric, strict monoidal
category.
\end{theorem}

\noindent For some TQFT functors only the endofunctor $Q$ on $\mathcal{S}(%
\mathbf{D})$ will play a role. For others we need a corresponding
endofunctor on $\mathcal{D}$.

\begin{definition}
\label{Def311}Given an assignment:\newline
for any $V\in \func{Ob}(\mathbf{D})$ and any $x\in G(V)$, an element $%
x_{Q(V)}$ of $G(Q(V))$, such that:

\begin{enumerate}
\item[\emph{a)}]  for any $(f,x,y)\in \limfunc{Mor}_{\mathcal{D}%
}((V,x),(W,y))$%
\[
G(Q(f))(x_{Q(V)})=y_{Q(W)}\text{,}
\]

\item[\emph{b)}]  for any $V,W\in \func{Ob}(\mathbf{D})$, $x\in G(V)$, $y\in
G(W)$%
\[
G(\theta _{2(V,W)})(x_{Q(V)}\otimes y_{Q(W)})=(x\otimes y)_{Q(V\otimes W)}%
\text{,}
\]

\item[\emph{c)}]  $G(\theta _{0})(i)=i_{Q(I)}$,
\end{enumerate}

\noindent we define a strong monoidal endofunctor $(\mathbf{Q},\mathbf{%
\theta }_{2},\mathbf{\theta }_{0})$, called an \emph{extension} of $Q$ on $%
\mathcal{D}$, as follows:\newline
on objects:
\[
\mathbf{Q}(V,x)=(Q(V),x_{Q(V)})\text{,}
\]
on morphisms:
\[
\mathbf{Q}(f,x,y)=(Q(f),x_{Q(V)},y_{Q(W)})\text{,}
\]
and the natural isomorphism $\mathbf{\theta }_{2}$ and the isomorphism $%
\mathbf{\theta }_{0}$ are given by:
\[
\mathbf{\theta }_{2((V,x),(W,y))}=(\theta _{2(V,W)},x_{Q(V)}\otimes
y_{Q(W)},(x\otimes y)_{Q(V\otimes W)})\text{,}
\]
\[
\mathbf{\theta }_{0}=(\theta _{0},i,i_{Q(I)})\text{.}
\]
\end{definition}

\begin{theorem}
$(\mathbf{Q},\mathbf{\theta }_{2},\mathbf{\theta }_{0})$ defines a
symmetric, strong monoidal endofunctor on $\mathcal{D}$.
\end{theorem}

\proof%
Conditions a), b) and c) guarantee that the relevant morphisms preserve
elements, i.e. belong to $\limfunc{Mor}(\mathcal{D})$.%
\endproof%
\medskip

{\noindent \bf Example: }%
\label{pag33}For our example the objects, morphisms and monoidal structure
of $\mathcal{D}$ are clear from the general theory. An endofunctor $\mathbf{Q%
}$ on $\mathcal{D}$ is given by the assignment:
\[
x_{Q(V)}:=k_{V}(x)\text{.}
\]
It is easy to check that the conditions for the extension apply:\newline
for $f:(V,x)\rightarrow (W,y)$,
\[
Q(f)(x_{Q(V)})=(k_{W}\circ f\circ k_{V}^{-1})(k_{V}(x))=k_{W}(y)=y_{Q(W)}
\]
and, for the structural morphisms of $\mathbf{Q}$,
\begin{eqnarray*}
\theta _{2(V,W)}(x_{\overline{V}}\otimes y_{\overline{W}}) &=&k_{V\otimes
W}(k_{V}^{-1}\otimes k_{W}^{-1})(k_{V}(x)\otimes k_{W}(y)) \\
&=&k_{V\otimes W}(x\otimes y) \\
&=&(x\otimes y)_{\overline{V\otimes W}}
\end{eqnarray*}
and
\[
\theta _{0}(1)=k_{\Bbb{C}}(\overline{1})=k_{\Bbb{C}}(1)=1_{\overline{\Bbb{C}}%
}\text{.}\blacktriangle
\]
\medskip

We are now able to give the formal definition of the two types of algebraic
category which will appear as target categories for the TQFT functors.

\begin{definition}
Let $(\mathbf{D},G,Q)$ be an algebraic starting category and $(e_{V})_{V\in
\func{Ob}(\mathbf{D})}$ be a choice of evaluations on $(\mathbf{D},G,Q)$.
The \emph{algebraic category} for these data is the pair $(\mathcal{D},Q)$,
where $\mathcal{D}$ denotes the symmetric, strict monoidal category $(%
\mathcal{D},\otimes ,\widetilde{I},c)$ and $Q$ is the monoidal endofunctor
on $\mathcal{S}(\mathbf{D})$ defined in Definition \ref{Def37}.

Given, in addition, an assignment $(V,x)\mapsto x_{Q(V)}$ satisfying
conditions a)-c) of Definition \ref{Def311}, the \emph{full algebraic
category }for these data is the pair $(\mathcal{D},\mathbf{Q})$, where $%
\mathcal{D}$ is as before and $\mathbf{Q}$ is the monoidal endofunctor on $%
\mathcal{D}$ defined in Definition \ref{Def311}.\medskip
\end{definition}

Having constructed these two types of algebraic category, corresponding to
the two types of topological category defined at the end of Section \ref
{section2}, we are now ready to proceed to the definition of TQFT functors
in the next section. However, since the example we have been using in this
section is not the most generic, we will conclude with another, more
generic, example of an algebraic category based on the category of hermitian
spaces.

\begin{example}
\label{example37}Let $\mathbf{D}$ be the category whose objects are
hermitian linear spaces, i.e. pairs $(V,h)$, where $V$ is a finite
dimensional vector space over the field $\Bbb{K}$ with an involution $j$ and
$h:V\times V\rightarrow \Bbb{K}$ is a non-degenerate hermitian form on $V$
(linear in the first argument and $j$-semilinear in the second). The
morphisms are just linear maps. The monoidal product of two objects $(V,h)$
and $(W,h^{\prime })$ is defined to be the pair $(V\otimes W,h\otimes
h^{\prime })$, where $h\otimes h^{\prime }$ is the hermitian form given by:
\[
(h\otimes h^{\prime })(x_{1}\otimes y_{1},x_{2}\otimes
y_{2})=h(x_{1},x_{2})h^{\prime }(y_{1},y_{2})\text{.}
\]
The unit object is the pair $(\Bbb{K},h_{1})$, where $h_{1}:\Bbb{K}\times
\Bbb{K}\rightarrow \Bbb{K}$ is defined by $(x,y)\mapsto xj(y)$. Clearly $%
\mathbf{D}$ is a symmetric monoidal category with the obvious braiding $c$.
The monoidal endofunctor $(Q,\theta _{2},\theta _{0})$ is defined as
follows. On objects we have $Q(V,h)=(V^{j},h^{j})$, where $V^{j}$ is the
involution vector space (i.e. the same additive group as $V$ and with scalar
multiplication $\cdot _{j}$ given by $\alpha \cdot _{j}x=j(\alpha )x$, where
scalar multiplication in $V$ is denoted by juxtaposition), and $h^{j}$ is
related to $h$ by
\[
h(x,y)=h^{j}(k_{V}(y),k_{V}(x))\text{,}
\]
where, for each $V$, $k_{V}:V\rightarrow V^{j}$ is the $j$-semilinear
`identity' map, satisfying for each $\alpha \in \Bbb{K}$ and $x\in V$, $%
k_{V}(\alpha x)=j(\alpha )\cdot _{j}x$ (see \cite[p. 63]{Evans-KawBook} for
an equivalent condition for Hilbert spaces). On morphisms $Q$ acts as
follows: for $f:(V,h)\rightarrow (W,h^{\prime })$ we have $%
Q(f):(V^{j},h^{j})\rightarrow (W^{j},h^{\prime j})$, with
\[
Q(f)=k_{W}\circ f\circ k_{V}^{-1}\text{.}
\]
The isomorphisms $\theta _{2}$ and $\theta _{0}$ are given by:
\[
\theta _{2((V,h),(W,h^{\prime }))}=k_{V\otimes W}\circ (k_{V}^{-1}\otimes
k_{W}^{-1})
\]
and
\[
\theta _{0}(\alpha )=k_{\Bbb{K}}(j(\alpha ))\text{, for all }\alpha \in \Bbb{%
K}\text{.}
\]
Then $(Q,\theta _{2},\theta _{0})$ defines a monoidal endofunctor on $%
\mathbf{D}$ in the same way as in the previous example, where $\Bbb{K}$ was $%
\Bbb{C}$ and $j$ was complex conjugation.\newline
To define the category $\mathcal{S}(\mathbf{D})$ we must first specify
evaluations $e_{(V,h)}$, for each object $(V,h)$. These may be given in
terms of the hermitian structures by:
\[
e_{(V,h)}(x\otimes y)=h(y,k_{V}^{-1}(x))\text{.}
\]
This choice of evaluations indeed satisfies the multiplication and
conjugation axioms. The multiplication axiom follows directly from the
definition of $(V,h)\otimes (W,h^{\prime })$. The conjugation axiom a) is
derived by acting on $x\otimes k_{V}(y)\in V\otimes V^{j}$:
\begin{eqnarray*}
Q(e_{(V,h)})\circ \theta _{2((V^{j},h^{j}),(V,h))}(x\otimes k_{V}(y))
&=&Q(e_{(V,h)})(k_{V^{j}\otimes V}(k_{V}(x)\otimes y)) \\
&=&k_{\Bbb{K}}(e_{(V,h)}(k_{V}(x)\otimes y)) \\
&=&k_{\Bbb{K}}(h(y,x)) \\
&=&\theta _{0}(h(x,y)) \\
&=&\theta _{0}(h^{j}(k_{V}(y),k_{V}(x))) \\
&=&\theta _{0}\circ e_{(V^{j},h^{j})}(x\otimes k_{V}(y))\text{.}
\end{eqnarray*}
The conjugation axiom b) reads:
\[
e_{Q((V,h)\otimes (W,h^{\prime }))}\circ (Q(\theta _{2((V,h),(W,h^{\prime
}))})\otimes \theta _{2((V,h),(W,h^{\prime }))})=e_{(V^{j},h^{j})\otimes
(W^{j},h^{\prime j})}\text{.}
\]
This equality is proved by applying both sides to $k_{V^{j}\otimes
W^{j}}(k_{V}(x)\otimes k_{W}(y))\otimes k_{V}(z)\otimes k_{W}(w)$ (using the
multiplicativity axiom on the right hand side) to obtain $h(x,z)h^{\prime
}(y,w)$ in both cases. Details are left to the reader.\newline
The conjugation axiom c)
\[
e_{Q(\Bbb{K},h_{1})}\circ (Q(\theta _{0})\otimes \theta _{0})=e_{(\Bbb{K}%
,h_{1})}
\]
is derived by acting on $k_{\Bbb{K}}(\alpha )\otimes \beta \in \Bbb{K}%
^{j}\otimes \Bbb{K}$:
\begin{eqnarray*}
e_{(\Bbb{K}^{j},h_{1}^{j})}\circ (Q(\theta _{0})\otimes \theta _{0})(k_{\Bbb{%
K}}(\alpha )\otimes \beta ) &=&e_{(\Bbb{K}^{j},h_{1}^{j})}(j(\alpha )\otimes
k_{\Bbb{K}}(j(\beta ))) \\
&=&h_{1}^{j}(k_{\Bbb{K}}(j(\beta )),k_{\Bbb{K}}(j(\alpha ))) \\
&=&h_{1}(j(\alpha )\otimes j(\beta )) \\
&=&h_{1}(\beta ,\alpha ) \\
&=&e_{_{(\Bbb{K},h_{1})}}(k_{\Bbb{K}}(\alpha )\otimes \beta )\text{.}
\end{eqnarray*}
The morphisms of $\mathcal{S}(\mathbf{D})$ are the isomorphisms of $\mathbf{D%
}$ which preserve the evaluations: for $f:(V,h)\rightarrow (W,h^{\prime })$
we have $e_{(W,h^{\prime })}\circ (Q(f)\otimes f)=e_{(V,h)}$. This condition
translates to the condition that $f$ preserves the hermitian structures:
\[
h^{\prime }(f(x),f(y))=h(x,y)\text{,}
\]
for all $(x,y)\in V\times V$. The construction of $\mathcal{D}$ and an
endofunctor $\mathbf{Q}$ on $\mathcal{D}$ proceed in an analogous fashion to
the previous example. As a final remark we note that this example could
easily be adapted to categories of real or complex vector spaces endowed
with inner products, giving rise to a consistent choice of evaluations in
the same way as above.$\blacktriangle $
\end{example}

\section{TQFT functors\label{section4}}

In this section we will be constructing a class of functors from the
topological category to the algebraic category (not yet the respective full
categories, which will be done in the next section). The first stage is to
define a so-called pre-TQFT functor from the category of subobjects to the
evaluation-preserving category, which can be thought of as an algebraic
representation of topological isomorphisms amongst subobjects. A TQFT
functor from the topological category to the algebraic category is a
specific type of extension of the pre-TQFT functor. Its action on
isomorphisms is determined by the pre-TQFT functor, and its action on gluing
morphisms is given essentially by evaluations on the algebraic side. To take
account of the extra structure of the endofunctors in both categories, the
definition also involves a natural isomorphism which interpolates between
the endofunctors and the pre-TQFT functor. This natural isomorphism makes it
possible for some features on the algebraic side to be richer than on the
topological side (in contrast with the usual purpose of a functor to
simplify things). Thus, for instance, in the example of $2$-dimensional
TQFT's, $P$ acts trivially on the empty set, whereas $Q$ acts non-trivially
on $\Bbb{C}$, and the natural isomorphism provides the additional
flexibility for this to be so. In our analysis of the example of $2$%
-dimensional TQFT's at the end of the section we prove a theorem which
characterizes a class of TQFT functors for this case.

We now proceed with the details of the construction.

\begin{definition}
Given a category of subobjects $(\mathcal{S}(\mathbf{C}),P)$ and an
evaluation--preserving category $(\mathcal{S}(\mathbf{D}),Q)$, as introduced
in Definitions \ref{DefN27} and \ref{Def37} respectively, a\emph{\ pre-TQFT
functor} from $(\mathcal{S}(\mathbf{C}),P)$ to $(\mathcal{S}(\mathbf{D}),Q)$
is a pair $(Z^{\prime },\eta )$, where $Z^{\prime }:\mathcal{S}(\mathbf{C}%
)\rightarrow \mathcal{S}(\mathbf{D})$ is a strict symmetric, monoidal
functor, denoted on objects by $Z^{\prime }(A)=V_{A}$ and on morphisms by $%
Z^{\prime }(\alpha )=Z_{\alpha }^{\prime }$, and $\eta $ is a monoidal
natural isomorphism
\[
\eta :(Z^{\prime }P,Z_{\pi _{2}}^{\prime },Z_{\pi _{0}}^{\prime
})\rightarrow (QZ^{\prime },\theta _{2},\theta _{0})
\]
from the monoidal functor $Z^{\prime }P$ to the monoidal functor $QZ^{\prime
}$, i.e. $\eta $ is a natural isomorphism between the functors $Z^{\prime }P$
and $QZ^{\prime }$ such that the following two diagrams, for each pair $%
(A,B) $ of objects in $\mathcal{S}(\mathbf{C)}$,%
\begin{center}
\setsqparms[1`1`1`1;1400`700]
\square[Z'P(A)\otimes Z'P(B)`QZ'(A)\otimes QZ'(B)`Z'P(A\sqcup B)`QZ'(A\sqcup B);\eta _{A}\otimes \eta _{B}`Z'_{\pi_{2}}`\theta_{2}`\eta _{A\sqcup B}]
\end{center}%
and%
\begin{center}
\settriparms[1`1`1;700]
\qtriangle[I`Z'P(E)`QZ'(E);Z'_{\pi_{0}}`\theta _{0}`\eta_{E}]
\end{center}%
commute in $\mathcal{S}(\mathbf{D})$.
\end{definition}

\begin{thm_remark}
Thus the pre-TQFT functor $Z^{\prime }$ respects the endofunctors $P$ and $Q$
only up to monoidal natural isomorphism. The need for this weakening will be
clear from the following example.\smallskip
\end{thm_remark}

{\noindent \bf Example: }%
\label{examplepg37}In our example a pre-TQFT functor $Z^{\prime }$ is given
by the following assignments on objects:
\[
V_{C_{+}}=V\text{,}\Bbb{\hspace{0.2in}}V_{C_{-}}=\overline{V}\text{,}\Bbb{%
\hspace{0.2in}}V_{\emptyset }=\Bbb{C\hspace{0.2in}}\text{and}\Bbb{\hspace{%
0.2in}}V_{A\sqcup B}=V_{A}\otimes V_{B}\text{.}
\]
The last two assignments express the strictness of $Z^{\prime }$. On
(iso)morphisms $\alpha :C_{+}\rightarrow C_{+}$, $Z_{\alpha }^{\prime
}:V\rightarrow V$ is a morphism of $\mathcal{S}(\mathbf{D})$, i.e. given by
a matrix $B_{\alpha }$ with respect to the basis $(e_{i})_{i}$ of $V$,
satisfying $\overline{B}_{\alpha }^{T}AB_{\alpha }=A$, where $A$ corresponds
to the evaluation on $V$ (see the example on page \pageref{examplep33}).
Similar statements hold for $\alpha :C_{+}\rightarrow C_{-}$, $%
C_{-}\rightarrow C_{+}$ and $C_{-}\rightarrow C_{-}$, e.g., in the first
case, $Z_{\alpha }^{\prime }$ is given by $B_{\alpha }$ satisfying $%
\overline{B}_{\alpha }^{T}\overline{A}B_{\alpha }=A$.

There are some further restrictions on the matrices $B_{\alpha }$, e.g. from
functoriality we have
\[
B_{\alpha \circ \beta }=B_{\alpha }B_{\beta }\text{\hspace{0.2in}and\hspace{%
0.2in}}B_{\limfunc{id}\nolimits_{C_{\pm }}}=I\text{.}
\]
Also, for $\alpha :C_{+}\rightarrow C_{-}$ given by $\alpha (z)=\overline{z}$%
, we have $P(\alpha )=\alpha ^{-1}$ and hence the corresponding matrix
satisfies $\overline{B}_{\alpha }B_{\alpha }=I$. However, we will not study
the many possibilities in depth, since the extension to a TQFT functor will
impose strong extra conditions on the choices for $Z_{\alpha }^{\prime }$.
As regards monoidal products and braidings, since $Z^{\prime }$ is a strict
symmetric, monoidal functor we have:
\[
Z_{\alpha \sqcup \beta }^{\prime }=Z_{\alpha }^{\prime }\otimes Z_{\beta
}^{\prime }\hspace{0.2in}\text{and}\hspace{0.2in}Z_{c_{(A,B)}}^{\prime
}=c_{(V_{A},V_{B})}\text{.}
\]

The monoidal natural isomorphism $\eta $ is given as follows. For $\emptyset
$%
\[
\eta _{\Bbb{\emptyset }}:Z^{\prime }P(\emptyset )=\Bbb{C\rightarrow }%
QZ^{\prime }(\emptyset )=\overline{\Bbb{C}}\text{,\hspace{0.2in}}\eta _{\Bbb{%
\emptyset }}=\theta _{0}\text{,}
\]
for irreducible subobjects
\[
\eta _{C_{\pm }}:Z^{\prime }P(C_{\pm })\rightarrow QZ^{\prime }(C_{\pm })%
\text{,\hspace{0.2in}}\eta _{C_{+}}=\limfunc{id}\nolimits_{\overline{V}}%
\text{\hspace{0.2in}and\hspace{0.2in}}\eta _{C_{-}}=\limfunc{id}\nolimits_{V}%
\text{,}
\]
and for monoidal products $A\sqcup B$ with $A$ and $B$ irreducible
\[
\eta _{A\sqcup B}:Z^{\prime }P(A\sqcup B)=V_{P(A)}\otimes
V_{P(B)}\rightarrow QZ^{\prime }(A\sqcup B)=\overline{V}_{A\sqcup B}
\]
is given by:
\[
\eta _{A\sqcup B}=\theta _{2(V_{A},V_{B})}\circ (\eta _{A}\otimes \eta _{B})%
\text{.}
\]
Thus $\eta $ plays the role of interpolating between the monoidal features
of ``change of orientation'' on the topological side (where $\pi _{2}$ and $%
\pi _{0}$ are trivial), and ``passing to the conjugate module'' on the
algebraic side (where $\theta _{2}$ and $\theta _{0}$ are non-trivial).$%
\blacktriangle $\medskip

We now come to our main definition, namely the definition of a TQFT functor,
which is a suitable extension of a pre-TQFT functor to a functor from $%
\mathcal{C}$ to $\mathcal{D}$.

\begin{definition}
Given a topological category $(\mathcal{C},P)$, an algebraic category $(%
\mathcal{D},Q)$ and a pre-TQFT functor $(Z^{\prime },\eta )$ from $(\mathcal{%
S}(\mathbf{C}),P)$ to $(\mathcal{S}(\mathbf{D}),Q)$, a \emph{TQFT functor}
from $(\mathcal{C},P)$ to $(\mathcal{D},Q)$ extending $(Z^{\prime },\eta )$
is a pair $(Z,\eta )$, where $\eta $ is as before and $Z$ is a symmetric,
strict monoidal functor from $\mathcal{C}$ to $\mathcal{D}$ determined by
assignments:\label{assign}

\begin{enumerate}
\item[$Z1)$]  to each object $(X,A,m)$ of $\mathcal{C}$ an object $%
Z(X,A,m)=(V_{A},Z_{(X,A,m)})$ of $\mathcal{D}$, where $V_{A}=Z^{\prime }(A)$%
, and such that, for any pair of objects $(X,A,m)$, $(X^{\prime },A^{\prime
},m^{\prime })$ of $\mathcal{C}$%
\[
Z_{(X\sqcup X^{\prime },A\sqcup A^{\prime },m\sqcup m^{\prime
})}=Z_{(X,A,m)}\otimes Z_{(X^{\prime },A^{\prime },m^{\prime })}\text{,}
\]

\item[$Z2)$]  to each isomorphism $(f,\alpha ):(X,A,m)\rightarrow (X^{\prime
},A^{\prime },m^{\prime })$ of $\mathcal{C}$ the isomorphism $Z_{(f,\alpha
)}:(V_{A},Z_{(X,A,m)})\rightarrow (V_{A^{\prime }},Z_{(X^{\prime },A^{\prime
},m^{\prime })})$, given by:
\[
Z_{(f,\alpha )}=Z_{\alpha }^{\prime }
\]
and to each gluing morphism of $\mathcal{C}$, $(f,\varphi ,\alpha
):(X,A,m)\rightarrow (X^{\prime },A^{\prime },m^{\prime })$, where $\varphi
:A_{I}\rightarrow P(A)_{J}$ and $\alpha :A_{R}\rightarrow A^{\prime }$, the
morphism $Z_{(f,\varphi ,\alpha )}:(V_{A},Z_{(X,A,m)})\rightarrow
(V_{A^{\prime }},Z_{(X^{\prime },A^{\prime },m^{\prime })})$, given by:
\[
Z_{(f,\varphi ,\alpha )}=(e_{V_{A_{J}}}\otimes Z_{\alpha }^{\prime })\circ
p^{\prime }\circ ((\eta _{A_{J}}\circ Z_{\varphi }^{\prime })\otimes
\limfunc{id}\nolimits_{V_{A_{J}}\otimes V_{A_{R}}})\circ p\text{,}
\]
where $p:V_{A}\rightarrow V_{A_{I}}\otimes V_{A_{J}}\otimes V_{A_{R}}$ and $%
p^{\prime }:Q(V_{A})_{J}\otimes V_{A_{J}}\otimes V_{A_{R}}\rightarrow
(\otimes _{j\in J}(Q(V_{A_{j}})\otimes V_{A_{j}}))\otimes V_{A_{R}}$ are the
appropriate permuting isomorphisms, and denoting $Q(V_{A})_{J}:=\otimes
_{j\in J}Q(V_{A_{j}})$, $e_{V_{A_{J}}}:=\otimes _{j\in J}e_{V_{A_{j}}}$ and $%
\eta _{A_{J}}:=\otimes _{j\in J}\eta _{A_{j}}$.

\item[$Z3)$]  to the unit object $(E,E,\limfunc{id}\nolimits_{E})$, the
object
\[
Z(E,E,\limfunc{id}\nolimits_{E})=(I,i)\text{.}
\]
\end{enumerate}
\end{definition}

\begin{thm_remark}
A TQFT functor is then given by choosing the elements $Z_{(X,A,m)}$
consistently. Apart from the monoidal restrictions in $Z1)$ and $Z3)$, every
morphism of $\mathcal{C}$ gives rise to an equation for the elements $%
Z_{(X,A,m)}$ through the condition $Z2)$. Neither $Z_{(f,\alpha )}$ nor $%
Z_{(f,\varphi ,\alpha )}$ depend explicitly on $f$, apart from the fact that
$f$ makes up a valid morphism of $\mathcal{C}$, together with $\alpha $, or $%
\varphi $ and $\alpha $. The formula for $Z_{(f,\varphi ,\alpha )}$ may be
loosely described by saying that gluing on the topological side corresponds
to evaluation on the algebraic side.
\end{thm_remark}

\begin{theorem}
A choice of assignments satisfying $Z1)$-$Z3)$ determines a symmetric,
strict monoidal functor $Z$ from $\mathcal{C}$ to $\mathcal{D}$.
\end{theorem}

\proof%
We show that $Z2)$ is functorial for the composition of two gluing
morphisms. The other cases of composition can easily be checked. To simplify
the proof we consider the special case of the following gluing morphisms $%
(X,A,m)\stackrel{(f,\varphi ,\alpha )}{\rightarrow }(X^{\prime },A^{\prime
},m^{\prime })\stackrel{(g,\psi ,\beta )}{\rightarrow }(X^{\prime \prime
},A^{\prime \prime },m^{\prime \prime })$, where $A=\sqcup _{i=1,\ldots
,5}A_{i}$, $A^{\prime }=\sqcup _{i=3,4,5}A_{i}^{\prime }$, $\varphi
:A_{1}\rightarrow P(A_{2})$, $\psi :A_{3}^{\prime }\rightarrow
P(A_{4}^{\prime })$, $\alpha =\sqcup _{i=3,4,5}\alpha _{i}$, with $\alpha
_{i}:A_{i}\rightarrow A_{i}^{\prime }$ and $\beta :A_{5}^{\prime
}\rightarrow A^{\prime \prime }$. We redefine the respective objects by $%
V_{A_{1}}=V$, $V_{A_{2}}=W$, $V_{A_{3}}=X$, $V_{A_{4}}=Y$, $V_{A_{5}}=T$, $%
V_{A_{3}^{\prime }}=X^{\prime }$, $V_{A_{4}^{\prime }}=Y^{\prime }$, $%
V_{A_{5}^{\prime }}=T^{\prime }$ and $V_{A_{5}^{\prime \prime }}=T^{\prime
\prime }$, and write $\theta $ instead of $\theta _{2}$. In the following
calculations composition is denoted by juxtaposition, and compositions are
performed before the monoidal product inside any bracketed expression:
\begin{eqnarray*}
Z_{(g,\psi ,\beta )}Z_{(f,\varphi ,\alpha )} &=&(e_{Y^{\prime }}\otimes
Z_{\beta }^{\prime })(\eta _{A_{4}^{\prime }}Z_{\psi }^{\prime }\otimes
\limfunc{id}\nolimits_{Y^{\prime }\otimes T^{\prime }})(e_{W}\otimes
Z_{\alpha }^{\prime })(\eta _{A_{2}}Z_{\varphi }^{\prime }\otimes \limfunc{id%
}\nolimits_{W\otimes X\otimes Y\otimes T}) \\
&=&(e_{Y^{\prime }}\otimes Z_{\beta }^{\prime })(\eta _{A_{4}^{\prime
}}Z_{\psi }^{\prime }\otimes \limfunc{id}\nolimits_{Y^{\prime }\otimes
T^{\prime }})(e_{W}(\eta _{A_{2}}Z_{\varphi }^{\prime }\otimes \limfunc{id}%
\nolimits_{W})\otimes Z_{\alpha }^{\prime }) \\
&=&e_{W}(\eta _{A_{2}}Z_{\varphi }^{\prime }\otimes \limfunc{id}%
\nolimits_{W})\otimes e_{Y^{\prime }}(\eta _{A_{4}^{\prime }}Z_{\psi
}^{\prime }\otimes \limfunc{id}\nolimits_{Y^{\prime }})(Z_{\alpha
_{3}}^{\prime }\otimes Z_{\alpha _{4}}^{\prime })\otimes Z_{\beta \alpha
_{5}}^{\prime } \\
&=&e_{W}(\eta _{A_{2}}Z_{\varphi }^{\prime }\otimes \limfunc{id}%
\nolimits_{W})\otimes e_{Y^{\prime }}(\eta _{A_{4}^{\prime }}Z_{\psi \alpha
_{3}}^{\prime }\otimes Z_{\alpha _{4}}^{\prime })\otimes Z_{\beta \alpha
_{5}}^{\prime } \\
&=&e_{W}(\eta _{A_{2}}Z_{\varphi }^{\prime }\otimes \limfunc{id}%
\nolimits_{W})\otimes e_{Y^{\prime }}(\eta _{A_{4}^{\prime }}Z_{P(\alpha
_{4})}^{\prime }Z_{\widetilde{\psi }}^{\prime }\otimes Z_{\alpha
_{4}}^{\prime })\otimes Z_{\beta \alpha _{5}}^{\prime } \\
&=&e_{W}(\eta _{A_{2}}Z_{\varphi }^{\prime }\otimes \limfunc{id}%
\nolimits_{W})\otimes e_{Y^{\prime }}(Q(Z_{\alpha _{4}}^{\prime })\eta
_{A_{4}}Z_{\widetilde{\psi }}^{\prime }\otimes Z_{\alpha _{4}}^{\prime
})\otimes Z_{\beta \alpha _{5}}^{\prime } \\
&=&e_{W}(\eta _{A_{2}}Z_{\varphi }^{\prime }\otimes \limfunc{id}%
\nolimits_{W})\otimes e_{Y^{\prime }}(Q(Z_{\alpha _{4}}^{\prime })\otimes
Z_{\alpha _{4}}^{\prime })(\eta _{A_{4}}Z_{\widetilde{\psi }}^{\prime
}\otimes \limfunc{id}\nolimits_{Y})\otimes Z_{\beta \alpha _{5}}^{\prime } \\
&=&e_{W}(\eta _{A_{2}}Z_{\varphi }^{\prime }\otimes \limfunc{id}%
\nolimits_{W})\otimes e_{Y}(\eta _{A_{4}}Z_{\widetilde{\psi }}^{\prime
}\otimes \limfunc{id}\nolimits_{Y})\otimes Z_{\beta \alpha _{5}}^{\prime } \\
&=&(e_{W}\otimes e_{Y}\otimes Z_{\beta \alpha _{5}}^{\prime })(\eta
_{A_{2}}Z_{\varphi }^{\prime }\otimes \limfunc{id}\nolimits_{W}\otimes \eta
_{A_{4}}Z_{\widetilde{\psi }}^{\prime }\otimes \limfunc{id}%
\nolimits_{Y\otimes T})\text{.}
\end{eqnarray*}

The first equality is the definition of the gluing morphism, the second,
third and fourth are the interchange law, the fifth is the equation $%
P(\alpha _{4})\circ \widetilde{\psi }=\psi \circ \alpha _{3}$, the sixth is
the naturality of $\eta $, and the seventh, eight and ninth are the
interchange law.

For the other side, we have:
\begin{eqnarray*}
Z_{(gf,\varphi \sqcup \widetilde{\psi },\beta \alpha _{5})} &=&(e_{W}\otimes
e_{Y}\otimes Z_{\beta \alpha _{5}}^{\prime })(\limfunc{id}%
\nolimits_{W}\otimes c_{(Q(Y),W)}\otimes \limfunc{id}\nolimits_{Y\otimes T})
\\
&&((\eta _{A_{2}}\otimes \eta _{A_{4}})Z_{\varphi \sqcup \widetilde{\psi }%
}^{\prime }\otimes \limfunc{id}\nolimits_{W\otimes Y\otimes T})(\limfunc{id}%
\nolimits_{V}\otimes c_{(W,X)}\otimes \limfunc{id}\nolimits_{Y\otimes T}) \\
&=&(e_{W}\otimes e_{Y}\otimes Z_{\beta \alpha _{5}}^{\prime })(\limfunc{id}%
\nolimits_{W}\otimes c_{(Q(Y),W)}\otimes \limfunc{id}\nolimits_{Y\otimes T})
\\
&&(\eta _{A_{2}}Z_{\varphi }^{\prime }\otimes \eta _{A_{4}}Z_{\widetilde{%
\psi }}^{\prime }\otimes \limfunc{id}\nolimits_{W\otimes Y\otimes T})(%
\limfunc{id}\nolimits_{V}\otimes c_{(W,X)}\otimes \limfunc{id}%
\nolimits_{Y\otimes T}) \\
&=&(e_{W}\otimes e_{Y}\otimes Z_{\beta \alpha _{5}}^{\prime })(\eta
_{A_{2}}Z_{\varphi }^{\prime }\otimes \limfunc{id}\nolimits_{W}\otimes \eta
_{A_{4}}Z_{\widetilde{\psi }}^{\prime }\otimes \limfunc{id}%
\nolimits_{Y\otimes T})\text{.}
\end{eqnarray*}

The first equality is the definition of the gluing morphism, the second is
the interchange law, and the third is the naturality and the symmetry of $c$%
. The functor $Z$ preserves identity morphisms, since by $Z2)$
\[
Z(\limfunc{id}\nolimits_{X},\limfunc{id}\nolimits_{A})=Z^{\prime }(\limfunc{%
id}\nolimits_{A})=\limfunc{id}\nolimits_{V_{A}}\text{.}
\]
Therefore we conclude that $Z$ is a functor.

The functor $Z$ is strict monoidal on objects from $Z1)$ and $Z3)$. To
simplify the proof that $Z$ is monoidal on morphisms, we just consider the
monoidal product of the following two gluing morphisms $(X,\sqcup
_{i=1,2,3}A_{i},m)\stackrel{(f,\varphi ,\alpha )}{\rightarrow }(X^{\prime
},A_{3}^{\prime },m^{\prime })$ and $(Y,\sqcup _{i=1,2,3}B_{i},n)\stackrel{%
(g,\psi ,\beta )}{\rightarrow }(Y^{\prime },B_{3}^{\prime },n^{\prime })$,
where $\varphi :A_{1}\rightarrow P(A_{2})$ , $\psi :B_{1}\rightarrow
P(B_{2}) $, $\alpha :A_{3}\rightarrow A_{3}^{\prime }$ and $\beta
:B_{3}\rightarrow B_{3}^{\prime }$, and we set $V_{A_{1}}=V$, $V_{A_{2}}=W$,
$V_{A_{3}}=S$, $V_{B_{1}}=X$, $V_{B_{2}}=Y$ and $V_{B_{3}}=T$.

We have:
\begin{eqnarray*}
Z_{(f,\varphi ,\alpha )\sqcup (g,\psi ,\beta )} &=&Z_{(f\sqcup g,\varphi
\sqcup \psi ,\alpha \sqcup \beta )} \\
&=&(e_{W}\otimes e_{Y}\otimes Z_{\alpha \sqcup \beta }^{\prime })p^{\prime
}((\eta _{A_{2}}\otimes \eta _{B_{2}})(Z_{\varphi }^{\prime }\otimes Z_{\psi
}^{\prime })\otimes \limfunc{id}\nolimits_{W\otimes Y\otimes S\otimes T})p \\
&=&(e_{W}\otimes e_{Y}\otimes Z_{\alpha }^{\prime }\otimes Z_{\beta
}^{\prime })p^{\prime }(\eta _{A_{2}}Z_{\varphi }^{\prime }\otimes \eta
_{B_{2}}Z_{\psi }^{\prime }\otimes \limfunc{id}\nolimits_{W\otimes Y\otimes
S\otimes T})p \\
&=&(e_{W}\otimes e_{Y}\otimes Z_{\alpha }^{\prime }\otimes Z_{\beta
}^{\prime })\widetilde{p}(\eta _{A_{2}}Z_{\varphi }^{\prime }\otimes
\limfunc{id}\nolimits_{W\otimes S}\otimes \eta _{B_{2}}Z_{\psi }^{\prime
}\otimes \limfunc{id}\nolimits_{Y\otimes T}) \\
&=&(e_{W}\otimes Z_{\alpha }^{\prime }\otimes e_{Y}\otimes Z_{\beta
}^{\prime })(\eta _{A_{2}}Z_{\varphi }^{\prime }\otimes \limfunc{id}%
\nolimits_{W\otimes S}\otimes \eta _{B_{2}}Z_{\psi }^{\prime }\otimes
\limfunc{id}\nolimits_{Y\otimes T}) \\
&=&(e_{W}\otimes Z_{\alpha }^{\prime })(\eta _{A_{2}}Z_{\varphi }^{\prime
}\otimes \limfunc{id}\nolimits_{W\otimes S})\otimes (e_{Y}\otimes Z_{\beta
}^{\prime })(\eta _{B_{2}}Z_{\psi }^{\prime }\otimes \limfunc{id}%
\nolimits_{Y\otimes T}) \\
&=&Z_{(f,\varphi ,\alpha )}\otimes Z_{(g,\psi ,\beta )}
\end{eqnarray*}
where, in the second, third and fourth equalities $p$, $p^{\prime }$ and $%
\widetilde{p}$ are the permuting isomorphisms
\begin{eqnarray*}
p &=&\limfunc{id}\nolimits_{V}\otimes (c_{(W,X)}\otimes \limfunc{id}%
\nolimits_{Y\otimes S})(\limfunc{id}\nolimits_{W}\otimes c_{(S,X\otimes
Y)})\otimes \limfunc{id}\nolimits_{T}\text{,} \\
p^{\prime } &=&\limfunc{id}\nolimits_{Q(W)}\otimes c_{(Q(Y),W)}\otimes
\limfunc{id}\nolimits_{Y\otimes S\otimes T}\text{,} \\
\widetilde{p} &=&\limfunc{id}\nolimits_{Q(W)\otimes W}\otimes
c_{(S,Q(Y)\otimes Y)}\otimes \limfunc{id}\nolimits_{T}\text{,}
\end{eqnarray*}
and we use naturality of the braiding in the third and fourth equalities,
and the interchange law in the second and fifth equalities. Finally, since $%
Z^{\prime }$ is symmetric $Z$ is also symmetric:
\begin{eqnarray*}
Z(c_{((X,A,m),(X^{\prime },A^{\prime },m^{\prime }))}) &=&Z((c_{(X,X^{\prime
})},c_{(A,A^{\prime })})) \\
&=&Z^{\prime }(c_{(A,A^{\prime })}) \\
&=&c_{(Z^{\prime }(A),Z^{\prime }(A^{\prime }))} \\
&=&c_{((Z^{\prime }(A),Z_{(X,A,m)}),(Z^{\prime }(A^{\prime }),Z_{(X^{\prime
},A^{\prime },m^{\prime })}))} \\
&=&c_{(Z_{(X,A,m)},Z_{(X^{\prime },A^{\prime },m^{\prime })})}\text{.}
\end{eqnarray*}
This completes the proof.%
\endproof%
\medskip

{\noindent \bf Example: }%
The possibility of extending a pre-TQFT functor $Z^{\prime }$ to a TQFT
functor $Z$ imposes strong restrictions on $Z^{\prime }$ as we will now
discuss.

We consider first the gluing morphism
\[
(f,\varphi ,\alpha ):(D_{-+},C_{-+},m)\rightarrow (S_{+},\emptyset
,m^{\prime })
\]
corresponding to the gluing together of two disks to make a sphere, with $%
\varphi =\limfunc{id}\nolimits_{C_{-}}$, described on page \pageref
{examplepg10}. Set $Z_{(D_{+},C_{+},m_{+})}=d_{i}e_{i}\in V$ (again, and
also in what follows, using the summation convention over repeated indices)
and $Z_{(D_{-},C_{-},m_{-})}=\delta _{j}\rule{0.05in}{0in}\overline{\cdot }%
\rule{0.05in}{0in}\overline{e}_{j}\in \overline{V}$, where $m=m_{+}\sqcup
m_{-}$. Then by $Z1)$ and $Z2)$%
\begin{eqnarray*}
Z_{(S_{+},\emptyset ,m^{\prime })} &=&Z_{(f,\varphi ,\alpha )}(\delta _{j}%
\rule{0.05in}{0in}\overline{\cdot }\rule{0.05in}{0in}\overline{e}_{j}\otimes
d_{i}e_{i}) \\
&=&e_{V}(Z_{\limfunc{id}\nolimits_{C_{-}}}^{\prime }\otimes \limfunc{id}%
\nolimits_{V})(\delta _{j}\rule{0.05in}{0in}\overline{\cdot }\rule%
{0.05in}{0in}\overline{e}_{j}\otimes d_{i}e_{i}) \\
&=&e_{V}(\delta _{j}\rule{0.05in}{0in}\overline{\cdot }\rule{0.05in}{0in}%
\overline{e}_{j}\otimes d_{i}e_{i})=\delta _{j}d_{i}a_{ji}\text{,}
\end{eqnarray*}
where $a_{ij}$ are the entries of the matrix corresponding to $e_{V}$ - see
page \pageref{examplep30}. Now suppose we glue the two disks together ``with
a twist'', i.e. with a different identification along their boundaries. The
result is still a sphere and the process is described by a gluing morphism
of the form:
\[
(g,\psi ,\beta ):(D_{-+},C_{-+},m)\rightarrow (S_{+},\emptyset ,m^{\prime })%
\text{,}
\]
where $\psi :C_{-}\rightarrow C_{-}$ is no longer $\limfunc{id}%
\nolimits_{C_{-}}$ and $g$ is chosen to be compatible with $\psi $. This
gives rise to the equation:
\begin{eqnarray*}
Z_{(S_{+},\emptyset ,m^{\prime })} &=&Z_{(g,\psi ,\beta )}(\delta _{j}\rule%
{0.05in}{0in}\overline{\cdot }\rule{0.05in}{0in}\overline{e}_{j}\otimes
d_{i}e_{i}) \\
&=&e_{V}(Z_{\psi }^{\prime }(\delta _{j}\rule{0.05in}{0in}\overline{\cdot }%
\rule{0.05in}{0in}\overline{e}_{j})\otimes d_{i}e_{i}) \\
&=&b_{kj}\delta _{j}d_{i}a_{ki}\text{,}
\end{eqnarray*}
where $B_{\psi }=[b_{ij}]$ is the matrix representing $Z_{\psi }^{\prime }$
with respect to the bases $(\overline{e}_{i})_{i}$. More generally, for the
gluing together of any two objects of the form $(X,C_{-},m_{1})$ and $%
(Y,C_{+},m_{2})$, we have, setting
\[
Z_{(X,C_{-},m_{1})}=x_{j}\rule{0.05in}{0in}\overline{\cdot }\rule%
{0.05in}{0in}\overline{e}_{j}\hspace{0.15in}\text{and}\hspace{0.15in}%
Z_{(Y,C_{+},m_{2})}=y_{i}e_{i}\text{,}
\]
the equation
\[
b_{kj}x_{j}y_{i}a_{ki}=x_{k}y_{i}a_{ki}\text{,}
\]
for any $B_{\psi }=[b_{ij}]$. Therefore we will impose the following
requirement on $Z^{\prime }$:
\[
\text{for all \hspace{0.1in}}\psi \in \limfunc{Mor}\nolimits_{\mathcal{S}(%
\mathbf{C})}(C_{-},C_{-})\text{,\hspace{0.2in}}Z_{\psi }^{\prime }=\limfunc{%
id}\nolimits_{\overline{V}}\text{,}
\]
which in turn implies, setting $\varphi =P(\psi ):C_{+}\rightarrow C_{+}$
and using $Z_{\varphi }^{\prime }=\overline{Z_{\psi }^{\prime }}$,
\[
\text{for all \hspace{0.1in}}\varphi \in \limfunc{Mor}\nolimits_{\mathcal{S}(%
\mathbf{C})}(C_{+},C_{+})\text{,\hspace{0.2in}}Z_{\varphi }^{\prime }=%
\limfunc{id}\nolimits_{V}\text{.}
\]

Before considering constraints on $Z_{\varphi }^{\prime }$ for the remaining
morphisms of $\mathcal{S}(\mathbf{C})$, we will derive a constraint on $%
e_{V} $, by considering gluing two disks together with their order swapped,
described by the gluing morphism:
\[
(h,\chi ,\alpha ):(D_{+-},C_{+-},m_{+}\sqcup m_{-})\rightarrow
(S_{+},\emptyset ,m^{\prime })\text{,}
\]
where $h(z,1)=f(z,2)$ and $h(z,2)=f(z,1)$ (see above and page \pageref
{examplepg10}) and $\chi =\limfunc{id}\nolimits_{C_{+}}$. Applying the TQFT
functor we get the equation:
\begin{eqnarray*}
Z_{(S_{+},\emptyset ,m^{\prime })} &=&Z_{(h,\chi ,\alpha
)}(d_{i}e_{i}\otimes \delta _{j}\rule{0.05in}{0in}\overline{\cdot }\rule%
{0.05in}{0in}\overline{e}_{j}) \\
&=&e_{\overline{V}}(d_{i}e_{i}\otimes \delta _{j}\rule{0.05in}{0in}\overline{%
\cdot }\rule{0.05in}{0in}\overline{e}_{j}) \\
&=&d_{i}\delta _{j}\overline{a}_{ij}\text{,}
\end{eqnarray*}
using the relation between $e_{V}$ and $e_{\overline{V}}$ described on page
\pageref{page22} in the last equality. Comparing with the previous equation
for $Z_{(S_{+},\emptyset ,m^{\prime })}$, it is natural to impose
\[
A=\overline{A}^{T}\text{,}
\]
for the matrix $A$ corresponding to the evaluation.

Now, given any $\alpha \in \limfunc{Mor}\nolimits_{\mathcal{S}(\mathbf{C}%
)}(C_{+},C_{-})$ and $\beta \in \limfunc{Mor}\nolimits_{\mathcal{S}(\mathbf{C%
})}(C_{-},C_{+})$, we have from the requirements above:
\[
Z_{\alpha }^{\prime }\circ Z_{\beta }^{\prime }=\limfunc{id}\nolimits_{%
\overline{V}}\hspace{0.15in}\text{and}\hspace{0.15in}Z_{\beta }^{\prime
}\circ Z_{\alpha }^{\prime }=\limfunc{id}\nolimits_{V}\text{.}
\]
Thus $Z_{\alpha }^{\prime }$ and $Z_{\beta }^{\prime }$ are independent of $%
\alpha $ and $\beta $, respectively, and represented by the matrices $B$ and
$B^{-1}$, respectively, with respect to the bases $(e_{i})_{i}$ and $(%
\overline{e}_{i})_{i}$ of $V$ and $\overline{V}$. Furthermore, from the
example starting on page \pageref{examplepg37}, we have the condition $%
\overline{B}B=I$, and so a natural choice is
\[
B=I\text{,}
\]
assuming that this choice satisfies the condition for $Z_{\alpha }^{\prime }$
to belong to $\limfunc{Mor}(\mathcal{S}(\mathbf{D}))$, namely $\overline{B}%
^{T}\overline{A}B=A$. Thus we must have:
\[
\overline{A}=A\text{,}
\]
whilst we already had $\overline{A}^{T}=A$, so that $A$ has to be symmetric
and real. A natural choice, which also ensures non-degeneracy of the
evaluation on $V$, is $A=I$.

Summarizing the above considerations, we will from now on restrict our
attention to TQFT functors $(Z,\eta )$ which extend the pre-TQFT functor $%
(Z^{\prime },\eta )$, where $Z^{\prime }$ satisfies:
\begin{eqnarray}
\text{for all } &&\alpha \in \limfunc{Mor}\nolimits_{\mathcal{S}(\mathbf{C}%
)}(C_{+},C_{+})\text{,\hspace{0.25in}}Z_{\alpha }^{\prime }=\limfunc{id}%
\nolimits_{V}  \nonumber \\
\text{for all } &&\alpha \in \limfunc{Mor}\nolimits_{\mathcal{S}(\mathbf{C}%
)}(C_{-},C_{-})\text{,\hspace{0.25in}}Z_{\alpha }^{\prime }=\limfunc{id}%
\nolimits_{\overline{V}}  \nonumber \\
\text{for all } &&\alpha \in \limfunc{Mor}\nolimits_{\mathcal{S}(\mathbf{C}%
)}(C_{+},C_{-})\text{,\hspace{0.25in}}Z_{\alpha }^{\prime }(e_{i})=\overline{%
e}_{i}  \label{4Zequations} \\
\text{for all } &&\alpha \in \limfunc{Mor}\nolimits_{\mathcal{S}(\mathbf{C}%
)}(C_{-},C_{+})\text{,\hspace{0.25in}}Z_{\alpha }^{\prime }(\overline{e}%
_{i})=e_{i}\text{,}  \nonumber
\end{eqnarray}
where $\eta $ (introduced in the example starting on page \pageref
{examplepg37}) satisfies:
\begin{eqnarray}
\eta _{\emptyset } &=&\theta _{0}\text{,}  \nonumber \\
\eta _{C_{+}} &=&\limfunc{id}\nolimits_{\overline{V}}\hspace{0.2in}\text{and}%
\hspace{0.2in}\eta _{C_{-}}=\limfunc{id}\nolimits_{V}\text{,}  \nonumber \\
\eta _{A\sqcup B} &=&\theta _{2(V_{A},V_{B})}\circ (\eta _{A}\otimes \eta
_{B})\text{,\label{Z2equation}}
\end{eqnarray}
and for which the evaluation in the algebraic category is given by:
\begin{equation}
e_{V}(\overline{e}_{i}\otimes e_{j})=\delta _{ij}\text{,\label{Zequation}}
\end{equation}
where $\delta _{ij}$ is the Kronecker symbol.

Having made these choices, we investigate some relations arising from
topological isomorphisms. Consider first the isomorphism
\[
(f,\alpha ):(S_{+},\emptyset ,m)\rightarrow (S_{-},\emptyset ,m^{\prime })%
\text{,}
\]
where $m$ and $m^{\prime }$ are the respective empty maps, $\alpha =\limfunc{%
id}\nolimits_{\emptyset }$, and $f$ is given by $(x,y,z)\mapsto (-x,y,z)$.
Under the TQFT functor $(f,\alpha )$ goes to the morphism of $\mathcal{D}$:
\[
Z_{(f,\alpha )}:(\Bbb{C},Z_{(S_{+},\emptyset ,m)})\rightarrow (\Bbb{C}%
,Z_{(S_{-},\emptyset ,m^{\prime })})\text{.}
\]
Since $Z_{(f,\alpha )}=Z_{\alpha }^{\prime }=\limfunc{id}\nolimits_{\Bbb{C}}$%
, this implies the equation:
\[
Z_{(S_{+},\emptyset ,m)}=Z_{(S_{-},\emptyset ,m^{\prime })}\text{.}
\]
A similar argument applied to the torus (and indeed to any manifold with
empty boundary, as we will see shortly) gives:
\[
Z_{(T_{+},\emptyset ,m)}=Z_{(T_{-},\emptyset ,m^{\prime })}\text{.}
\]

Next we will consider four objects in $\mathcal{C}$ associated with the disk
$D$:
\[
(D_{+},C_{+},m_{1})\text{,\hspace{0.15in}}(D_{-},C_{-},m_{2})\text{,\hspace{%
0.15in}}(D_{+},C_{-},m_{3})\text{,\hspace{0.15in}}(D_{-},C_{+},m_{4})\text{,}
\]
where $m_{1}(z)=m_{2}(z)=z$ and $m_{3}(z)=m_{4}(z)=\overline{z}$. Let $%
Z_{(D_{+},C_{+},m_{1})}=d_{i}e_{i}\in V$. The isomorphism
\[
(\limfunc{id}\nolimits_{D},\alpha ):(D_{+},C_{+},m_{1})\rightarrow
(D_{+},C_{-},m_{3})\text{,}
\]
where $\alpha :C_{+}\rightarrow C_{-}$ is given by $\alpha (z)=\overline{z}$%
, yields the equation:
\[
Z_{(D_{+},C_{-},m_{3})}=Z_{\alpha }^{\prime }(d_{i}e_{i})=d_{i}\rule%
{0.05in}{0in}\overline{\cdot }\rule{0.05in}{0in}\overline{e}_{i}\in
\overline{V}
\]
(using the third equation of (\ref{4Zequations}) above for $Z_{\alpha
}^{\prime }$). The reverse map described in the example on page \pageref
{examplepg8}:
\[
(r,\alpha ):(D_{+},C_{+},m_{1})\rightarrow (D_{-},C_{-},m_{2})\text{,}
\]
with $r(z)=\alpha (z)=\overline{z}$, gives rise to the equation:
\[
Z_{(D_{-},C_{-},m_{2})}=d_{i}\rule{0.05in}{0in}\overline{\cdot }\rule%
{0.05in}{0in}\overline{e}_{i}\in \overline{V}\text{.}
\]
Similarly
\[
Z_{(D_{-},C_{+},m_{4})}=d_{i}e_{i}\in V\text{.}
\]
Thus it is enough to fix $Z$ of one of the objects in order to determine $Z$
of the other three.

Also, taking different monomorphisms $m$ does not introduce anything new.
For instance, consider the object $(D_{+},C_{+},\widetilde{m})$, where $%
\widetilde{m}\neq m_{1}$. Then there exists an automorphism $\alpha
:C_{+}\rightarrow C_{+}$ such that $\widetilde{m}=m_{1}\circ \alpha $, i.e.
we have an isomorphism
\[
(\limfunc{id}\nolimits_{D},\alpha ):(D_{+},C_{+},\widetilde{m})\rightarrow
(D_{+},C_{+},m_{1})\text{.}
\]
Since $Z_{(\limfunc{id}\nolimits_{D},\alpha )}=Z_{\alpha }^{\prime }=%
\limfunc{id}\nolimits_{V}$ by the first equation of (\ref{4Zequations})
above, we have
\[
Z_{(D_{+},C_{+},\widetilde{m})}=Z_{(D_{+},C_{+},m_{1})}=d_{i}e_{i}\in V\text{%
,}
\]
for any $\widetilde{m}$.

The annulus can appear as an object of $\mathcal{C}$ in sixteen different
versions, by varying the orientation of the annulus, the orientations of the
two subobject circles and the two ways of associating the subobject circles
to the boundary circles of the annulus. These are all related by
isomorphisms, which means that after choosing $Z$ of one of them, $Z$ of all
the other combinations is determined. For instance, the isomorphism:
\[
(\limfunc{id}\nolimits_{A_{+}},c_{(C_{-},C_{+})}):(A_{+},C_{-+},m)%
\rightarrow (A_{+},C_{+-},m^{\prime })\text{,}
\]
where $m(z,1)=z=m^{\prime }(z,2)$ and $m(z,2)=2z=m^{\prime }(z,1)$, implies
the equation:
\[
Z_{(A_{+},C_{+-},m^{\prime })}=c_{(\overline{V},V)}(Z_{(A_{+},C_{-+},m)})
\]
(since $Z^{\prime }$ and $Z$ are symmetric monoidal). To pass to the other
orientation of the annulus, we have the isomorphism:
\[
(f,c_{(C_{-},C_{+})}):(A_{+},C_{-+},m)\rightarrow (A_{-},C_{+-},m^{\prime
\prime })\text{,}
\]
where $f(re^{i\theta })=(3-r)e^{i\theta }$, and $m$ is given by $%
m(z,1)=z=m^{\prime \prime }(z,1)$ and $m(z,2)=2z=m^{\prime \prime }(z,2)$,
which implies the equation:
\[
Z_{(A_{-},C_{+-},m^{\prime \prime })}=c_{(\overline{V}%
,V)}(Z_{(A_{+},C_{-+},m)})\text{.}
\]

Finally we give a couple more examples of equations arising from gluing
morphisms before giving a general result. The gluing morphism from the
annulus to the torus
\[
(g,\psi ,\beta ):(A_{+},C_{-+},m)\rightarrow (T_{+},\emptyset ,m^{\prime })%
\text{,}
\]
described in example c) on page \pageref{examplepg17}, leads to the
relation:
\[
Z_{(T_{+},\emptyset ,m^{\prime })}=Z_{(g,\psi ,\beta )}(Z_{(A_{+},C_{-+},m)})%
\text{.}
\]
Writing
\[
Z_{(A_{+},C_{-+},m)}=c_{ij}\overline{e}_{i}\otimes e_{j}\in \overline{V}%
\otimes V
\]
and using Condition $Z2)$ and Equation (\ref{Zequation}) this determines $%
Z_{(T_{+},\emptyset ,m^{\prime })}$ in terms of $Z_{(A_{+},C_{-+},m)}$:
\[
Z_{(T_{+},\emptyset ,m^{\prime })}=c_{ii}\in \Bbb{C}\text{.}
\]
The gluing together of two annuli to make an annulus, may be described by
the gluing morphism:
\[
(f,\varphi ,\alpha ):(A_{++},C_{-+-+},m\sqcup m)\rightarrow (A_{+},C_{-+},m)%
\text{,}
\]
with $m$ as above, $f$ given by
\[
f(re^{i\theta },1)=\frac{r+1}{2}e^{i\theta }\text{\hspace{0.15in}and\hspace{%
0.15in}}f(re^{i\theta },2)=\frac{r+2}{2}e^{i\theta }\text{,}
\]
$I=\left\{ 2\right\} $, $J=\left\{ 3\right\} $, $\varphi =\limfunc{id}%
\nolimits_{C_{+}}$ and $\alpha =\limfunc{id}\nolimits_{C_{-+}}$. The
corresponding morphism in $\mathcal{D}$ implies the condition:
\[
Z_{(f,\varphi ,\alpha )}(c_{ik}c_{lj}\overline{e}_{i}\otimes e_{k}\otimes
\overline{e}_{l}\otimes e_{j})=c_{ij}\overline{e}_{i}\otimes e_{j}\text{,}
\]
which, using Condition $Z2)$ and the equation $e_{\overline{V}}(e_{k}\otimes
\overline{e}_{l})=\overline{\delta }_{kl}=\delta _{kl}$ (coming from
Equation (\ref{Z2equation}) and the example on page \pageref{examplep30})
leads to the well-known constraint on the components $c_{ij}$:
\[
c_{ik}c_{kj}=c_{ij}\text{.}
\]

To conclude the example in this section we will now proceed to characterize
the TQFT functors $Z$ under the Conditions (\ref{4Zequations})-(\ref
{Zequation}) above. First of all we show how various different topological
gluing scenarios give the same algebraic result. Suppose we have a gluing
morphism
\[
(f,\varphi ,\alpha ):(X,C_{-+},m)\rightarrow (X^{\prime },\emptyset
,m^{\prime })\text{,}
\]
where $X$ is a $2$-manifold with two circles as its boundary, and $I=\left\{
1\right\} $, $J=\left\{ 2\right\} $, $\varphi =\limfunc{id}\nolimits_{C_{-}}$
and $\alpha =\limfunc{id}\nolimits_{\emptyset }$. Fixing $f$, we can change
the orientation of the boundary circles of the domain object:
\[
(f,\psi ,\alpha ):(X,C_{+-},\overline{m})\rightarrow (X^{\prime },\emptyset
,m^{\prime })\text{,}
\]
with $I=\left\{ 1\right\} $, $J=\left\{ 2\right\} $, $\psi =P(\varphi )=%
\limfunc{id}\nolimits_{C_{+}}$, $\overline{m}(z,1)=m(\overline{z},1)$ and $%
\overline{m}(z,2)=m(\overline{z},2)$, or permute the boundary components in
either orientation
\[
(f,\chi ,\alpha ):(X,C_{+-},m_{c})\rightarrow (X^{\prime },\emptyset
,m^{\prime })\text{,}
\]
where $m_{c}=m\circ c_{(C_{+},C_{-})}$, $I=\left\{ 2\right\} $, $J=\left\{
1\right\} $, $\chi =\limfunc{id}\nolimits_{C_{-}}$ and
\[
(f,\rho ,\alpha ):(X,C_{-+},\overline{m}_{c})\rightarrow (X^{\prime
},\emptyset ,m^{\prime })\text{,}
\]
where $\overline{m}_{c}=\overline{m}\circ c_{(C_{-},C_{+})}$, $I=\left\{
2\right\} $, $J=\left\{ 1\right\} $, $\rho =\limfunc{id}\nolimits_{C_{+}}$.
As we have seen above, isomorphisms between the different objects built from
$X$ imply relations. If $Z_{(X,C_{-+},m)}=x_{ij}\overline{e}_{i}\otimes
e_{j} $ then
\begin{eqnarray*}
Z_{(X,C_{+-},\overline{m})} &=&x_{ij}e_{i}\otimes \overline{e}_{j}\text{,} \\
Z_{(X,C_{+-},m_{c})} &=&x_{ij}e_{j}\otimes \overline{e}_{i}\text{,} \\
Z_{(X,C_{-+},\overline{m}_{c})} &=&x_{ij}\overline{e}_{j}\otimes e_{i}\text{.%
}
\end{eqnarray*}
However, applying the TQFT functor to the previous gluing morphisms, these
elements are all mapped to $Z_{(X^{\prime },\emptyset ,m^{\prime })}=x_{ii}$%
, since $e_{V}(\overline{e}_{i}\otimes e_{j})=\delta _{ij}=e_{\overline{V}%
}(e_{i}\otimes \overline{e}_{j})$.\newline
The same independence of the description holds for objects $X$ with more
than two boundary circles when we glue two of the boundary components
together.

We will shortly be needing the following result on decomposing surfaces. Let
$X$ be a surface of genus $g$ with $n$ boundary circles. A marking on $X$ is
a set of $3g-3+n$ non-contractible pairwise non-isotopic circles on $X$,
considered up to isotopy. Cutting along these circles decomposes $X$ into $%
2g-2+n$ pairs-of-pants. The only surfaces for which there is no such pants
decomposition are the sphere $S$, the disk $D$, the annulus $A$ and the
torus $T$. In \cite[Lemma 1.2]{Kohno} Kohno proved that any two markings of $%
X$ can be obtained from each other by a finite sequence of moves of the
following two types:

\begin{figure}[h]
\centerline{\psfig{figure=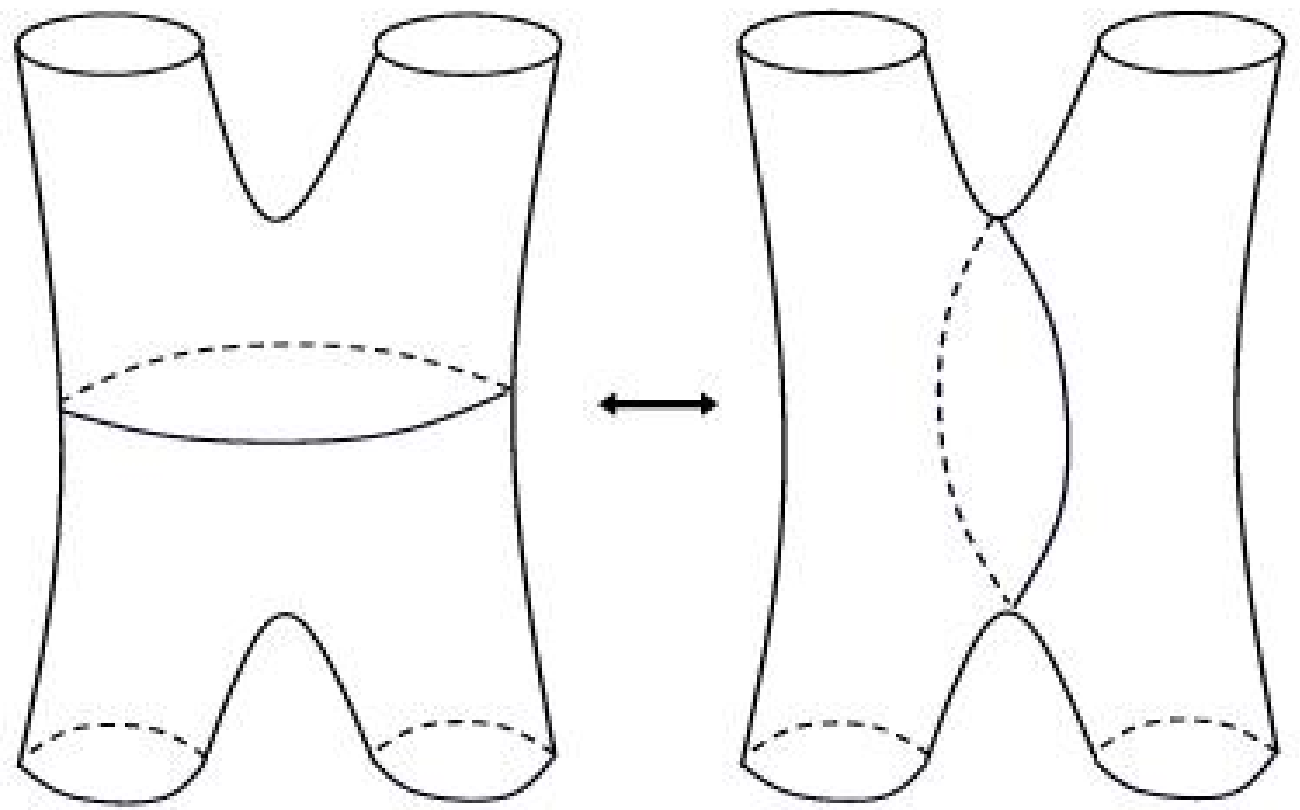,height=4cm,
}}
\caption{Type I move on markings}
\label{K1}
\end{figure}

\begin{figure}[h]
\centerline{\psfig{figure=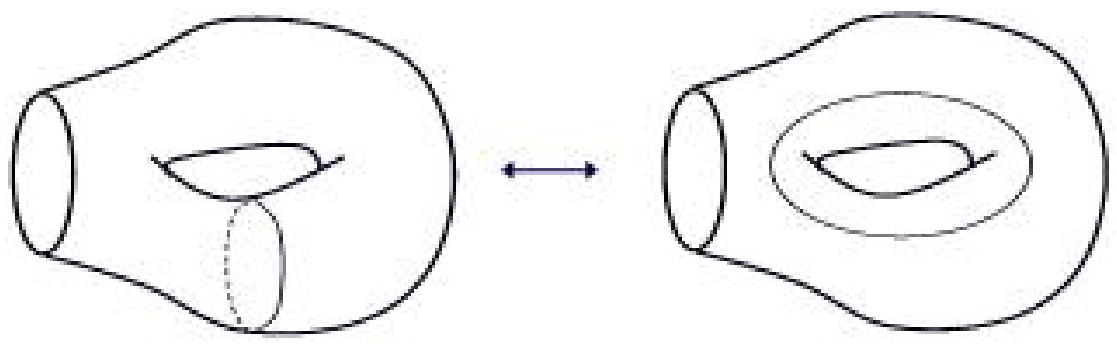,height=3cm,
}}
\caption{Type II move on markings}
\label{K2}
\end{figure}

Now we can state our general result:

\begin{theorem}
\label{thm44}TQFT functors $(Z,\eta )$, which extend $(Z^{\prime },\eta )$
for which the Conditions (\ref{4Zequations})-(\ref{Zequation}) hold, are in
1-1 correspondence with pairs of complex-valued tensors $d_{i}$, $p_{ijk}$
satisfying the relations:

\begin{enumerate}
\item  \label{Relation1}$p_{ijk}p_{klm}=p_{ilk}p_{kjm}$,

\item  \label{relation2}$p_{ijk}=p_{jik}=p_{ikj}$,

\item  \label{relation3}$d_{k}p_{kij}d_{j}=d_{i}$,

\item  \label{relation4}$d_{k}p_{kij}p_{jlm}=p_{ilm}$.
\end{enumerate}
\end{theorem}

\proof%
Given such a TQFT functor $(Z,\eta )$, we define $d_{i}$, $p_{ijk}$ by
\begin{eqnarray*}
Z_{(D_{+},C_{+},m)} &=&d_{i}e_{i}\text{,} \\
Z_{(P_{+},C_{+++},m^{\prime })} &=&p_{ijk}e_{i}\otimes e_{j}\otimes e_{k}%
\text{.}
\end{eqnarray*}
Relation \ref{Relation1} comes from looking at the two different ways of
obtaining the ``quaternion'' surface $Q$ in the type $I$ move above by
gluing two pairs-of-pants. Let us label the boundary circles of $Q$ by $i$, $%
j$, $m$ and $l$, starting from the top left hand component and going round
anticlockwise. Consider the object $(Q_{+},C_{++++},m_{ijlm})$, where the
notation $m_{ijlm}$ means that the first circle get mapped to the boundary
component labelled $i$, and so on. This object is obtained from a gluing
morphism, say $(f,\varphi ,\alpha )$, with domain $%
(P_{++},C_{++-+++},m_{ijknlm})$ and via a different gluing morphism , say $%
(g,\psi ,\beta )$, from $(P_{++},C_{-+++++},m_{knijlm})$.\newline
Now, applying the TQFT functor, we get the equations
\begin{eqnarray*}
Z_{(Q_{+},C_{++++},m_{ijlm})} &=&Z_{(f,\varphi ,\alpha
)}(p_{ijk}p_{nlm}e_{i}\otimes e_{j}\otimes \overline{e}_{k}\otimes
e_{n}\otimes e_{l}\otimes e_{m}) \\
&=&p_{ijk}p_{klm}e_{i}\otimes e_{j}\otimes e_{l}\otimes e_{m}
\end{eqnarray*}
and
\begin{eqnarray*}
Z_{(Q_{+},C_{++++},m_{ijlm})} &=&Z_{(g,\psi ,\beta )}(p_{ilk}p_{njm}%
\overline{e}_{k}\otimes e_{n}\otimes e_{i}\otimes e_{j}\otimes e_{l}\otimes
e_{m}) \\
&=&p_{ilk}p_{kjm}e_{i}\otimes e_{j}\otimes e_{l}\otimes e_{m}\text{,}
\end{eqnarray*}
and thus we have Relation \ref{Relation1}.

The first equality in Relation \ref{relation2} follows by applying $Z$ to
the isomorphism
\[
(\limfunc{id}\nolimits_{P_{+}},c_{(C_{+},C_{+})}\sqcup \limfunc{id}%
\nolimits_{C_{+}}):(P_{+},C_{+++},m_{ijk})\rightarrow
(P_{+},C_{+++},m_{jik})
\]
giving
\[
(c_{(V,V)}\otimes \limfunc{id}\nolimits_{V})(p_{ijk}e_{i}\otimes
e_{j}\otimes e_{k})=p_{jik}e_{j}\otimes e_{i}\otimes e_{k}
\]
i.e.
\[
p_{ijk}=p_{jik}\text{.}
\]
The second equality in Relation \ref{relation2} is proved in identical
fashion.

\begin{figure}[h]
\centerline{\psfig{figure=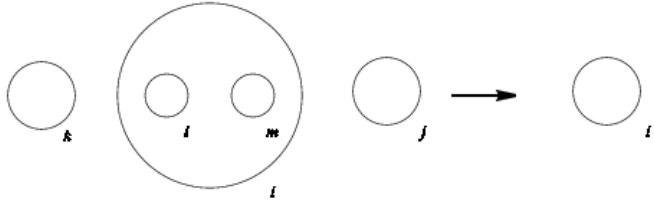,height=3cm,
}}
\caption{Relation 3}
\label{rel3}
\end{figure}

The third relation is shown by considering a gluing morphism
\[
(f,\varphi ,\alpha ):(P_{+}\sqcup D_{++},C_{-+-++},m_{lkmji})\rightarrow
(D_{+},C_{+},m_{i})\text{,}
\]
with the boundary circles labelled as in Fig. \ref{rel3} and with $I=\left\{
1,3\right\} $, $J=\left\{ 2,4\right\} $.

Applying $Z$ gives the equation:
\[
Z_{(f,\varphi ,\alpha )}(d_{k}p_{lim}d_{j}\overline{e}_{l}\otimes
e_{k}\otimes \overline{e}_{m}\otimes e_{j}\otimes e_{i})=d_{i}e_{i}
\]
i.e.
\[
d_{k}p_{kij}d_{j}=d_{i}\text{.}
\]

\begin{figure}[h]
\centerline{\psfig{figure=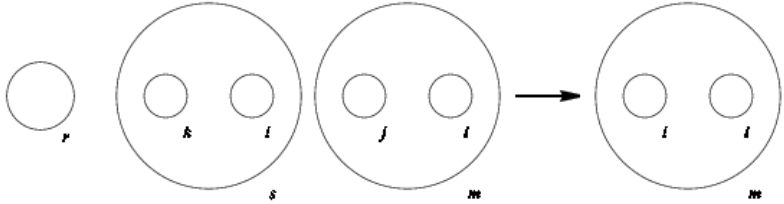,height=3cm,
}}
\caption{Relation 4}
\label{rel4}
\end{figure}

Finally the fourth relation is shown by considering a gluing morphism
\[
(f,\varphi ,\alpha ):(D_{+}\sqcup P_{++},C_{-+-++++},m_{krjsilm})\rightarrow
(P_{+},C_{+++},m)\text{,}
\]
with the boundary circles labelled as in Fig. \ref{rel4} and with $I=\left\{
1,3\right\} $, $J=\left\{ 2,4\right\} $.

Applying $Z$ gives:
\[
Z_{(f,\varphi ,\alpha )}(d_{r}p_{kis}p_{jlm}\overline{e}_{k}\otimes
e_{r}\otimes \overline{e}_{j}\otimes e_{s}\otimes e_{i}\otimes e_{l}\otimes
e_{m})=p_{ilm}e_{i}\otimes e_{l}\otimes e_{m}
\]
i.e.
\[
d_{k}p_{kij}p_{jlm}=p_{ilm}\text{.}
\]

Conversely, given tensors $d_{i}$, $p_{ijk}$, satisfying Relations 1-4, we
define a TQFT functor $Z$ as follows. Any object of $\mathcal{C}$ is of the
form
\[
(X_{\pm },C_{\pm \pm \pm \cdots },m_{ijk\cdots })\text{,}
\]
where the boundary circles have been marked $i$, $j$, $k$, $\ldots $ etc..
Set
\[
Z_{(X_{\pm },C_{\pm \pm \pm \cdots },m_{ijk\cdots })}=x_{ijk\cdots }\stackrel%
{\tiny (\_)}{e}_{i}\otimes \stackrel{\tiny (\_)}{e}_{j}\otimes \stackrel%
{\tiny (\_)}{e}_{k}\cdots \text{,}
\]
where $\stackrel{(\_)}{e}_{i}$ denotes $e_{i}$ or $\overline{e}_{i}$
depending on the orientation of the corresponding subobject circle $C_{\pm }$%
. The tensors $x_{ijk\cdots }$ are given as follows:
\[
\text{%
\begin{tabular}{ll}
empty set $\emptyset $ & $1$ \\
sphere $S$ & $d_{i}d_{i}$ \\
disk $D$ & $d_{i}$ \\
annulus $A$ & $d_{k}p_{kij}$ \\
torus $T$ & $d_{k}p_{kii}$%
\end{tabular}
}
\]
and for any other connected object $X$, by taking a pants decomposition of $%
X $, labelling all circles in the decomposition, assigning a tensor to each
labelled pair-of-pants and summing over all repeated indices (circles which
are not in the boundary of $X$). This is well-defined because of the Kohno
result above. For objects with more than one connected component we multiply
the tensors associated to each connected component.

\begin{figure}[h]
\centerline{\psfig{figure=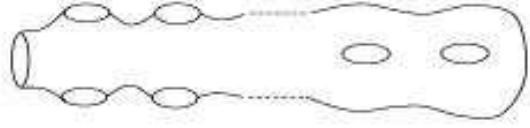,height=3cm,
}}
\caption{Reflection of $X$ in the
horizontal plane}
\label{refl}
\end{figure}

It remains to show that these assignments are consistent with the $Z2)$
axiom for $Z$. For isomorphisms which change a subobject component from $%
C_{+}$ to $C_{-}$, or vice-versa, such as $(f,\alpha
):(D_{+},C_{+},m)\rightarrow (D_{+},C_{-},\overline{m})$, where $m(z)=z$, $%
\overline{m}(z)=\overline{z}$, $f=\limfunc{id}\nolimits_{D_{+}}$ and $\alpha
(z)=\overline{z}$, the components of the tensor are unchanged (since $%
Z_{\alpha }^{\prime }(e_{i})=\overline{e}_{i}$). For isomorphisms which are
the identity on $X$ and permute the boundary components, the consistency
condition is the complete symmetry of the tensor $x_{ijk\cdots }$ under
interchanges of indices. This is clear for the annulus and the pair-of-pants
because of Relation \ref{relation2}, and for the general case one can show
symmetry under the interchange of any pair of indices corresponding to
boundary circles in the same connected component of $X$ by repeated
application of Relations \ref{Relation1} and \ref{relation2}, until the
indices both belong to the same pants tensor in the expression. For any $X$
one can construct an isomorphism $(X,A,m)\rightarrow (P(X),A^{\prime
},m^{\prime })$ by reflecting in a plane of symmetry (see Fig. \ref{refl})
so that the tensors for $X_{+}$ and $X_{-}$ are the same.

To show consistency with Condition $Z2)$ for gluing morphisms, it is enough
to consider gluing morphisms which glue two boundary circles together, since
any gluing morphism involving more than two circles can be written as a
composition of such gluing morphisms. We will consider the six combinations
of gluing together two objects, both of which are the disk, the annulus or $%
X $, where $X$ denotes a surface with non-empty boundary admitting a pants
decomposition.

\begin{enumerate}
\item  $disk+disk\rightarrow sphere$\newline
The consistency condition is an identity:
\[
d_{i}d_{i}=d_{i}d_{i}\text{.}
\]

\item  $disk+annulus\rightarrow disk$\newline
The consistency condition is
\[
d_{j}d_{k}p_{kij}=d_{i}\text{,}
\]
which holds because of Relation \ref{relation3}.

\item  $annulus+annulus\rightarrow annulus$\newline
The consistency condition is
\[
d_{k}p_{kij}d_{l}p_{ljm}=d_{k}p_{kim}\text{,}
\]
which holds because of Relation \ref{relation4}.

\item  $disk+X_{g,n}\rightarrow X_{g,n-1}$, where the indices $g$ and $n$
denote the genus and the number of boundary components, respectively.\newline
If $X=X_{0,3}=P$, then $X_{0,2}=A$ and the consistency condition is an
identity:
\[
d_{k}p_{kij}=d_{k}p_{kij}\text{.}
\]
Likewise, if $X$ is $X_{1,1}$, then $X_{1,0}=T$ and we have an identity:
\[
d_{k}p_{kii}=d_{k}p_{kii}\text{.}
\]
Otherwise $X$ has a pants decomposition into two or more pairs-of-pants and
consistency follows from Relation \ref{relation4} (see Fig. \ref{diskandx}).%

\begin{figure}[h]
\centerline{\psfig{figure=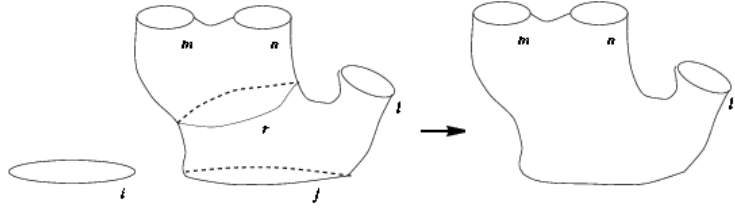,height=3cm,
}}
\caption{Gluing a disk and $X$}
\label{diskandx}
\end{figure}

\item  $annulus+X_{g,n}\rightarrow X_{g,n\text{.}}$\newline
Consistency again follows from Relation \ref{relation4} since suppose the
annulus is attached to a boundary circle labelled $j$ of $X_{g,n}$,
belonging to a pair-of-pants labelled with $j$, $l$ and $m$ in some chosen
pants decomposition. Then
\[
d_{k}p_{kij}p_{jlm}=p_{ilm}\text{.}
\]

\item  $X_{g,n}+X_{g^{\prime },n^{\prime }}^{\prime }\rightarrow
X_{g+g^{\prime },n+n^{\prime }-2}^{\prime \prime }$.\newline
Here consistency is immediate since pants decompositions of $X$ and $%
X^{\prime }$ induce a pants decomposition of $X^{\prime \prime }$ and the
formulae for the corresponding tensors are clearly compatible.\endproof
\end{enumerate}

\section{Hermitian TQFT's\label{section5}}

The TQFT functors of the previous section did not involve the full
topological and algebraic categories, i.e. in terms of the internal
structure on the categories, they only took account of the endofunctors $P$
and $Q$ defined on $\mathcal{S}(\mathbf{C})$ and $\mathcal{S}(\mathbf{D})$,
respectively. In this final section we will define TQFT functors from a full
topological category to a full algebraic category, which involve extended
endofunctors $\mathbf{P}$ and $\mathbf{Q}$ defined on $\mathcal{C}$ and $%
\mathcal{D}$, respectively. This is achieved by extending the natural
isomorphism $\eta :Z^{\prime }P\rightarrow QZ^{\prime }$ to a natural
isomorphism $\mathbf{\eta }:Z\mathbf{P}\rightarrow \mathbf{Q}Z$.

Indeed we have seen in Section \ref{section2} that the monoidal endofunctor $%
P$ defined on $\mathbf{C}$ and $\mathcal{S}(\mathbf{C})$ extends to a
monoidal endofunctor $\mathbf{P}$ on $\mathcal{C}$ (Definition \ref{Def218}%
), and in Section \ref{section3} that, under certain conditions, the
monoidal endofunctor $Q$ defined on $\mathbf{D}$ and $\mathcal{S}(\mathbf{D}%
) $ extends to a monoidal endofunctor $\mathbf{Q}$ on $\mathcal{D}$
(Definition \ref{Def311}). When $Q$ extends in this way, we can consider
extending the natural isomorphism $\eta :Z^{\prime }P\rightarrow QZ^{\prime
} $ to a natural isomorphism $\mathbf{\eta }:Z\mathbf{P}\rightarrow \mathbf{Q%
}Z $.

\begin{definition}
Given an extension $\mathbf{Q}$ of $Q$, the \emph{extension} of $\eta $
corresponding to $\mathbf{Q}$ (if it exists) is the monoidal natural
isomorphism $\mathbf{\eta }:Z\mathbf{P}\rightarrow \mathbf{Q}Z$ satisfying,
for every $(X,A,m)\in \func{Ob}(\mathcal{C})$,
\[
\mathbf{\eta }_{(X,A,m)}:=\eta _{A}\text{.}
\]
\end{definition}

For $\mathbf{\eta }_{(X,A,m)}$ to belong to $\limfunc{Mor}(\mathcal{D})$, we
have the necessary condition,
\begin{equation}
\eta _{A}(Z_{\mathbf{P}(X,A,m)})=(Z_{(X,A,m)})_{Q(V_{A})}\text{.\label%
{eqHermit}}
\end{equation}
This condition is also sufficient:

\begin{theorem}
If $\eta _{A}(Z_{\mathbf{P}(X,A,m)})=(Z_{(X,A,m)})_{Q(V_{A})}$ holds for
every $(X,A,m)\in \func{Ob}(\mathcal{C})$, then $\mathbf{\eta }$ as defined
above is a monoidal natural isomorphism from $Z\mathbf{P}$ to $\mathbf{Q}Z$.
\end{theorem}

\proof%
For isomorphisms $(f,\alpha ):(X,A,m)\rightarrow (X^{\prime },A^{\prime
},m^{\prime })$ the naturality of $\mathbf{\eta }$ follows from the
naturality of $\eta :Z^{\prime }P\rightarrow QZ^{\prime }$ for $\alpha $,
since $Z_{\mathbf{P}(f,\alpha )}$ and $\mathbf{Q}(Z_{(f,\alpha )})$ are
morphisms of $\mathcal{D}$, i.e. preserve elements, and $\mathbf{\eta }%
_{(X,A,m)}$, $\mathbf{\eta }_{(X^{\prime },A^{\prime },m^{\prime })}$
preserve elements by Equation (\ref{eqHermit}).

For gluing morphisms $(f,\varphi ,\alpha ):(X,A,m)\rightarrow (X^{\prime
},A^{\prime },m^{\prime })$ the naturality of $\mathbf{\eta }$ is given by
\[
\mathbf{Q}(Z_{(f,\varphi ,\alpha )})\circ \mathbf{\eta }_{(X,A,m)}=\mathbf{%
\eta }_{(X^{\prime },A^{\prime },m^{\prime })}\circ Z_{P(f,\varphi ,\alpha )}%
\text{,}
\]
where again all morphisms preserve elements. For simplicity we will take $%
A=\sqcup _{i=1,2,3}A_{i}$, $\varphi :A_{1}\rightarrow P(A_{2})$ and $\alpha
:A_{3}\rightarrow A^{\prime }$, set $V_{A_{1}}=V$, $V_{A_{2}}=W$, $%
V_{A_{3}}=Y$, $V_{A^{\prime }}=Y^{\prime }$, $V_{P(A_{2})}=W_{P}$ and $%
V_{P(A_{3})}=Y_{P}$, and denote $\pi _{2}$, $\theta _{2}$ by $\pi $ and $%
\theta $, respectively. We prove the equivalent equation:
\[
\mathbf{Q}(Z_{(f,\varphi ,\alpha )})\circ \mathbf{\eta }_{(X,A,m)}\circ
Z_{\pi _{(A_{1},A_{2},A_{3})}}^{\prime }=\mathbf{\eta }_{(X^{\prime
},A^{\prime },m^{\prime })}\circ Z_{P(f,\varphi ,\alpha )}\circ Z_{\pi
_{(A_{1},A_{2},A_{3})}}^{\prime }\text{,}
\]
where
\[
Z_{\pi _{(A_{1},A_{2},A_{3})}}^{\prime }=Z_{\pi _{(A_{1}\sqcup
A_{2},A_{3})}}^{\prime }\circ (Z_{\pi _{(A_{1},A_{2})}}^{\prime }\otimes
\limfunc{id}\nolimits_{Y_{P}})\text{.}
\]
We have (writing composition as juxtaposition and denoting $Z_{\pi
_{(A_{1},A_{2},A_{3})}}^{\prime }$ by $Z_{\pi }^{\prime }$):
\begin{eqnarray*}
\mathbf{Q}(Z_{(f,\varphi ,\alpha )})\mathbf{\eta }_{(X,A,m)}Z_{\pi }^{\prime
} &=&Q((e_{W}\otimes Z_{\alpha }^{\prime })(\eta _{A_{2}}Z_{\varphi
}^{\prime }\otimes \limfunc{id}\nolimits_{W\otimes Y}))\eta _{A_{1}\sqcup
A_{2}\sqcup A_{3}}Z_{\pi }^{\prime } \\
&=&Q((e_{W}\otimes Z_{\alpha }^{\prime })(\eta _{A_{2}}Z_{\varphi }^{\prime
}\otimes \limfunc{id}\nolimits_{W\otimes Y}))\theta _{(V\otimes W,Y)}(\theta
_{(V,W)}\otimes \limfunc{id}\nolimits_{Q(Y)}) \\
&&(\eta _{A_{1}}\otimes \eta _{A_{2}}\otimes \eta _{A_{3}}) \\
&=&Q(e_{W}\otimes Z_{\alpha }^{\prime })Q(\eta _{A_{2}}Z_{\varphi }^{\prime
}\otimes \limfunc{id}\nolimits_{W\otimes Y}))\theta _{(V\otimes W,Y)} \\
&&(\theta _{(V,W)}(\eta _{A_{1}}\otimes \eta _{A_{2}})\otimes \eta _{A_{3}})
\\
&=&Q(e_{W}\otimes Z_{\alpha }^{\prime })\theta _{(Q(W)\otimes W,Y)}(Q(\eta
_{A_{2}}Z_{\varphi }^{\prime }\otimes \limfunc{id}\nolimits_{W}\otimes
\limfunc{id}\nolimits_{Q(Y)}) \\
&&(\theta _{(V,W)}(\eta _{A_{1}}\otimes \eta _{A_{2}})\otimes \eta _{A_{3}})
\\
&=&\theta _{(I,Y^{\prime })}(Q(e_{W})\otimes Q(Z_{\alpha }^{\prime
}))(Q(\eta _{A_{2}}Z_{\varphi }^{\prime }\otimes \limfunc{id}%
\nolimits_{W}\otimes \limfunc{id}\nolimits_{Q(Y)}) \\
&&(\theta _{(V,W)}(\eta _{A_{1}}\otimes \eta _{A_{2}})\otimes \eta _{A_{3}})
\\
&=&\theta _{(I,Y^{\prime })}(Q(e_{W}(\eta _{A_{2}}Z_{\varphi }^{\prime
}\otimes \limfunc{id}\nolimits_{W}))\otimes Q(Z_{\alpha }^{\prime })) \\
&&(\theta _{(V,W)}(\eta _{A_{1}}\otimes \eta _{A_{2}})\otimes \eta _{A_{3}})
\\
&=&\theta _{(I,Y^{\prime })}(Q(e_{W}(\eta _{A_{2}}Z_{\varphi }^{\prime
}\otimes \limfunc{id}\nolimits_{W}))\theta _{(V,W)}(\eta _{A_{1}}\otimes
\eta _{A_{2}})\otimes Q(Z_{\alpha }^{\prime })\eta _{A_{3}} \\
&=&\theta _{(I,Y^{\prime })}(Q(e_{W}(\eta _{A_{2}}Z_{\varphi }^{\prime
}\otimes \limfunc{id}\nolimits_{W}))\theta _{(V,W)}(\eta _{A_{1}}\otimes
\eta _{A_{2}})\otimes \eta _{A^{\prime }}Z_{P(\alpha )}^{\prime } \\
&=&\theta _{(I,Y^{\prime })}(\theta _{0}e_{W_{P}}(\eta
_{P(A_{2})}Z_{P(\varphi )}^{\prime }\otimes \limfunc{id}\nolimits_{W_{P}})%
\otimes \eta _{A^{\prime }}Z_{P(\alpha )}^{\prime } \\
&=&\theta _{(I,Y^{\prime })}(\theta _{0}\otimes \eta _{A^{\prime
}})(e_{W_{P}}(\eta _{P(A_{2})}Z_{P(\varphi )}^{\prime }\otimes \limfunc{id}%
\nolimits_{W_{P}})\otimes Z_{P(\alpha )}^{\prime }) \\
&=&\theta _{(I,Y^{\prime })}(\theta _{0}\otimes \eta _{A^{\prime
}})(e_{W_{P}}\otimes Z_{P(\alpha )}^{\prime })(\eta _{P(A_{2})}Z_{P(\varphi
)}^{\prime }\otimes \limfunc{id}\nolimits_{W_{P}\otimes Y_{P}}) \\
&=&\eta _{A^{\prime }}(e_{W_{P}}\otimes Z_{P(\alpha )}^{\prime })(\eta
_{P(A_{2})}Z_{P(\varphi )}^{\prime }\otimes \limfunc{id}\nolimits_{W_{P}%
\otimes Y_{P}}) \\
&=&\mathbf{\eta }_{(X^{\prime },A^{\prime },m^{\prime })}Z_{P(f,\varphi
,\alpha )}\text{.}
\end{eqnarray*}
The first and the last equalities are definitions, the second and ninth are
shown in Lemma \ref{lemma3}, the third, sixth, seventh, tenth and eleventh
are the interchange law, the fourth and fifth are the naturality of $\theta $%
, the eighth is the naturality of $\eta $, and the twelfth is part of the
definition of $Q$ being a monoidal functor. Finally, since $\eta :Z^{\prime
}P\rightarrow QZ^{\prime }$ is monoidal we also have the corresponding
equations for $\mathbf{\eta }:Z\mathbf{P}\rightarrow \mathbf{Q}Z$:
\[
\mathbf{\theta }_{(Z(X,A,m),Z(X^{\prime },A^{\prime },m^{\prime }))}\circ (%
\mathbf{\eta }_{(X,A,m)}\otimes \mathbf{\eta }_{(X^{\prime },A^{\prime
},m^{\prime })})=\mathbf{\eta }_{(X,A,m)\sqcup (X^{\prime },A^{\prime
},m^{\prime })}\circ Z_{\pi _{((X,A,m),(X^{\prime },A^{\prime },m^{\prime
}))}}
\]
and the equation
\[
\mathbf{\eta }_{(E,E,\limfunc{id}\nolimits_{E})}Z_{\mathbf{\pi }_{0}}=%
\mathbf{\theta }_{0}
\]
in which all the morphisms preserve elements.%
\endproof%

We now prove the required lemma:

\begin{lemma}
\label{lemma3}Under the previous hypotheses, we have:
\[
Q(e_{W}(\eta _{A_{2}}Z_{\varphi }^{\prime }\otimes \limfunc{id}%
\nolimits_{W}))\theta _{(V,W)}(\eta _{A_{1}}\otimes \eta _{A_{2}})=\theta
_{0}e_{W_{P}}(\eta _{P(A_{2})}Z_{P(\varphi )}^{\prime }\otimes \limfunc{id}%
\nolimits_{W_{P}})
\]
and
\[
\eta _{A_{1}\sqcup A_{2}\sqcup A_{3}}Z_{\pi _{(A_{1},A_{2},A_{3})}}^{\prime
}=\theta _{(V\otimes W,Y)}(\theta _{(V,W)}\otimes \limfunc{id}%
\nolimits_{Q(Y)})(\eta _{A_{1}}\otimes \eta _{A_{2}}\otimes \eta _{A_{3}})%
\text{.}
\]
\end{lemma}

\proof%
To show the first equation we compose both sides with $(\eta
_{P(A_{2})}Z_{P(\varphi )}^{\prime })^{-1}\otimes \limfunc{id}%
\nolimits_{W_{P}}$. Then we have:
\begin{eqnarray*}
\theta _{0}e_{W_{P}} &=&\theta _{0}e_{Q(W)}(Q(\eta _{A_{2}})\otimes \eta
_{A_{2}}) \\
&=&Q(e_{W})\theta _{(Q(W),W)}(Q(\eta _{A_{2}})\otimes \eta _{A_{2}}) \\
&=&Q(e_{W})\theta _{(Q(W),W)}(Q(\eta _{A_{2}}Z_{\varphi }^{\prime })\otimes
\limfunc{id}\nolimits_{W})(Q(Z_{\varphi ^{-1}}^{\prime })\otimes \eta
_{A_{2}}) \\
&=&Q(e_{W}(\eta _{A_{2}}Z_{\varphi }^{\prime }\otimes \limfunc{id}%
\nolimits_{W}))\theta _{(V,W)}(Q(Z_{\varphi ^{-1}}^{\prime })\otimes \eta
_{A_{2}}) \\
&=&Q(e_{W}(\eta _{A_{2}}Z_{\varphi }^{\prime }\otimes \limfunc{id}%
\nolimits_{W}))\theta _{(V,W)}(\eta _{A_{1}}(\eta _{P(A_{2})}Z_{P(\varphi
)}^{\prime })^{-1}\otimes \eta _{A_{2}}) \\
&=&Q(e_{W}(\eta _{A_{2}}Z_{\varphi }^{\prime }\otimes \limfunc{id}%
\nolimits_{W}))\theta _{(V,W)}(\eta _{A_{1}}\otimes \eta _{A_{2}})(\eta
_{P(A_{2})}Z_{P(\varphi )}^{\prime })^{-1}\otimes \limfunc{id}%
\nolimits_{W_{P}})\text{.}
\end{eqnarray*}
The first is the fact that $\eta _{A_{2}}$ is a morphism of $\mathcal{S}(%
\mathbf{D})$, the second is the conjugation axiom of the evaluation map, the
third is the interchange law, the fourth is the naturality of $\theta $, the
fifth is the naturality of $\eta $ and the sixth is the interchange law.%
\newline
Now we prove the second equation:
\begin{eqnarray*}
\eta _{A_{1}\sqcup A_{2}\sqcup A_{3}}Z_{\pi _{(A_{1},A_{2},A_{3})}}^{\prime
} &=&\theta _{(V\otimes W,Y)}(\eta _{A_{1}\sqcup A_{2}}\otimes \eta
_{A_{3}})(Z_{\pi _{(A_{1},A_{2})}}^{\prime }\otimes \limfunc{id}%
\nolimits_{Y_{P}}) \\
&=&\theta _{(V\otimes W,Y)}(\eta _{A_{1}\sqcup A_{2}}Z_{\pi
_{(A_{1},A_{2})}}^{\prime }\otimes \eta _{A_{3}}) \\
&=&\theta _{(V\otimes W,Y)}(\theta _{(V,W)}(\eta _{A_{1}}\otimes \eta
_{A_{2}})\otimes \eta _{A_{3}}) \\
&=&\theta _{(V\otimes W,Y)}(\theta _{(V,W)}\otimes \limfunc{id}%
\nolimits_{Q(Y)})(\eta _{A_{1}}\otimes \eta _{A_{2}}\otimes \eta _{A_{3}})
\end{eqnarray*}
using the monoidal property of $\eta $ in the first and third equalities and
the interchange law in the second and fourth equalities.%
\endproof%
\medskip

{\noindent \bf Example: }%
Recall from Section \ref{section3} that in our example an extension of $Q$
to $\mathcal{D}$ is given by:
\[
\mathbf{Q}(V,x)=(\overline{V},k_{V}(x))\text{,}
\]
and a corresponding condition on morphisms. Consider an object $(X,C_{+},m)$
of $\mathcal{C}$, with a single boundary component. Then the Condition (\ref
{eqHermit}) for this object is:
\[
\eta _{C_{+}}(Z_{(P(X),C_{-},P(m))})=k_{V}(Z_{(X,C_{+},m)})\text{.}
\]
Now $\eta _{C_{+}}=\limfunc{id}\nolimits_{\overline{V}}$, so imposing the
Conditions (\ref{4Zequations})-(\ref{Zequation}) as in Section \ref{section4}%
, and bearing in mind the relation between $Z_{(X,C_{+},m)}$ and $%
Z_{(P(X),C_{-},P(m))}$ (see the discussion relating to Fig. \ref{refl}), we
have, setting $Z_{(X,C_{+},m)}=x_{i}e_{i}$,
\[
x_{i}\rule{0.05in}{0in}\overline{\cdot }\rule{0.05in}{0in}%
e_{i}=k_{V}(x_{i}e_{i})=\overline{x}_{i}\rule{0.05in}{0in}\overline{\cdot }%
\rule{0.05in}{0in}e_{i}\text{,}
\]
i.e. the tensor $x_{i}$ is real.

For an object with empty boundary $(X,\emptyset ,m)$ the Condition (\ref
{eqHermit}) is:
\[
\eta _{\emptyset }(Z_{(P(X),\emptyset ,P(m))})=k_{\Bbb{C}}(Z_{(X,\emptyset
,m)})\text{.}
\]
Setting $Z_{(X,\emptyset ,m)}=x$ we have $Z_{(P(X),\emptyset ,P(m))}=x$
also, and using $\eta _{\emptyset }=\theta _{0}$ we get:
\[
\theta _{0}(x)=k_{\Bbb{C}}(x)
\]
i.e.
\[
k_{\Bbb{C}}(\overline{x})=k_{\Bbb{C}}(x)
\]
so that $x$ has to be real.

Finally, we will consider an object $(X,C_{-+},m)$ with two boundary
components (but the reasoning can be extended to any number of components).
The Condition (\ref{eqHermit}) for this object reads:
\[
\eta _{C_{-+}}(Z_{(P(X),C_{+-},P(m))})=(Z_{(X,C_{-+},m)})_{Q(\overline{V}%
\otimes V)}\text{.}
\]
Setting $Z_{(X,C_{-+},m)}=x_{ij}\overline{e}_{i}\otimes e_{j}$ we have $%
Z_{(P(X),C_{+-},P(m))}=x_{ij}e_{i}\otimes \overline{e}_{j}$. On the left
hand side we use:
\[
\eta _{C_{-+}}=\theta _{2(\overline{V},V)}(\eta _{C_{-}}\otimes \eta
_{C_{+}})=\theta _{2(\overline{V},V)}
\]
and on the right hand side:
\begin{eqnarray*}
(x_{ij}\overline{e}_{i}\otimes e_{j})_{Q(\overline{V}\otimes V)} &=&((x_{ij}%
\rule{0.05in}{0in}\overline{\cdot }\rule{0.05in}{0in}\overline{e}%
_{i})\otimes e_{j})_{Q(\overline{V}\otimes V)} \\
&=&\theta _{2(\overline{V},V)}((x_{ij}\rule{0.05in}{0in}\overline{\cdot }%
\rule{0.05in}{0in}\overline{e}_{i})_{Q(\overline{V})}\otimes (e_{j})_{Q(V)})
\\
&=&\theta _{2(\overline{V},V)}((\overline{x}_{ij}e_{i})\otimes \overline{e}%
_{j}) \\
&=&\theta _{2(\overline{V},V)}(\overline{x}_{ij}e_{i}\otimes \overline{e}%
_{j})\text{.}
\end{eqnarray*}
Thus we get
\[
\theta _{2(\overline{V},V)}(x_{ij}e_{i}\otimes \overline{e}_{j})=\theta _{2(%
\overline{V},V)}(\overline{x}_{ij}e_{i}\otimes \overline{e}_{j})\text{,}
\]
i.e. the components $x_{ij}$ are again real.$\blacktriangle \medskip $

We now introduce some terminology to describe TQFT's for which $\eta $
extends to $\mathbf{\eta }$. First we define the notion of an algebraic
category with hermitian structure, based on the example in Section \ref
{section3}.

\begin{definition}
We say that the full algebraic category $(\mathcal{D},\mathbf{Q})$ has a
hermitian structure if:

\begin{enumerate}
\item  $K$ has an involution, $j:K\rightarrow K$,

\item  The endofunctor $Q$ on $\mathbf{D}$ acts on the underlying $K$%
-modules and their morphisms as follows: for objects $V$ of $\mathbf{D}$%
\[
G(Q(V))=G(V)^{j}\text{,}
\]
where $G(V)^{j}$ denotes the $K$-module with the same underlying set and
addition as $G(V)$ and scalar multiplication $\cdot _{j}$ given by $\alpha
\cdot _{j}x=j(\alpha )x$, where scalar multiplication in $G(V)$ is denoted
by juxtaposition. For morphisms $f:V\rightarrow W$ of $\mathbf{D}$:
\[
G(Q(f))=k_{G(W)}\circ G(f)\circ k_{G(V)}^{-1}\text{,}
\]
where $k_{G(V)}:G(V)\rightarrow G(V)^{j}$ is the $j$-semilinear map given by
$k_{G(V)}(x)=x$, for all $x\in G(V)$,

\item  for every object $V$ of $\mathcal{S}(\mathbf{D})$ the evaluation is
non-degenerate, in the sense that the associated adjoint homomorphisms
\[
G(Q(V))\rightarrow G(V)^{*}\text{ and }G(V)\rightarrow G(Q(V))^{*}\text{,}
\]
are isomorphisms, where $G(V)^{*}=\limfunc{Hom}(G(V),K)$,

\item  $\mathbf{Q}:\mathcal{D}\rightarrow \mathcal{D}$ is given, in terms of
Definition \ref{Def311}, by
\[
x_{Q(V)}=k_{G(V)}(x)\text{.}
\]
\end{enumerate}
\end{definition}

\begin{definition}
Let $(Z,\eta ):(\mathcal{C},P)\rightarrow (\mathcal{D},Q)$ be a TQFT functor
and $\mathbf{Q}$ be an extension of $Q$. If $\eta $ extends to $\mathbf{\eta
}:Z\mathbf{P}\rightarrow \mathbf{Q}Z$, we say that the pair $(Z,\mathbf{\eta
})$ is a \emph{full TQFT functor} from $(\mathcal{C},\mathbf{P})$ to $(%
\mathcal{D},\mathbf{Q})$.\newline
If in addition, $\mathcal{D}$ has a hermitian structure, $(Z,\mathbf{\eta })$
is said to be a \emph{hermitian TQFT}.\newline
A hermitian TQFT for which $K=\Bbb{C}$, $j$ is complex-conjugation and the
evaluation is positive definite for every object $\mathcal{S}(\mathbf{D})$
(i.e. the bilinear map $Q(V)\times V\rightarrow K$, corresponding to the
evaluation, is positive definite) is called a \emph{unitary} TQFT.
\end{definition}

{\noindent \bf Example: }%
The category $\mathcal{D}$ has hermitian structure both in our main example
(see page \pageref{pag33}), if we ensure that $e_{V}$ is non-degenerate, and
in the Example \ref{example37} of hermitian linear spaces.

For the $2$-dimensional TQFT's considered in the previous section, our
discussion above about the conditions for $\eta $ to extend to $\mathbf{\eta
}$ gives rise to the following characterization of unitary TQFT's:

\begin{theorem}
\label{thm56}Unitary TQFT's of the type considered in Theorem \ref{thm44}
are in 1-1 correspondence with pairs of real-valued tensors $d_{i}$, $%
p_{ijk} $ satisfying the Relations 1-4 there.$\blacktriangle $
\end{theorem}

\begin{thm_remark}
In \cite{Atiyah} Atiyah defined $(d+1)$-dimensional TQFT's over $\Bbb{C}$,
for which $X$ is a $(d+1)$-dimensional oriented differentiable manifold with
boundary $A$, and $V_{A}$ is a finite-dimensional complex vector space, with
$Z_{X}\in V_{A}$. Our Condition (\ref{eqHermit}) corresponds to Atiyah's
hermitian axiom
\[
Z_{-X}=\overline{Z}_{X}\text{,}
\]
where $-X$ denotes $X$ with opposite orientation, and $\overline{Z}_{X}$
denotes the complex conjugate of $Z_{X}$ when $A=\emptyset $, i.e. $V_{A}=%
\Bbb{C}$, and $k_{V_{A}}(Z_{X})$ otherwise.\newline
For a hermitian TQFT the evaluation on $V_{A}$ gives rise to a hermitian
form on $V_{A}$ and induces an isomorphism
\[
V_{-A}\cong V_{A}^{*}\text{,}
\]
where $-A$ denotes $A$ with opposite orientation, which is the involutory
axiom of Atiyah's paper.
\end{thm_remark}

\section{Final Comments}

The example of $2$-dimensional TQFT's was useful for illustrating the
formalism, but describes a rather straightforward topological setup. In
future work we hope to use our framework to gain new insight into more
substantial examples, in particular Stallings manifolds \cite{Stallings} in $%
3$-dimensional topology, which are obtained by self-gluing from manifolds of
the form $\Sigma \otimes I$, where $\Sigma $ is $2$-dimensional, and
parallel transport for gerbes \cite{Mackaay-Picken,Bunke-Turner-Willerton}.
An interesting avenue on the algebraic side would be to try and incorporate
infinite-dimensional vector spaces into the approach.\medskip

\bibliographystyle{unsrt}
\bibliography{art,artTQFT,Cat,LivrosAlgTop,LivrosCatHomAlg,Teses}

\begin{thebibliography}{10}

\bibitem{Witten2}
{E. Witten}.
\newblock {Topological Quantum Field Theory}.
\newblock {\em {Commun. Math. Phys.}}, 117:353--386, 1988.

\bibitem{Atiyah}
M.~Atiyah.
\newblock {Topological Quantum Field Theories}.
\newblock {\em {Publ. Math. Inst. Hautes Etudes Sci.}}, 68:175--186, 1989.

\bibitem{Segal}
{G. Segal}.
\newblock {Two-Dimensional Conformal Field Theories and Modular Functors}.
\newblock In {\em {Proceedings of the XIth International Conference on
  Mathematical Physics}}, pages 22--37, Swansea, 1988. Adam Hilger, Bristol,
  1989.

\bibitem{Witten1}
E.~Witten.
\newblock {Quantum Field Theory and the Jones Polynomial}.
\newblock {\em {Commun. Math. Phys.}}, 121:351--399, 1989.

\bibitem{Reshetikhin-Turaev}
N.~Reshetikhin and V.~G. Turaev.
\newblock {Invariants of 3-manifolds via link polynomials and quantum groups}.
\newblock {\em Invent. Math.}, 103:547--597, 1991.

\bibitem{Turaev-Viro}
V.~Turaev and O.~Viro.
\newblock {State sum invariants of three-manifolds and quantum $6j$-symbols}.
\newblock {\em Topology}, 31:865--902, 1992.

\bibitem{Crane-Yetter}
{L. Crane and D. Yetter}.
\newblock {A categorical construction of $4$D topological quantum field
  theories}.
\newblock In {\em {Quantum Topology}}, number~3 in Ser. Knots Everything, pages
  120--130. World Sci. Publishing, River Edge, NJ, 1993.

\bibitem{Dijkgraaf-Witten}
R.~Dijkgraaf and E.~Witten.
\newblock {Topological Gauge Theories and Group Cohomology}.
\newblock {\em Commun. Math. Phys.}, 129:393--429, 1990.

\bibitem{Sawin2}
S.~Sawin.
\newblock {Links, Quantum Groups and TQFT's}.
\newblock {\em {Bull. Amer. Math. Soc.}}, 33(4):413--445, 1996.

\bibitem{Birman}
{J. Birman}.
\newblock {\em {Braids, Links and Mapping Class Groups}}, volume~82 of {\em
  Annals Math. Studies}.
\newblock Princeton Univ. Press, Princeton, 1975.

\bibitem{Turaev}
{V. Turaev}.
\newblock {Operator invariants of tangles and $R$-matrices}.
\newblock {\em {Math. USSR-Izv.}}, 35(2):411--444, 1990.

\bibitem{TuraevHQFT1}
{V. Turaev}.
\newblock {Homotopy field theory in dimension 2 and group-algebras}.
\newblock 1999.
\newblock math.QA/9910010.

\bibitem{TuraevHQFT2}
V.~Turaev.
\newblock {Homotopy field theory in dimension 3 and crossed group-categories}.
\newblock 2000.
\newblock math.GT/0005291.

\bibitem{Brightwell-Turner}
{M. Brightwell and P. Turner}.
\newblock {Representations of the homotopy surface category of a simply
  connected space}.
\newblock {\em {J. Knot Theory Ramifications}}, 9(7):855--864, 2000.

\bibitem{Goncalo}
{G. Rodrigues}.
\newblock {Homotopy quantum field theories and the homotopy cobordism category
  in dimension $1+1$}.
\newblock 2001.
\newblock to appear in J. Knot Theory Ramifications, math.QA/0105018.

\bibitem{Barrett-Westbury}
{J. Barrett and B. Westbury}.
\newblock {Spherical categories}.
\newblock {\em {Adv. Math.}}, 143(2):357--375, 1999.

\bibitem{Crane-Frenkel}
L.~Crane and I.~Frenkel.
\newblock {Four-dimensional Topological Quantum Field Theory, Hopf Categories,
  and the Canonical Bases}.
\newblock {\em {J. Math. Phys.}}, 35(10):5136--5154, 1994.

\bibitem{Crane-Yetter2}
{L. Crane and D. Yetter}.
\newblock {On algebraic structures implicit in topological quantum field
  theories}.
\newblock {\em {J. Knot Theory Ramifications}}, 8(2):125--163, 1999.

\bibitem{Mackaay1}
{M. Mackaay}.
\newblock {Spherical $2$-Categories and $4$-Manifold Invariants}.
\newblock {\em {Adv. Math.}}, 143:288--348, 1999.

\bibitem{Mackaay2}
{M. Mackaay}.
\newblock {Finite Groups, Spherical $2$-Categories, and $4$-Manifold
  Invariants}.
\newblock {\em {Adv. Math.}}, 153:353--390, 2000.

\bibitem{Evans-KawBook}
{D. Evans and Y. Kawahigashi}.
\newblock {\em Quantum Symmetries on Operator Algebras}.
\newblock Oxford Mathematical Monographs. Oxford Science Publications, 1998.

\bibitem{Baez-Dolan}
{J. C. Baez and J. Dolan}.
\newblock {Higher-dimensional Algebra and Topological Quantum Field Theory}.
\newblock {\em {J. Math. Phys.}}, 36:6073--6105, 1995.

\bibitem{Baez-LangfordIV}
{J. Baez and L. Langford}.
\newblock {Higher-dimensional Algebra IV: $2$-Tangles}.
\newblock to appear in Adv. Math.

\bibitem{Bunke-Turner-Willerton}
{U. Bunke, P. Turner and S. Willerton}.
\newblock {Gerbes and homotopy quantum field theories}.
\newblock 2002.
\newblock math.AT/0201116.

\bibitem{QuinnInBook}
F.~Quinn.
\newblock {\em Lectures on Axiomatic Topological Quantum Field Theory},
  volume~I of {\em Math. Series}, pages 323--459.
\newblock Amer. Math. Soc., IAS/Park City, 1995.

\bibitem{Quinn2}
F.~Quinn.
\newblock Group categories and their field theories.
\newblock In {\em Proceedings of the Kirbyfest (Berkeley, CA, 1998)}, number~2
  in Geom. Topol. Monogr., pages 407--453 (electronic). Geom. Topol. Publ.,
  Coventry, 1999.

\bibitem{TuraevBook}
V.~Turaev.
\newblock {\em Quantum Invariants of Knots and 3-Manifolds}.
\newblock Walter de Gruyter, 1994.

\bibitem{Roger}
{R. F. Picken}.
\newblock {Reflections on Topological Quantum Field Theory}.
\newblock {\em {Rep. Math. Phys.}}, 40(2):295--303, 1997.

\bibitem{LowellAbrams}
L.~Abrams.
\newblock {Two-dimensional Topological Quantum Field Theories and Frobenius
  Algebras}.
\newblock {\em {J. Knot Theory Ramifications}}, 5(5):569--587, 1996.

\bibitem{Sawin1}
S.~Sawin.
\newblock {Direct Sum Decompositions and Indecomposable TQFT's}.
\newblock {\em {J. Math. Phys.}}, 36(12):6673--6680, 1995.

\bibitem{RuthLawrence}
{R. J. Lawrence}.
\newblock {An Introduction to Topological Field Theory}.
\newblock In {\em {Proceedings of Symposia in Applied Mathematics}}, volume~51,
  pages 89--128, 1996.

\bibitem{PickSemiao}
{R. Picken and P. Semi{\~a}o}.
\newblock {TQFT - a new direction in algebraic topology}.
\newblock In {\em {Proceedings of the Meeting ``New Developments in Algebraic
  Topology''}}, pages 85--101, {Faro, Portugal}, July 13-14 1998. {Yu. Kubyshin
  et al. eds.}, {Univ. do Algarve, 2000}.
\newblock math.QA/9912085.

\bibitem{PSemiaoPhD}
{P. Semi{\~a}o}.
\newblock {\em {Teorias Topol{\'o}gicas Qu{\^a}nticas do Campo no Contexto da Topologia
  Alg{\'e}brica}}.
\newblock {Ph.D. Dissertation}, Univ. do Algarve, \'Area Departamental de
  Matem\'atica, November 2001.

\bibitem{MacLane1}
S.~Mac Lane.
\newblock {\em Categories for the {W}orking {M}athematician}.
\newblock Springer {V}erlag, 2nd edition, 1998.

\bibitem{Kassel}
C.~Kassel.
\newblock {\em Quantum Groups}.
\newblock Springer Verlag, 1995.

\bibitem{Adamek}
{J. Ad{\'a}mek, H. Herrlich and G. Strecker}.
\newblock {\em {Abstract and Concrete Categories}}.
\newblock John {W}iley, 1990.

\bibitem{Kohno}
T.~Kohno.
\newblock {Topological invariants for 3-manifolds using representations of
  mapping class groups I}.
\newblock {\em Topology}, 31(2):203--230, 1992.

\bibitem{Stallings}
J.~Stallings.
\newblock {On Fibering Certain 3-Manifolds}.
\newblock In {\em {Topology of 3-Manifolds and related topics, Proceedings of
  the University of Georgia Institute, 1961}}, pages 95--100. {M. K. Fort, Jr.
  Ed.}, {Prentice-Hall, New Jersey, 1962}.

\bibitem{Mackaay-Picken}
{M. Mackaay and R. Picken}.
\newblock {Holonomy and parallel transport for abelian gerbes}.
\newblock {\em {Adv. Math.}}, 170(7):287--339, 2002.

\end{thebibliography}

\nocite{Adamek}

\end{document}